\documentclass[12pt]{amsart}
\usepackage[left=1in,right=1in,top=1.25in,bottom=1.25in]{geometry}
\usepackage[pdfencoding=auto]{hyperref}
\usepackage{amssymb,stmaryrd,bbm}
\usepackage{bm}
\usepackage{tikz}
\usepackage{tikz-cd}
\usepackage{pgfplots}
\pgfplotsset{compat=1.5}
\usetikzlibrary{pgfplots.fillbetween,positioning}
\usepackage{tabmac}
\usepackage{ytableau}
\usepackage{mathtools}
\title{LLT polynomials in the Schiffmann algebra}

   \author{J. Blasiak}
   \author{M. Haiman}
   \author{J. Morse}
   \author{A. Pun}
   \author{G. H. Seelinger}

    \address[Blasiak]{
    Dept.\ of Mathematics\\
    Drexel University\\
    Philadelphia, PA}
    \email{jblasiak@gmail.com}

   \address[Haiman]{Dept.\ of Mathematics\\
            University of California\\
            Berkeley, CA}
   \email{mhaiman@math.berkeley.edu}

   \address[Morse]{
   Dept.\ of Mathematics\\
            University of Virginia\\
            Charlottesville, VA}
   \email{morsej@virginia.edu}

   \address[Pun] {Dept.\ of Mathematics\\
            University of Virginia\\
            Charlottesville, VA}
   \email{ayp6e@virginia.edu}

   \address[Seelinger]{Dept.\ of Mathematics\\
   University of Michigan\\
   Ann Arbor, MI}
   \email{ghseeli@umich.edu}

   \thanks{Authors were supported by NSF Grants DMS-1855784 (J.~B.)
    and DMS-1855804 (J.~M.).}

   \date{\today}

\newtheorem{thm}{Theorem}[subsection]
\newtheorem{lemma}[thm]{Lemma}
\newtheorem{prop}[thm]{Proposition}
\newtheorem{cor}[thm]{Corollary}

\theoremstyle{definition}
\newtheorem{defn}[thm]{Definition}

\theoremstyle{remark}
\newtheorem{example}[thm]{Example}

\newtheorem{remark}[thm]{Remark}

\DeclareFontFamily{U}{mathx}{\hyphenchar\font45}
\DeclareFontShape{U}{mathx}{m}{n}{
      <5> <6> <7> <8> <9> <10>
      <10.95> <12> <14.4> <17.28> <20.74> <24.88>
      mathx10
      }{}
\DeclareSymbolFont{mathx}{U}{mathx}{m}{n}
\DeclareFontSubstitution{U}{mathx}{m}{n}
\DeclareMathAccent{\widecheck}{0}{mathx}{"71}

\newcommand{\kk}{{\mathbbm k}} 

\newcommand{\gfrak}{\mathfrak{g}}
\newcommand{\gl}{\mathfrak{gl}}

\newcommand{\ww}{{\mathbf w}}
\newcommand{\xx}{{\mathbf x}}
\newcommand{\yy}{{\mathbf y}}
\newcommand{\zz}{{\mathbf z}}

\newcommand{\sigmabold}{{\boldsymbol \sigma }}
\newcommand{\nubold}{{\boldsymbol \nu }}
\newcommand{\Abold}{{\mathbf A}}
\newcommand{\obold}{{\mathbf o}}
\newcommand{\bbold}{{\mathbf b}}

\newcommand{\Acal}{{\mathcal A}}
\newcommand{\Ecal}{{\mathcal E}}
\newcommand{\Gcal}{{\mathcal G}}
\newcommand{\Rcal}{{\mathcal R}}
\newcommand{\Scal}{{\mathcal S}}
\newcommand{\Tcal}{{\mathcal T}}

\newcommand{\NN}{{\mathbb N}}
\newcommand{\QQ}{{\mathbb Q}}
\newcommand{\ZZ}{{\mathbb Z}}

\newcommand{\ctild}{\tilde{c}}
\newcommand{\Htild}{\tilde{H}}
\newcommand{\bboldtild}{\widetilde{\mathbf{b}} }

\newcommand{\tGamma}{{\widecheck{\Gamma}}} 
\newcommand{\hGamma}{{\widehat{\Gamma}}}  
\newcommand{\ctens}{\mathbin{\widehat{\otimes}}} 

\DeclareMathOperator{\inv}{inv}
\DeclareMathOperator{\pol}{pol}
\DeclareMathOperator{\GL}{GL}

\DeclareMathOperator{\SSYT}{SSYT}
\DeclareMathOperator{\cub}{cub}

\newcommand{\barRq}{R_{+\backslash q}}
\newcommand{\barRqt}{R_{t\backslash qt}}
\newcommand{\barRt}{R_{q\backslash t}}

\newcommand{\nb}[1]{\nubold(#1)}
\newcommand{\nbm}[1]{\nubold ^{m}(#1)}
\newcommand{\diagD}[1]{\operatorname{diag}(#1)}

\newcommand{\defeq}{\mathbin{\overset{\text{{\rm def}}}{=}}}

\newcommand{\Asum}{\ensuremath{\mathcal{A}}}
\newcommand{\stre}{{\mathrm{stretch}}}
\newcommand{\north}{{\mathrm{north}}}
\newcommand{\south}{{\mathrm{south}}}
\newcommand{\east}{{\mathrm{east}}}
\newcommand{\west}{{\mathrm{west}}}

\def\TTiny{\fontsize{6.6pt}{6pt}\selectfont}

\newcommand{\mymidgray}{black!47}  

\newlength{\cellsize}
\newcommand\mytableau[1]{
\vcenter{
\let\\=\cr
\baselineskip=-16000pt \lineskiplimit=16000pt \lineskip=0pt
\halign{&\mytableaucell{##}\cr#1\crcr}}}

\newcommand{\mytableaucell}[1]{{%
\def \arg{#1}\def \void{}%
\ifx \void \arg
\vbox to \cellsize{\vfil \hrule width \cellsize height 0pt}%
\else \unitlength=\cellsize
\begin{picture}(1,1)
\put(0,0){\makebox(1,1){$#1\vphantom{-#1}$}}
\put(0,0){\line(1,0){1}}
\put(0,1){\line(1,0){1}}
\put(0,0){\line(0,1){1}}
\put(1,0){\line(0,1){1}}
\end{picture}%
\fi}}

\newcommand{\partition}[1]{{\setlength{\cellsize}{1ex}  \mytableau{#1}}}

\begin{document}


\begin{abstract}
We identify certain combinatorially defined rational functions which,
under the shuffle to Schiffmann algebra isomorphism, map to LLT
polynomials in any of the distinguished copies $\Lambda
(X^{m,n})\subset \Ecal $ of the algebra of symmetric functions
embedded in the elliptic Hall algebra $\Ecal $ of Burban and
Schiffmann.  As a corollary, we deduce an explicit raising operator
formula for the $\nabla$ operator applied to any LLT polynomial.  In
particular, we obtain a formula for $\nabla ^m s_\lambda$ which serves
as a starting point for our proof of the Loehr-Warrington conjecture
in a companion paper to this one.
\end{abstract}

\maketitle

\section{Introduction}
\setcounter{subsection}{1}
\label{s intro}

In this paper we introduce {\it Catalanimals}---symmetric rational
functions in variables $\zz = z_{1},\ldots,z_{l}$ defined by
\begin{equation}\label{e d bigraded catalan}
H(R_q,R_t,R_{qt},\lambda) = \sum_{w \in S_l} w\bigg(\frac{\zz
^\lambda \prod_{\alpha \in R_{qt}} (1-q\, t\, \zz ^\alpha )}
{\prod_{\alpha \in R_+} (1-\zz ^{-\alpha}) \prod_{\alpha \in R_q} (1-q
\, \zz ^\alpha) \prod_{\alpha \in R_t} (1-t \, \zz ^\alpha)}\bigg),
\end{equation}
depending on a weight $\lambda \in \ZZ^l$ and subsets $R_q,R_t,R_{qt}$
of the set of positive roots $R_+ = \{\epsilon_i-\epsilon_j \mid 1 \le
i < j \le l\}$ for $\GL _{l}$.  Using Catalanimals, we give explicit
formulas for elements of the elliptic Hall algebra $\Ecal$ of Burban
and Schiffmann \cite{BurbSchi12} corresponding to arbitrary LLT
polynomials $\Gcal _{\nubold }(X;q)$.  This generalizes a formula of
Negut \cite{Negut14} for elements corresponding to ribbon shaped skew
Schur functions.  The key to our more general result is the use of
Catalanimals to combinatorialize Negut's shuffle algebra tools.  As a
corollary, we also obtain a raising operator formula for $\nabla
\Gcal_\nubold $.

The Schiffmann algebra $\Ecal$ is generated by subalgebras
$\Lambda(X^{m,n})$ isomorphic to the ring of symmetric functions over
$\kk =\QQ(q,t)$, one for each coprime pair $(m,n) \in \ZZ^2$, along
with an additional central subalgebra.  The `right half-plane'
subalgebra $\Ecal ^{+}\subseteq \Ecal $ generated by $\Lambda
(X^{m,n})$ for $m>0$ is known to be isomorphic to a graded algebra
$\Scal_{\hGamma} \subseteq \bigoplus_{l}
\kk(z_{1},\ldots,z_{l})^{S_{l}}$, called the shuffle algebra, whose
degree $l$ component consists of certain symmetric rational functions
in $l$ variables.  We denote this isomorphism by $\psi_{\hGamma}
\colon \Scal_{\hGamma} \rightarrow \Ecal ^{+}$.

Our main result, Theorem~\ref{t mn LLT cuddly} (see also
Remark~\ref{rem:intro-form}), is the construction of Catalanimals
$H^{m,n}_{\nubold^m }(\zz )$ such that
\begin{equation}
\label{e:main-preview} \psi _{\hGamma}(H^{m,n}_{\nubold ^m}(\zz )) =
c^{m,n}_{\nubold ^m} \Gcal_{\nubold }[-MX^{m,n}],
\end{equation}
where $c^{m,n}_{\nubold ^m}\in \pm q^{\ZZ } t^{\ZZ }$, $M =
(1-q)(1-t)$ and the square brackets denote plethystic substitution
(see \S \ref{ss:symmetric-functions}).  Here $\Gcal_{\nubold } = \Gcal
_{\nubold}(X;q)$ is the `attacking inversions' LLT polynomial
(Definition \ref{def:G-nu}) indexed by a tuple of skew shapes $\nubold
= (\nu_{(1)}, \dots,\nu_{(k)})$; it is a $q$-analog of the product of
skew Schur functions $s_{\nu_{(1)}} \cdots s_{\nu_{(k)}}$.

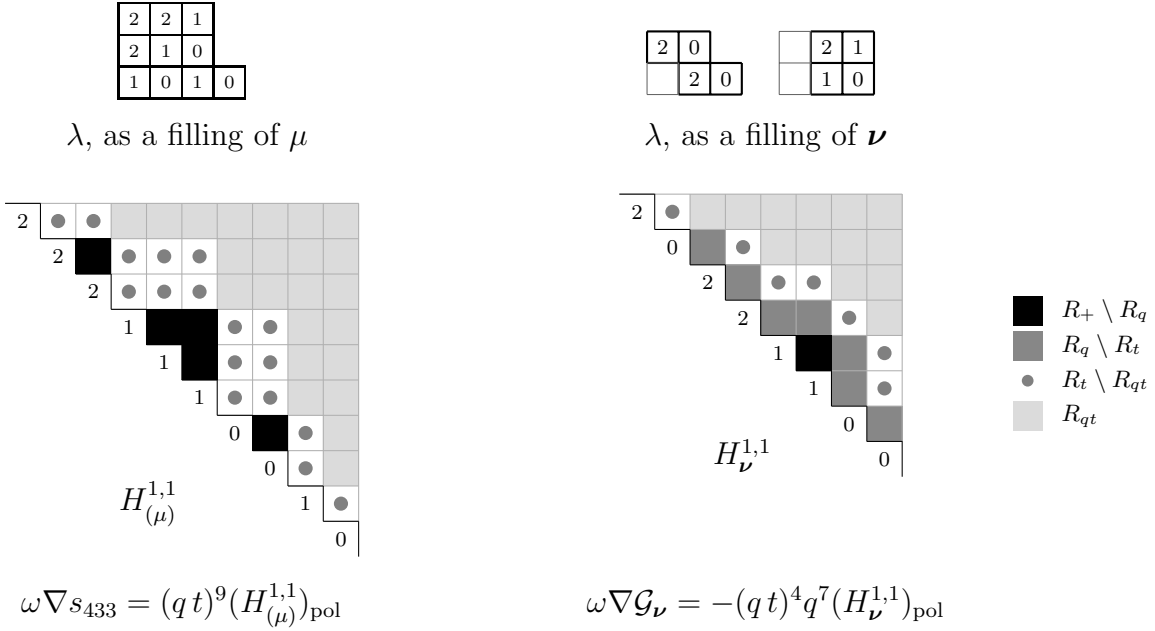
\begin{figure}
\[
\begin{array}{@{\hspace{1.2mm}}c@{\qquad\qquad\quad\quad\quad}cr}
{\fontsize{7pt}{6pt}\selectfont \tableau{2&2&1\\2&1&0\\1&0&1& 0\\} } 
&
\raisebox{-5.4mm}{
\begin{tikzpicture}
\begin{scope}[scale = .84*.5]
\draw[thick] (0,1) grid (2,2);
\draw[thick] (1,0) grid (3,1);
\draw[thin, black!60]  (0,0) grid (1,1);
\node at (0.5, 1.5) { \tiny 2};
\node at (1.5, 1.5) { \tiny 0};
\node at (1.5, 0.5) { \tiny 2};
\node at (2.5, 0.5) { \tiny 0};
\end{scope}
\begin{scope}[xshift = 50, scale = .84*.5]
\draw[thick] (1,0) grid (3,2);
\draw[thin, black!60] (0,0) grid (1,2);
\node at (1.5, 1.5) {\tiny 2};
\node at (2.5, 1.5) { \tiny 1};
\node at (1.5, 0.5) { \tiny 1};
\node at (2.5, 0.5) { \tiny 0};
\end{scope}
\end{tikzpicture} }
& \\[7mm]
\text{  $\lambda$, as a filling of $\mu$}
& 
\text{ $\lambda$, as a filling of $\nubold $}
& \\[5mm]
\begin{tikzpicture}[scale = .47]
\begin{scope}
\draw[draw = none, fill = black!14] (4,-1) rectangle (5,-2);
 \draw[draw = none, fill = black!14] (5,-1) rectangle (6,-2);
 \draw[draw = none, fill = black!14] (6,-1) rectangle (7,-2);
 \draw[draw = none, fill = black!14] (7,-1) rectangle (8,-2);
 \draw[draw = none, fill = black!14] (7,-2) rectangle (8,-3);
 \draw[draw = none, fill = black!14] (7,-3) rectangle (8,-4);
 \draw[draw = none, fill = black!14] (8,-1) rectangle (9,-2);
 \draw[draw = none, fill = black!14] (8,-2) rectangle (9,-3);
 \draw[draw = none, fill = black!14] (8,-3) rectangle (9,-4);
 \draw[draw = none, fill = black!14] (9,-1) rectangle (10,-2);
 \draw[draw = none, fill = black!14] (9,-2) rectangle (10,-3);
 \draw[draw = none, fill = black!14] (9,-3) rectangle (10,-4);
 \draw[draw = none, fill = black!14] (9,-4) rectangle (10,-5);
 \draw[draw = none, fill = black!14] (9,-5) rectangle (10,-6);
 \draw[draw = none, fill = black!14] (9,-6) rectangle (10,-7);
 \draw[draw = none, fill = black!14] (10,-1) rectangle (11,-2);
 \draw[draw = none, fill = black!14] (10,-2) rectangle (11,-3);
 \draw[draw = none, fill = black!14] (10,-3) rectangle (11,-4);
 \draw[draw = none, fill = black!14] (10,-4) rectangle (11,-5);
 \draw[draw = none, fill = black!14] (10,-5) rectangle (11,-6);
 \draw[draw = none, fill = black!14] (10,-6) rectangle (11,-7);
 \draw[draw = none, fill = black!14] (10,-7) rectangle (11,-8);
 \draw[draw = none, fill = black!14] (10,-8) rectangle (11,-9);
 \draw[draw = none, fill = gray!100] (2+0.5, -1-0.5) circle (.2);
\draw[draw = none, fill = gray!100] (3+0.5, -1-0.5) circle (.2);
\draw[draw = none, fill = gray!100] (4+0.5, -2-0.5) circle (.2);
\draw[draw = none, fill = gray!100] (5+0.5, -2-0.5) circle (.2);
\draw[draw = none, fill = gray!100] (6+0.5, -2-0.5) circle (.2);
\draw[draw = none, fill = gray!100] (4+0.5, -3-0.5) circle (.2);
\draw[draw = none, fill = gray!100] (5+0.5, -3-0.5) circle (.2);
\draw[draw = none, fill = gray!100] (6+0.5, -3-0.5) circle (.2);
\draw[draw = none, fill = gray!100] (7+0.5, -4-0.5) circle (.2);
\draw[draw = none, fill = gray!100] (8+0.5, -4-0.5) circle (.2);
\draw[draw = none, fill = gray!100] (7+0.5, -5-0.5) circle (.2);
\draw[draw = none, fill = gray!100] (8+0.5, -5-0.5) circle (.2);
\draw[draw = none, fill = gray!100] (7+0.5, -6-0.5) circle (.2);
\draw[draw = none, fill = gray!100] (8+0.5, -6-0.5) circle (.2);
\draw[draw = none, fill = gray!100] (9+0.5, -7-0.5) circle (.2);
\draw[draw = none, fill = gray!100] (9+0.5, -8-0.5) circle (.2);
\draw[draw = none, fill = gray!100] (10+0.5, -9-0.5) circle (.2);
\draw[thin, black!31] (1,-1) -- (11,-1);
\draw[thin, black!31] (2,-1) -- (2,-1);
\draw[thin, black!31] (2,-2) -- (11,-2);
\draw[thin, black!31] (3,-2) -- (3,-1);
\draw[thin, black!31] (3,-3) -- (11,-3);
\draw[thin, black!31] (4,-3) -- (4,-1);
\draw[thin, black!31] (4,-4) -- (11,-4);
\draw[thin, black!31] (5,-4) -- (5,-1);
\draw[thin, black!31] (5,-5) -- (11,-5);
\draw[thin, black!31] (6,-5) -- (6,-1);
\draw[thin, black!31] (6,-6) -- (11,-6);
\draw[thin, black!31] (7,-6) -- (7,-1);
\draw[thin, black!31] (7,-7) -- (11,-7);
\draw[thin, black!31] (8,-7) -- (8,-1);
\draw[thin, black!31] (8,-8) -- (11,-8);
\draw[thin, black!31] (9,-8) -- (9,-1);
\draw[thin, black!31] (9,-9) -- (11,-9);
\draw[thin, black!31] (10,-9) -- (10,-1);
\draw[thin, black!31] (10,-10) -- (11,-10);
\draw[thin, black!31] (11,-10) -- (11,-1);
\draw[draw = none, fill = black!100] (3,-2) rectangle (4,-3);
 \draw[draw = none, fill = black!100] (5,-4) rectangle (6,-5);
 \draw[draw = none, fill = black!100] (6,-4) rectangle (7,-5);
 \draw[draw = none, fill = black!100] (6,-5) rectangle (7,-6);
 \draw[draw = none, fill = black!100] (8,-7) rectangle (9,-8);
 \draw[thin] (1,-1) -- (2,-1);
\draw[thin] (2,-1) -- (2,-2);
\draw[thin] (2,-2) -- (3,-2);
\draw[thin] (3,-2) -- (3,-3);
\draw[thin] (3,-3) -- (4,-3);
\draw[thin] (4,-3) -- (4,-4);
\draw[thin] (4,-4) -- (5,-4);
\draw[thin] (5,-4) -- (5,-5);
\draw[thin] (5,-5) -- (6,-5);
\draw[thin] (6,-5) -- (6,-6);
\draw[thin] (6,-6) -- (7,-6);
\draw[thin] (7,-6) -- (7,-7);
\draw[thin] (7,-7) -- (8,-7);
\draw[thin] (8,-7) -- (8,-8);
\draw[thin] (8,-8) -- (9,-8);
\draw[thin] (9,-8) -- (9,-9);
\draw[thin] (9,-9) -- (10,-9);
\draw[thin] (10,-9) -- (10,-10);
\draw[thin] (10,-10) -- (11,-10);
\draw[thin] (11,-10) -- (11,-11);
\node at (3/2,-3/2) {\tiny $2 $};
\node at (5/2,-5/2) {\tiny $2 $};
\node at (7/2,-7/2) {\tiny $2 $};
\node at (9/2,-9/2) {\tiny $1 $};
\node at (11/2,-11/2) {\tiny $1 $};
\node at (13/2,-13/2) {\tiny $1 $};
\node at (15/2,-15/2) {\tiny $0 $};
\node at (17/2,-17/2) {\tiny $0 $};
\node (vv) at (19/2,-19/2) {\tiny $1 $};
\node at (21/2,-21/2) {\tiny $0 $};
\node [left=13.4mm of vv] {$H_{(\mu)}^{1,1}$};
\end{scope}
\end{tikzpicture}
&  
\raisebox{9mm}{
\begin{tikzpicture}[scale=.47]
\begin{scope}
\draw[draw = none, fill = black!14] (3,-1) rectangle (4,-2);
 \draw[draw = none, fill = black!14] (4,-1) rectangle (5,-2);
 \draw[draw = none, fill = black!14] (5,-1) rectangle (6,-2);
 \draw[draw = none, fill = black!14] (5,-2) rectangle (6,-3);
 \draw[draw = none, fill = black!14] (6,-1) rectangle (7,-2);
 \draw[draw = none, fill = black!14] (6,-2) rectangle (7,-3);
 \draw[draw = none, fill = black!14] (7,-1) rectangle (8,-2);
 \draw[draw = none, fill = black!14] (7,-2) rectangle (8,-3);
 \draw[draw = none, fill = black!14] (7,-3) rectangle (8,-4);
 \draw[draw = none, fill = black!14] (8,-1) rectangle (9,-2);
 \draw[draw = none, fill = black!14] (8,-2) rectangle (9,-3);
 \draw[draw = none, fill = black!14] (8,-3) rectangle (9,-4);
 \draw[draw = none, fill = black!14] (8,-4) rectangle (9,-5);
 \draw[draw = none, fill = gray!100] (2+0.5, -1-0.5) circle (.2);
\draw[draw = none, fill = gray!100] (4+0.5, -2-0.5) circle (.2);
\draw[draw = none, fill = gray!100] (5+0.5, -3-0.5) circle (.2);
\draw[draw = none, fill = gray!100] (6+0.5, -3-0.5) circle (.2);
\draw[draw = none, fill = gray!100] (7+0.5, -4-0.5) circle (.2);
\draw[draw = none, fill = gray!100] (8+0.5, -5-0.5) circle (.2);
\draw[draw = none, fill = gray!100] (8+0.5, -6-0.5) circle (.2);
\draw[draw = none, fill = \mymidgray] (3,-2) rectangle (4,-3);
 \draw[draw = none, fill = \mymidgray] (4,-3) rectangle (5,-4);
 \draw[draw = none, fill = \mymidgray] (5,-4) rectangle (6,-5);
 \draw[draw = none, fill = \mymidgray] (6,-4) rectangle (7,-5);
 \draw[draw = none, fill = \mymidgray] (7,-5) rectangle (8,-6);
 \draw[draw = none, fill = \mymidgray] (7,-6) rectangle (8,-7);
 \draw[draw = none, fill = \mymidgray] (8,-7) rectangle (9,-8);
 \draw[thin, black!31] (1,-1) -- (9,-1);
\draw[thin, black!31] (2,-1) -- (2,-1);
\draw[thin, black!31] (2,-2) -- (9,-2);
\draw[thin, black!31] (3,-2) -- (3,-1);
\draw[thin, black!31] (3,-3) -- (9,-3);
\draw[thin, black!31] (4,-3) -- (4,-1);
\draw[thin, black!31] (4,-4) -- (9,-4);
\draw[thin, black!31] (5,-4) -- (5,-1);
\draw[thin, black!31] (5,-5) -- (9,-5);
\draw[thin, black!31] (6,-5) -- (6,-1);
\draw[thin, black!31] (6,-6) -- (9,-6);
\draw[thin, black!31] (7,-6) -- (7,-1);
\draw[thin, black!31] (7,-7) -- (9,-7);
\draw[thin, black!31] (8,-7) -- (8,-1);
\draw[thin, black!31] (8,-8) -- (9,-8);
\draw[thin, black!31] (9,-8) -- (9,-1);
\draw[draw = none, fill = black!100] (6,-5) rectangle (7,-6);
 \draw[thin] (1,-1) -- (2,-1);
\draw[thin] (2,-1) -- (2,-2);
\draw[thin] (2,-2) -- (3,-2);
\draw[thin] (3,-2) -- (3,-3);
\draw[thin] (3,-3) -- (4,-3);
\draw[thin] (4,-3) -- (4,-4);
\draw[thin] (4,-4) -- (5,-4);
\draw[thin] (5,-4) -- (5,-5);
\draw[thin] (5,-5) -- (6,-5);
\draw[thin] (6,-5) -- (6,-6);
\draw[thin] (6,-6) -- (7,-6);
\draw[thin] (7,-6) -- (7,-7);
\draw[thin] (7,-7) -- (8,-7);
\draw[thin] (8,-7) -- (8,-8);
\draw[thin] (8,-8) -- (9,-8);
\draw[thin] (9,-8) -- (9,-9);
\node at (3/2,-3/2) {\tiny $2 $};
\node at (5/2,-5/2) {\tiny $0 $};
\node at (7/2,-7/2) {\tiny $2 $};
\node at (9/2,-9/2) {\tiny $2 $};
\node at (11/2,-11/2) {\tiny $1 $};
\node at (13/2,-13/2) {\tiny $1 $};
\node at (15/2,-15/2) {\tiny $0 $};
\node (vv) at (17/2,-17/2) {\tiny $0 $};
\node[left=11.4mm of vv] {\raisebox{2mm}{$H_{\nubold }^{1,1}$}};
\end{scope}
\end{tikzpicture} }
&
\ \ \
\raisebox{16mm}{
\begin{tikzpicture}[scale = .39]
\begin{scope}
\draw[draw = none, fill = black!100] (2,-1) rectangle (3,-2);
\node[anchor = west] at (3,-1.5) {\scriptsize  \  $R_+ \setminus R_q$};
\end{scope}
\begin{scope}[yshift = -34*1]
\draw[draw = none, fill = \mymidgray] (2,-1) rectangle (3,-2);
\node[anchor = west] at (3,-1.5) {\scriptsize  \  $R_q \setminus R_t$};
\end{scope}
\begin{scope}[yshift = -34*2]
\draw[thin, black!0] (2,-1) -- (3,-1);
\draw[thin, black!0] (2,-1) -- (2,-2);
\draw[thin, black!0] (2,-2) -- (3,-2);
\draw[thin, black!0] (3,-2) -- (3,-1);
\draw[draw = none, fill = gray!100] (2+0.5, -1-0.5) circle (.2);
\node[anchor = west] at (3,-1.5) {\scriptsize  \  $R_t \setminus R_{qt}$};
\end{scope}
\begin{scope}[yshift = -34*3]
\draw[draw = none, fill = black!14] (2,-1) rectangle (3,-2);
\node[anchor = west] at (3,-1.5) {\scriptsize \  $R_{qt}$};
\end{scope}
\end{tikzpicture}}
\\[2mm]
\omega\nabla s_{433}= (q\, t)^9 (H_{(\mu)}^{1,1})_{\pol} & \
\omega\nabla \Gcal_{\nubold }= -(q\, t)^4 q^7 (H_{\nubold
}^{1,1})_{\pol } &
\end{array}
\]
\caption{\label{fig:Catalanimals intro}%
(i) The Catalanimal $H_{(\mu)}^{1,1}$ for $\mu = (433)$.  (ii) The
Catalanimal $H_{\nubold}^{1,1}$ for $\nubold = ((32)/(10),(33)/(11))$.
These are illustrated by drawing the root sets in an $\ell\times \ell$
grid labeled by matrix-style coordinates, with the sets $R_+ \setminus
R_q$, $R_q \setminus R_t$, $R_t \setminus R_{qt}$, $R_{qt}$ specified
according to the legend on the right; the weight $\lambda$ is written
on the diagonal with $\lambda_1$ in the upper left.}
\end{figure}

The precise definition of $H_{\nubold ^m}^{m,n}$ is given in
\S\S\ref{ss:def LLT Catalanimals} and \ref{ss:mn LLT Catalanimal}, but
we give a flavor of our results here. For $m = n = 1$, the Catalanimal
$H_{\nubold }^{1,1}$ has root sets $R_+ \supseteq R_q \supseteq R_t
\supseteq R_{qt}$, determined as follows using the same attacking
inversion combinatorics as in the definition of the LLT polynomial
$\Gcal_{\nubold }$:
\begin{itemize}
\item $R_+ \setminus R_q  \xleftrightarrow{\ } $ pairs of boxes in the
same diagonal,
\item $R_q \setminus R_t  \, \,  \xleftrightarrow{\ } $ the attacking pairs,
\item $R_t \setminus R_{qt}  \,  \xleftrightarrow{\ } $ pairs going
between adjacent diagonals,
\end{itemize}
where the boxes of $\nubold $ are numbered $1,\dots, l$ in reading
order (see Example \ref{ex:LLT Catalanimals}).  The weight $\lambda$
is obtained by filling each diagonal $D$ of $\nubold $ with the value
\[
1+ \chi(\text{$D$ contains a row start}) - \chi( \text{$D$ contains a
row end}),
\]
where $\chi(P) = 1$ if $P$ is true
or $0$ if $P$ is false, and then reading this filling in the reading
order---see Figure \ref{fig:Catalanimals intro}.

Shuffle algebra representatives for elements of $\Lambda(X^{m,1})
\subseteq \Ecal^+$ carry information about the symmetric function
operator $\nabla$ from the theory of Macdonald polynomials.
Specifically, if a symmetric function $f(X)$ is related to a Catalanimal
$H \in \Scal_{\hGamma}$ by $\psi _{\hGamma}(H) = f[-MX^{m,1}]$, then
$\nabla^m f$ and $H$ are related by
\begin{align}
\label{eq:intro nabla}
(\omega \nabla^m f) (z_1,\dots, z_l) =  H_{\pol},
\end{align}
where $H_{\pol}$ is obtained by expanding $H$ as an infinite series of
$\GL_l$ characters and truncating to polynomial $\GL_l$ characters
(see \S\ref{ss:nabla}).

Combining \eqref{e:main-preview} and \eqref{eq:intro nabla}, we obtain
an explicit raising operator formula for $\nabla^m$ on any LLT
polynomial
\begin{equation} \label{e:Gnu-formula-prototype}
(\omega \nabla^m \Gcal_{\nubold })
(z_1,\dots, z_l) = (c_{\nubold^m }^{m,1})^{-1} (H_{\nubold }^{m,1})_{\pol}.
\end{equation}
In a companion paper
\cite{loehrwarringtonconj}, we 
use the case of formula~\eqref{e:Gnu-formula-prototype} where
the LLT polynomial is a Schur function to prove the Loehr-Warrington
conjecture \cite{LoehWarr08}, a
combinatorial formula for $\nabla^m s_\lambda$ in terms of LLT polynomials.
By the Schur positivity of LLT polynomials \cite{GrojHaim07},
this implies that $\nabla^m s_\lambda$ is Schur positive 
up to a global sign.

After giving background on symmetric functions and LLT polynomials in
\S \ref{S LLT}, we develop the Schiffmann and shuffle algebra tools
needed to prove identity \eqref{e:main-preview} in \S\S
\ref{s:schiffmann}--\ref{s:principal-specialization}.  Then, in \S\S
\ref{s:LLT-catalanimals}--\ref{s mn LLT Catalanimal}, we give the full
description of the Catalanimals $H^{m,n}_{\nubold ^m}$ and apply these
tools to complete the proof.

\section{Background on symmetric functions and LLT polynomials}
\label{S LLT}

\subsection{Symmetric functions and partition diagrams}
\label{ss:symmetric-functions}

Let $\Lambda = \Lambda _{\kk }(X)$ be the algebra of symmetric
functions in an infinite alphabet of variables $X =
x_{1},x_{2},\ldots$, with coefficients in the field $\kk = \QQ (q,t)$.
We follow the notation of Macdonald~\cite{Macdonald95} for the graded
bases of $\Lambda $.  We write $\omega \colon \Lambda \rightarrow
\Lambda $ for the involutory $\kk $-algebra automorphism determined by
$ \omega s_{\lambda } = s_{\lambda^*}$ for Schur functions $s_{\lambda
}$, where $\lambda^*$ denotes the conjugate partition of $\lambda$.

Given $f\in \Lambda $ and an expression $A$ involving indeterminates,
such as a polynomial, rational function, or formal series, the
plethystic evaluation $f[A]$ is defined by writing $f$ as a polynomial
in the power-sums $p_{k}$ and evaluating with $p_{k}\mapsto p_{k}[A]$,
where $p_{k}[A]$ is the result of substituting $a^{k}$ for every
indeterminate $a$ occurring in $A$.  The variables $q, t$ from our
ground field $\kk $ count as indeterminates.  We will often use the
fact that $f(x_{1},x_{2},\ldots) = f[X]$ with $X = x_{1}+x_{2}+\cdots
$.

The algebra $\Lambda $ is a Hopf algebra with coproduct given by
\begin{equation}\label{e:Lambda-coprod}
\Delta f = f[X+Y] \in \Lambda \otimes \Lambda = \Lambda (X)\otimes
\Lambda (Y).
\end{equation}
Here $f\in \Lambda $ and we use separate alphabets $X,Y$ to
distinguish the tensor factors.

We fix notation for the series
\begin{equation}\label{e:Omega}
\Omega = 1 + \sum _{k>0} h_{k} = \exp \sum _{k>0} \frac{p_{k}}{k},
\quad \text{or}\quad \Omega [a_{1}+a_{2}+\cdots -b_{1}-b_{2}-\cdots ]
= \frac{\prod_{i} (1-b_{i})}{\prod_{i} (1-a_{i})}
\end{equation}
and the quantities
\begin{equation}\label{e:M-hat}
M = (1-q)(1-t) \qquad \text{and} \qquad \widehat{M} =
\big(1-\frac{1}{q\, t} \big)M.
\end{equation}
An example of plethystic substitution using these quantities, which
will arise again later, is
\begin{equation}\label{eq:Omega-Mz}
\Omega[-w/y\, \widehat{M}] = \frac{(1 - q\, t\, w/y)(1 - q^{-1}w/y)(1
- t^{-1}w/y)}{(1 - (q\, t)^{-1}w/y)(1 - q\, w/y)(1 - t\, w/y)}\, .
\end{equation}

The (French style) diagram of a partition $\lambda$ is the set of
lattice points $\{(i,j)\mid 1\leq j\leq \ell(\lambda ),\; 1\leq i \leq
\lambda _{j}\}$, where $\ell(\lambda )$ is the length of $\lambda $.
We often identify $\lambda$ and its diagram with the set of lattice
squares, or {\em boxes}, with northeast corner at a point $(i,j)\in
\lambda$.  For $\mu\subseteq\lambda$, the {\it skew shape} $\lambda /
\mu$ is the set of boxes of $\lambda$ not contained in $\mu$.

For a box $a = (i,j) \in \ZZ^2$ (usually in some given skew shape), we
let $\south(a) = (i,j-1) \in \ZZ^2$ denote the box immediately south
of $a$, and define $\north(a) = (i,j+1), \west(a)= (i-1,j),$ and
$\east(a)= (i+1,j)$ similarly.

\subsection{LLT polynomials}
\label{ss LLT}

We recall the attacking inversions description of LLT polynomials
from~\cite{HaHaLoReUl05}.

Let $\nubold = (\nu_{(1)},\ldots,\nu_{(k)})$ be a tuple of skew
shapes.  We consider the set of boxes in $\nubold $ to be the disjoint
union of the sets of boxes in the $\nu _{(i)}$.  The {\em content} of
a box $a=(i,j)$ in row $j$, column $i$ of a skew diagram is
$c(a)=i-j$.  Fix $\epsilon >0$ small enough that $k\, \epsilon <1$.
The {\em adjusted content} of a box $a\in \nu_{(i)}$ is $\ctild (a) =
c(a)+i\, \epsilon $.  A {\em diagonal} of $\nubold $ is the set of
boxes of a fixed adjusted content, or, in other words, the set of
boxes of fixed content in one of the shapes $\nu _{(i)}$.  We write
$\diagD{a}$ for the diagonal containing a box $a$ of $\nubold$.

The {\em reading order} on $\nubold $ is the total ordering $<$ on the
boxes of $\nubold $ such that $a<b \Rightarrow \ctild (a)\leq \ctild
(b)$ and boxes on each diagonal increase from southwest to northeast
(see Example \ref{ex:LLT Catalanimals}).  We say that an ordered pair
of boxes $(a,b)$ in $\nubold $ is an {\em attacking pair} if $a<b$ in
reading order and $0<\ctild (b)-\ctild(a)<1$.

A {\em semistandard tableau} on the tuple $\nubold  $ is a map
$T\colon \nubold  \rightarrow \ZZ _{+}$ which restricts to a
semistandard Young tableau on each component $\nu_{(i)}$.  We write
$\SSYT (\nubold  )$ for the set of these.  An {\em attacking
inversion} in $T$ is an attacking pair $(a,b)$ such that $T(a)>T(b)$.
Let $\inv (T)$ denote the number of attacking inversions in $T$.

\begin{defn}\label{def:G-nu}
The {\em LLT polynomial} indexed by a tuple of skew diagrams $\nubold 
$ is the generating function, which is known to be symmetric
\cite{HaHaLoReUl05,LaLeTh97},
\begin{equation}\label{e:G-nu}
\Gcal_{\nubold  }(X;q) = \sum _{T\in \SSYT (\nubold 
)}q^{\inv (T)}\xx ^{T},
\end{equation}
where $\xx ^{T} = \prod _{a\in \nubold  } x_{T(a)}$.
\end{defn}

It follows directly from the definition that
\begin{equation}
\Gcal_{\nubold 
}(x_1,x_2,\dots,y_1,y_2,\dots; q) = \sum q^{A(\nubold '',\nubold ')}
\Gcal_{\nubold ' }(X; q) \Gcal_{\nubold '' }(Y; q),
\end{equation}
where the sum is over all partitions of $\nubold $ into a lower order
ideal $\nubold '$ and upper order ideal $\nubold ''$,
and
$A(\nubold '',\nubold ')$ is the number of attacking pairs $(a,b)$
with $a \in \nubold ''$, $b \in \nubold '$.
In other words, LLT polynomials satisfy the coproduct formula
\begin{equation}\label{e:LLT coprod}
\Gcal_{\nubold  }[X+Y] = \sum q^{A(\nubold '',\nubold ')}
\Gcal_{\nubold ' }[X] \Gcal_{\nubold '' }[Y].
\end{equation}
When we write LLT polynomials in plethystic notation, our convention is to
suppress the $q$.

We will also need a more general combinatorial formalism for LLT
polynomials involving a `signed' alphabet $\Acal = \Acal_{+} \coprod
\Acal _{-}$ with a {\em positive} letter $v\in \Acal _{+}$ and a {\em
negative} letter $\overline{v}\in \Acal _{-}$ for each $v\in \ZZ_{+}$,
with total ordering $1 < 2 < \cdots < \overline{1} < \overline{2} <
\cdots$.

A {\em super tableau} on a tuple of skew shapes $\nubold  $ is a map
$T\colon \nubold  \rightarrow \Acal $, weakly increasing along rows
and columns, with positive letters increasing strictly on columns and
negative letters increasing strictly on rows.

An attacking inversion in a super tableau is an attacking pair
$(a,b)$ such that either
$T(a)>T(b)$ in the ordering on $\Acal $, or $T(a)=T(b) = \overline{v}$
with $\overline{v}$ negative.  As before, $\inv (T)$ denotes the
number of attacking inversions.

\begin{lemma}[{\cite[(81--82) and Proposition 4.2]{HaHaLo05}}]
\label{lem:super-G}
We have the identity
\begin{equation}\label{e:super-G}
\omega _{Y} \Gcal _{\nubold  }[X+Y] = \sum _{T}
q^{\inv (T)} \xx ^{T_{+}} \yy ^{T_{-}},
\end{equation}
where the sum is over all super tableaux $T$ on $\nubold $, and
\begin{equation}\label{e:super-tableau-weight}
\xx ^{T_{+}} \yy ^{T_{-}} = \prod _{a\in \nubold }
\begin{cases}
x_{i},& T(a) = i \in \Acal _{+},\\
y_{i},& T(a) = \overline{i} \in \Acal _{-}.
\end{cases}
\end{equation}
\end{lemma}

\section{Catalanimals in the shuffle and Schiffmann algebras}
\label{s:schiffmann}

We briefly introduce and fix notation for the shuffle algebra of
Feigin et al. \cite{FeHaHoShYa09}, Feigin and Tsymbauliak
\cite{FeigTsym11}, and Negut \cite{Negut14}, and the elliptic Hall
algebra of Burban and Schiffmann \cite{BurbSchi12}, which we refer to
as the {\em Schiffmann algebra}.

We then give a preliminary description of how Catalanimals connect
with the shuffle and Schiffmann algebras (\S\ref{ss four ways to
view}) and how they relate to the $\nabla$ operator
(\S\ref{ss:nabla}).

\subsection{The shuffle algebra}
\label{ss:shuffle-algebra}

Let $\Gamma = \Gamma (w,y)$ be a non-zero rational function over
$\kk$.  The {\em large concrete shuffle algebra} is the graded
associative algebra with underlying space
\begin{equation}\label{e:big-shuffle}
\Rcal _{\Gamma } =\bigoplus _{l} \Rcal _{\Gamma }^{l} = \bigoplus _{l} \kk
(z_{1},\ldots,z_{l})^{S_{l}},
\end{equation}
equipped with the `shuffle' product whose graded component $\Rcal
_{\Gamma }^{k}\times \Rcal _{\Gamma }^{l-k}\rightarrow \Rcal _{\Gamma
}^{l}$ is defined by
\begin{equation}\label{e:shuffle-product}
f\cdot g = \sum _{w\in S_{l}/(S_{k}\times S_{l-k})} w
\bigl(f(z_{1},\ldots,z_{k})g(z_{k+1},\ldots,z_{l})
\prod_{i=1}^{k}\prod _{j=k+1}^{l}\Gamma (z_{i},z_{j})\bigr).
\end{equation}

Define symmetrization operators $\sigma _{\Gamma }^{l} \colon \kk
(z_{1},\ldots,z_{l}) \to \Rcal _\Gamma^l$ by
\begin{equation}\label{e:sigma-gamma}
\sigma _{\Gamma }^{l}(f) = \sum _{w\in S_{l}}
w\bigl(f(z_{1},\ldots,z_{l})\prod _{i<j}\Gamma (z_{i},z_{j})\bigr).
\end{equation}
The $\sigma _{\Gamma }^{l}$ are the components of a surjective graded
algebra homomorphism $\sigma _{\Gamma }\colon U\rightarrow \Rcal
_{\Gamma }$, where $U$ is the algebra with underlying space
\begin{equation}\label{e:concatenation-algebra}
U = \bigoplus_{l} U^l = \bigoplus_{l} \kk(z_{1},\ldots,z_{l})
\end{equation}
and the `concatenation' product whose component $U^{k}\times
U^{l-k}\rightarrow U^{l}$ is defined by
\begin{equation}\label{e:concatenation-product}
f\cdot g = f(z_{1},\ldots,z_{k})g(z_{k+1},\ldots,z_{l})\, .
\end{equation}
Let $I_{\Gamma }^{l} = \ker (\sigma _{\Gamma }^{l})$, so we have
\begin{equation}\label{e:shuffle-ideal}
\ker (\sigma _{\Gamma }) = I_{\Gamma } = \bigoplus _{l} I_{\Gamma
}^{l}\subseteq U
\end{equation}
and an induced isomorphism $\sigma _{\Gamma }\colon U/I_{\Gamma
}\xrightarrow {\simeq }\Rcal _{\Gamma }$.  We call $U/I_{\Gamma }$ the
{\em large abstract shuffle algebra}.

The subalgebra of $U$ consisting of Laurent polynomials,
\begin{equation}\label{e:tensor-algebra}
T = \bigoplus _{l} T^{l} = \bigoplus _{l} \kk[z_{1}^{\pm
1},\ldots,z_l^{\pm 1}] \subseteq U,
\end{equation} 
is generated by the basis elements $z_{1}^{a}$ of $T^{1} = \kk
[z_{1}^{\pm 1}]$ and is isomorphic to the tensor algebra on these
generators.

The {\em abstract shuffle algebra}, or just {\em shuffle algebra} for
short, is the image
\begin{equation}\label{e:small-shuffle}
S_{\Gamma } = T/(I_{\Gamma }\cap T) = (T+I_{\Gamma })/I_{\Gamma
}\subseteq U/I_{\Gamma }
\end{equation}
of $T$ in the large abstract shuffle algebra $U/I_{\Gamma }$.  The
isomorphism $U/I_{\Gamma }\cong \Rcal _{\Gamma }$ induced by $\sigma
_{\Gamma }$ restricts to an isomorphism
\begin{equation}\label{e:small-concrete-shuffle}
\sigma _{\Gamma }\colon S_{\Gamma } \xrightarrow{\simeq } \Scal
_{\Gamma } \defeq \sigma _{\Gamma }(T)\subseteq \Rcal _{\Gamma }.
\end{equation}
from the shuffle algebra $S_{\Gamma }$ to the subalgebra $\Scal
_{\Gamma }$ of $\Rcal _{\Gamma }$ generated by the elements
$z_{1}^{a}\in \Rcal _{\Gamma }^{1}$.  We call $\Scal _{\Gamma }$ the
{\em concrete shuffle algebra.}

The diagram below summarizes the relationships between these algebras.
\begin{equation}\label{e:shuffle-CD}
\begin{tikzcd}[/tikz/column 3/.append style={anchor=base west}]
T \ar[r,hookrightarrow]\ar[d,twoheadrightarrow]
  & U \ar[r,equal]\ar[d,twoheadrightarrow]
  &[-4ex]  \bigoplus_l \kk(z_1,\ldots,z_l)\\[-2ex]
S_{\Gamma } \ar[r,hookrightarrow]
            \ar[d,"\rotatebox{90}{$\simeq$}"',"\sigma _{\Gamma }"]
  & U/I_{\Gamma } \ar[d,"\rotatebox{90}{$\simeq$}"',"\sigma _{\Gamma }"]\\
\Scal _{\Gamma } \ar[r,hookrightarrow] & \Rcal _{\Gamma } \ar[r,equal] &
  \bigoplus_l \kk(z_1,\ldots,z_l)^{S_l}
\end{tikzcd}
\end{equation}

\begin{remark}\label{rem:Gamma-choices}
If $\Gamma'(w,y) = h(w,y)\Gamma (w,y)$, where $h(w,y) = h(y,w)\not =0$,
then $I_{\Gamma '} = I_{\Gamma }$ and $S_{\Gamma '} = S_{\Gamma }$.
In this way, different choices of $\Gamma $ give rise to different
concrete realizations $\Scal _{\Gamma }$ of the same shuffle algebra
$S_{\Gamma }$.  The algebra $S_{\Gamma }$ is thus more canonical than
$\Scal _\Gamma$ and has a simpler product (concatenation of Laurent
polynomials), but has non-trivial relations.  The algebra $\Scal
_{\Gamma }$ has the advantage that the symmetric rational function
$f(z_{1},\ldots,z_{l})$
representing an element of $\Scal _{\Gamma
}^{l}$ is unique.
\end{remark}

\begin{remark}\label{rem:non-Laurent-phi}
We can think of the abstract shuffle algebra $S_{\Gamma }$ either as a
quotient of $T$ or as the subalgebra $(T+I_{\Gamma })/I_{\Gamma }$ of
$U/I_{\Gamma }$.  From the latter point of view, any rational function
$\phi \in T^{l}+I_{\Gamma }^{l}$, that is, $\phi $ congruent modulo
$I_{\Gamma }^{l}$ to a Laurent polynomial, represents an element of
$S_{\Gamma }$.  An explicit Laurent polynomial $\eta $ such that $\eta
\equiv \phi \pmod{I_{\Gamma }^{l}}$ may be hard to compute, is not
unique, and need not have any simple form, even if $\phi $ does.
However, we can often avoid the need to construct $\eta $, since its
image $\sigma _{\Gamma }(\eta ) = \sigma _{\Gamma }(\phi )$ in the
concrete shuffle algebra $\Scal _{\Gamma }$ can be computed directly
from $\phi $.
\end{remark}

\subsection{The shuffle to Schiffmann algebra isomorphism}
\label{ss:shuffle-iso}

We use the same notation as in \cite{delta, paths} for the Schiffmann
algebra $\Ecal$ of \cite{BurbSchi12}.  In our notation, $\Ecal $ is
generated by subalgebras $\Lambda (X^{m,n})$ isomorphic to the algebra
$\Lambda$ of symmetric functions over $\kk = \QQ(q,t)$, one for each
pair of coprime integers $m,n$, and a central Laurent polynomial
subalgebra $\kk [c_{1}^{\pm 1}, c_{2}^{\pm 1}]$, subject to some
defining relations.  
A translation between this notation and that of
\cite{BurbSchi12,Schiffmann12,SchiVass13} can be found in
\cite[\S3.2]{paths}, and a presentation of the defining relations in
\cite[\S3]{delta}.

The `right half-plane' subalgebra $\Ecal ^{+}\subseteq \Ecal $ is
generated by $\Lambda (X^{m,n})$ for $m>0$, or equivalently (as a
consequence of the relations) by the elements
$p_{1}(X^{1,a})$. Schiffmann and Vasserot \cite[Theorem
10.1]{SchiVass13} showed that $\Ecal ^{+}$ is isomorphic to the
shuffle algebra $S_{\Gamma }$ for a suitable choice of $\Gamma $.
We use the following version of their theorem, modified the same way
as in \cite[Proposition 3.5.1]{paths}.

\begin{thm}[{\cite{SchiVass13}}]\label{thm:shuffle-isomorphism}
Let $\Gamma (w,y)$ be a rational function such that
\begin{equation}\label{e:Gamma-over-Gamma}
\Gamma(w,y)/\Gamma(y,w) = \Omega[-w/y\, \widehat{M}],
\end{equation}
where $\Omega [-w/y\, \widehat{M}]$ is given by
\eqref{eq:Omega-Mz}. Then
there is an algebra isomorphism $\psi \colon S_{\Gamma }\rightarrow
\Ecal ^{+}$ given on the generators
by $\psi (z_1^{a}) = p_{1}[-MX^{1,a}]$.
\end{thm}

If $\Gamma (w,y)$ and $\Gamma '(w,y)$ both satisfy
\eqref{e:Gamma-over-Gamma}, then they differ by a symmetric factor
$h(w,y) = \Gamma '(w,y)/\Gamma (w,y)$, as in
Remark~\ref{rem:Gamma-choices}.  Accordingly, from this point on, for
any $\Gamma (w,y)$ satisfying \eqref{e:Gamma-over-Gamma}, we fix
\begin{equation}\label{e:S-without-Gamma}
S = S_\Gamma,\quad S^l = S^l_\Gamma,\quad I = I_\Gamma,\quad I^l =
I^l_\Gamma\, .
\end{equation}
Although the abstract shuffle algebra $S$ does not depend on the
choice of $\Gamma$ in Theorem \ref{thm:shuffle-isomorphism}, the
concrete shuffle algebra $\Scal _{\Gamma }$ does.  The following two
choices of $\Gamma$ which satisfy \eqref{e:Gamma-over-Gamma}
turn out to be convenient.
\begin{align}
\label{e:Gamma hat}
\hGamma(w,y) &= \frac{1 - q\, t\, w/y}{(1-y/w)(1-q\, w/y)(1-t\, w/y)},  \\
\label{e:Gamma tilde}
\tGamma(w,y) &= (1 - w/y)(1-q\, y/w)(1-t\, y/w)(1-q\, t\, w/y).
\end{align}
For either of these we define $\psi _{\Gamma }$ to be the isomorphism
\begin{equation}\label{e:psi-gamma}
\psi _{\Gamma } = \psi \circ \sigma _{\Gamma }^{-1}\colon \Scal
_{\Gamma } \xrightarrow{\simeq } \Ecal ^{+},
\end{equation}
so we have a commutative diagram, with all arrows isomorphisms,
\begin{equation}\label{e:4-ways-lite}
\begin{tikzcd}
\Scal _{\hGamma } \ar[dr,"\psi _{\hGamma }"]\\
S \ar[r,"\psi "]\ar[u,"\sigma _{\hGamma }"']\ar[d,"\sigma _{\tGamma }"]
 & \Ecal ^{+}\, .\\
\Scal _{\tGamma } \ar[ur,"\psi _{\tGamma }"']
\end{tikzcd}
\end{equation}

\subsection{Grading}
\label{ss:grading}

The algebra $\Ecal$ has a $\ZZ^2$ grading in which $\kk [c_{1}^{\pm
1}, c_{2}^{\pm 1}]$ has degree $(0,0)$ and $f(X^{m,n})$ has degree
$(dm,dn)$ for $f \in \Lambda$ of degree $d$.  We denote by
$\Ecal^{(a,b)}$ the $(a,b)$-graded component.  The subalgebra
$\Ecal^+$ is an $\NN \times \ZZ$ graded subalgebra of $\Ecal$.  Set
$(\Ecal^+)^{(l,\bullet)} = \bigoplus_{d \in \ZZ} (\Ecal^+)^{(l,d)}$.

The algebra $T$ is $\NN \times \ZZ$ graded, where the component of
degree $(l,d)$ consists of Laurent polynomials $\phi \in T^l$
homogeneous of degree $d$.

If $\Gamma(w,y)$ is a function of $w/y$ satisfying
\eqref{e:Gamma-over-Gamma}, and in particular for $\Gamma = \hGamma $
or $\Gamma = \tGamma $, the symmetrization operators $\sigma _{\Gamma
}^{l}$ are degree preserving,
so $T \cap I=T \cap I_{\hGamma}
= T \cap I_\tGamma$ is an $\NN \times \ZZ$ graded ideal.  The shuffle
algebra $S$ therefore inherits an $\NN \times \ZZ$ grading from $T$,
and $\sigma _{\hGamma }$, $\sigma _{\tGamma }$ induce $\NN \times \ZZ$
gradings on $\Scal _{\hGamma }$ and $\Scal _{\tGamma }$ such that a
symmetric rational function $h(z_{1},\ldots,z_{l})$ in $\Scal
_{\hGamma }^{l}$ or $\Scal _{\tGamma }^{l}$ belongs to $\Scal
_{\hGamma }^{(l,d)}$ or $\Scal _{\tGamma }^{(l,d)}$ if and only if it
is homogeneous of degree $d$ in the variables $z_{i}$.

The following is clear from the definitions.
\begin{prop}\label{pr:grading}
The isomorphisms $\psi\colon S \rightarrow \Ecal^+$, $\psi_\hGamma
\colon \Scal _{\hGamma }\rightarrow \Ecal^+$, and $\psi_\tGamma \colon
\Scal _{\tGamma } \rightarrow \Ecal^+ $ preserve the $\NN \times \ZZ $
grading.
\end{prop}

\subsection{Catalanimals and their cubs}
\label{ss four ways to view}

Here and throughout, $R = \{\alpha _{ij} \mid 1\leq i,j\leq l,\; i\not
=j \}$ denotes the set of roots for $\GL _{l}$, where $\alpha _{ij} =
\epsilon_i - \epsilon_j \in \ZZ^l$, and $R_{+} = \{\alpha _{ij}\in R
\mid i<j\}$ the set of positive roots.  The number $l$ will usually be
understood from the context; otherwise we specify it by writing
$R(\GL_{l})$ or $R_{+}(\GL_{l})$.

\begin{defn}\label{def:Catalanimal}
Given a weight $\lambda \in \ZZ ^{l}$ and subsets $R_q, R_t, R_{qt}
\subseteq R_{+}$, we define the corresponding {\em Catalanimal} of
length $l$ to be the symmetric rational function
\begin{equation}\label{e:Catalanimal-def}
H(R_q,R_t,R_{qt},\lambda) \defeq \sum_{w \in S_l} w\bigg(\frac{\zz
^\lambda \prod_{\alpha \in R_{qt}} (1 - q\, t\, \zz ^\alpha )}
{\prod_{\alpha \in R_+} (1 - \zz ^{-\alpha}) \prod_{\alpha \in R_q} (1
- q\, \zz ^\alpha) \prod_{\alpha \in R_t} (1 - t\, \zz
^\alpha)}\bigg),
\end{equation}
in variables $\zz = z_1,\dots, z_l$, where $\zz ^\lambda$
stands for $z_{1}^{\lambda_{1}}\cdots z_{l}^{\lambda_{l}}$.
\end{defn}

We also define two related functions
\begin{gather}\label{eq:phi for Catalanimal}
\phi(R_q,R_t,R_{qt},\lambda) \, \defeq \, \frac{\zz ^\lambda
\prod_{\alpha \in R_+ \setminus R_q}(1-q\, \zz ^\alpha)\prod_{\alpha
\in R_+ \setminus R_t} (1-t\, \zz ^\alpha)}{\prod_{\alpha \in R_+
\setminus R_{qt}}(1-q\, t\, \zz ^\alpha)}, \\
\label{eq:laurent form g}
g(R_q,R_t,R_{qt},\lambda) \defeq \sum_{w \in S_l} w\bigl(\zz
^\lambda \! \prod_{\alpha \in R_+} \!\ (1 - \zz ^{\alpha}) \!\!\!
\prod_{\alpha \in R \setminus R_q } \!\!\!  (1 - q\, \zz ^\alpha)
\!\!\! \prod_{\alpha \in R \setminus R_t } \!\!\!  (1 - t\, \zz
^\alpha) \!  \prod_{\alpha \in R_{qt}} \!  (1 - q\, t\, \zz ^\alpha)
\bigr),
\end{gather}
so that
\begin{align}\label{eq:phi-and-g-to-H}
\sigma _{\hGamma }\big(\phi (R_{q}, R_{t}, R_{qt}, \lambda )\big)
 & = H(R_{q}, R_{t}, R_{qt}, \lambda ),\\
\sigma _{\tGamma }\big(\phi (R_{q}, R_{t}, R_{qt}, \lambda )\big)
 & = g(R_{q}, R_{t}, R_{qt}, \lambda ).
\end{align}

Note that the following conditions on a Catalanimal $H =
H(R_{q},R_{t},R_{qt},\lambda )$ of length $l$ and its associated
functions $\phi = \phi (R_{q},R_{t},R_{qt},\lambda )$ and $g =
g(R_{q},R_{t},R_{qt},\lambda )$ are equivalent:
\begin{enumerate}
\item [(i)] $H$ belongs to the concrete shuffle algebra $\Scal _{\hGamma }$;
\item [(ii)] $g$ belongs to the concrete shuffle algebra $\Scal _{\tGamma }$;
\item [(iii)] $\phi \in T^{l}+I^{l}$, hence $\phi $ represents an
element of the shuffle algebra $S$, as in Remark~\ref{rem:non-Laurent-phi}.
\end{enumerate}
When these conditions hold, there is a corresponding element of the
Schiffmann algebra
\begin{equation}\label{e:zeta-psi-H-g}
\zeta = \psi (\phi ) = \psi _{\hGamma }(H) = \psi _{\tGamma }(g)\in
\Ecal ^{+}\, .
\end{equation}
Our work here focuses on identifying certain Catalanimals that satisfy
the above conditions and have the further property that $\zeta \in
\Lambda (X^{m,n})$ for some $(m,n)$ with $m>0$.

\begin{defn}\label{def:cub}
Let $H= H(R_q,R_t, R_{qt}, \lambda)$ be a Catalanimal.  If $\psi
_{\hGamma }(H) \in \Lambda (X^{m,n})$,
we call the symmetric function
$f(X)$ such that $\psi _{\hGamma }(H) = f[-MX^{m,n}]$ its {\em cub},
and write $\cub(H) = f$.
\end{defn}

To summarize the discussion above in the setting of Definition
\ref{def:cub}, given a Catalanimal $H= H(R_q,R_t, R_{qt}, \lambda)$ of
length $l$ such that $\psi _{\hGamma }(H) = f[-MX^{m,n}]$, the four
objects $H$, $\phi (R_{q},R_{t},R_{qt},\lambda )$,
$g(R_{q},R_{t},R_{qt},\lambda )$, and $f[-MX^{m,n}]$ are related as
shown below.  We will frequently go back and forth between these
viewpoints.
\begin{equation*}
\begin{tikzcd}[/tikz/column 3/.append style={anchor=base west}]
&[1ex] H \in \Scal _{\hGamma }^l 
\ar[dr,"\psi _{\hGamma }", end anchor={[xshift=-2em]}]
\\[2ex]
\phi \in T^l+I^l
\ar[r,twoheadrightarrow]
\ar[ru,"\sigma_{\hGamma}^l", end anchor={[xshift=1.5em]}]
\ar[rd,swap,"\sigma_{\tGamma}^l",end anchor={[xshift=1.5em]}]
& S^l = (T^l+I^l)/I^l
\ar[r,"\psi "]
\ar[u,"\sigma _{\hGamma }^l",xshift=2ex]
\ar[d,"\sigma _{\tGamma}^l",swap,xshift=2ex]
& (\Ecal ^{+})^{(l,\bullet)}\ni f[-MX^{m,n}]
\\[2ex]
& g \in \Scal _{\tGamma }^l
\ar[ur,swap,"\psi _{\tGamma }",end anchor={[xshift=-2em]}]
\end{tikzcd}
\end{equation*}

\subsection{Catalanimals and the operator \texorpdfstring{$\nabla$}{nabla}}
\label{ss:nabla}

As touched on in the introduction, one nice consequence of having a
Catalanimal representative for an element $f[-MX^{m,1}]$ in the
Schiffmann algebra is a raising operator formula for $\nabla
^{m}f(X)$, where $\nabla$ is the linear operator introduced in
\cite{BeGaHaTe99}, which acts diagonally on the basis of modified
Macdonald polynomials $\Htild_{\mu }(X;q,t)$~\cite{GarsHaim96} by
$\nabla \Htild _{\mu } = t^{n(\mu )}q^{n(\mu ^{*})} \Htild _{\mu }$,
with $n(\mu ) = \sum _{i} (i-1)\mu _{i}.$

The following lemma relates the operator $\nabla$ to the action of $\Ecal $ on
$\Lambda $ constructed by Schiffmann
and Vasserot \cite{SchiVass13}.  Here we use the 
version of this action given by
\cite[Proposition 3.3.1]{paths}.  

\begin{lemma}\label{lem:nabla-m-f}
For any symmetric function $f$, the element $f[-M X^{m,1}]\in \Ecal $
acting on $1\in \Lambda (X)$ is given by
\begin{equation}\label{e:nabla-m-f}
f[-M X^{m,1}]\cdot 1 = \nabla^{m} f(X).
\end{equation}
\end{lemma}

\begin{proof}
By \cite[Lemma 3.4.1]{paths},
$f(X^{m,1})$ acts as $\nabla ^{m} f[-X/M]^{\bullet }\, \nabla ^{-m}$,
where $f[-X/M]^{\bullet }$ is the operator of multiplication by
$f[-X/M]$.  Since $\nabla (1) = 1$, the result follows.
\end{proof}

We denote the Weyl symmetrization operator for $\GL _{l}$ by
\begin{equation}\label{e:Weyl-symmetrization}
\sigmabold (f(z_{1},\ldots,z_{l})) = \sum _{w\in S_{l}} w \left(
\frac{f(\zz )}{\prod _{\alpha \in R_{+}}(1 - \zz ^{-\alpha })} \right).
\end{equation}
If  $\lambda \in \ZZ ^{l}$ is a
dominant weight, then
$\sigmabold (\zz ^{\lambda }) = \chi _{\lambda }$ is the corresponding
irreducible $\GL _{l}$ character. For any weight $\mu $, we
either have $\sigmabold (\zz ^{\mu }) = \pm \chi _{\lambda }$ for a suitable
$\lambda $, or $\sigmabold (\zz ^{\mu }) = 0$.

Let $\eta (z_{1},\ldots,z_{l})\in \kk [z_{1}^{\pm },\ldots,z_{l}^{\pm
1}]$ be a Laurent polynomial.  A {\em $q,t$ raising operator series}
is an expression of the form
\begin{equation}\label{e:raising-op-series}
h(\zz ) = \sigmabold \left(\frac{\eta (\zz )}{\prod_{\alpha \in R_{+}}
(1 - q\, \zz ^\alpha) \prod_{\alpha \in R_{+}} (1 - t\, \zz ^\alpha)}
\right),
\end{equation}
interpreted as an infinite formal linear combination of irreducible
$\GL _{l}$ characters $\chi _{\lambda }$ with coefficients in $\kk $
by expanding the factors $(1- q\, \zz^\alpha)^{-1} = 1+q\,
\zz^\alpha+\cdots $ and $(1- t\, \zz^\alpha)^{-1} = 1+t\,
\zz^\alpha+\cdots $ as geometric series before applying $\sigmabold $.
This makes sense since for each $\lambda $, the set of weights $\mu $
such that $\sigmabold (\zz ^{\mu }) = \pm \chi _{\lambda }$ is finite.

By writing $\eta (\zz )$ in the form $c(q,t)\eta '(\zz )$, where
$\eta' \in \ZZ [q,t][z_{1}^{\pm 1 },\ldots,z_{l}^{\pm 1}]$, we see that
$h(\zz )$ can be expressed as a scalar in $\kk $ times a power series
in $q$, $t$ over the ring of virtual $\GL _{l}$ characters $\ZZ
[z_{1}^{\pm 1},\ldots,z_{l}^{\pm 1}]^{S_{l}}$, and also that this
power series expansion is the same as the raising operator series
expansion.  It follows that $h(\zz )$ considered as a symmetric
rational function over $\kk $ determines its raising operator series
expansion (not every symmetric rational function is a raising operator
series, however).

In most of this paper, we regard Catalanimals merely as symmetric
rational functions.  However, rewriting \eqref{e:Catalanimal-def} as
\begin{equation}\label{e:Catalanimal-series}
H(R_{q},R_{t},R_{qt},\lambda ) = \sigmabold \left( \frac{\zz ^\lambda
\prod_{\alpha \in R_{qt}} (1 - q\, t\, \zz ^\alpha )}{\prod_{\alpha
\in R_q} (1 - q\, \zz ^\alpha) \prod_{\alpha \in R_t} (1 - t\, \zz
^\alpha)} \right) ,
\end{equation}
we see that every Catalanimal can also be viewed as a $q,t$ raising
operator series.

The {\em polynomial characters} of $\GL _{l}$ are the irreducible
characters $\chi _{\lambda }$ for which $\lambda \in \NN ^{l}$, that
is, $\lambda $ is a partition.  We define the {\em polynomial part}
$h(\zz )_{\pol }$ of a raising operator series $h(\zz )$ to be its
truncation to polynomial characters.  If $h(\zz )$ is homogeneous of
degree $d$, then all irreducible characters $\chi _{\lambda }$ in it
have 
$\lambda _{1} +\cdots + \lambda _{l} = d$, so $h(\zz )_{\pol }$
is a finite linear combination of polynomial $\GL _{l}$
characters---that is, a symmetric polynomial in $l$ variables over
$\kk $.

\begin{prop}\label{prop:nabla-on-cuddlyish}
Let $H= H(R_{q},R_{t},R_{qt},\lambda )$, be a Catalanimal of length $l$ 
such that $\psi _{\hGamma
}(H)\in \Lambda (X^{m,1})$, and let $f(X)$ be its cub.  Then
\begin{equation}\label{ecor:nabla on cub}
(\omega \nabla^m f) (z_1,\dots, z_l) = H_{\pol}.
\end{equation}
Moreover, this determines $\nabla^m f$, since the Schur expansion of
the symmetric function $\omega \nabla^m f(X)$ contains only terms
$s_{\lambda }$ with $\ell(\lambda) \leq l$.
\end{prop}

\begin{proof}
By Lemma~\ref{lem:nabla-m-f}, we can replace $\nabla^m f$ with
$f[-M X^{m,1}]\cdot 1$.  The result now follows from
\cite[Proposition 3.5.2]{paths} by taking $\phi $ in \cite[(48)]{paths} to be
a Laurent polynomial congruent modulo $I^{l}$ to
$\phi (R_{q},R_{t},R_{qt},\lambda )$ and noting that the right hand side of
\cite[(48)]{paths} then becomes $H_{\pol }$.
\end{proof}

\subsection{Shuffle algebra toolkit}
\label{ss:toolkit}

Negut \cite{Negut14} provides a useful toolkit for working with the
shuffle algebra, which we will use extensively in \S \S \ref{s
cuddly}--\ref{s:principal-specialization}, below.  
For the convenience of readers who wish to compare our versions
of results cited from \cite{Negut14} with the originals, we briefly
discuss how Negut's notation and conventions are related to ours.

For $m>0$, the elements of $\Ecal $ denoted $u_{km,kn}$ in
\cite{BurbSchi12} and \cite{Negut14} are $\omega p_{k}(X^{n,m})$ in
our notation.  Because we have switched the indices $m, n$, the
positive half subalgebra denoted $\Ecal ^{+}$ in \cite{Negut14} is
actually the upper half-plane subalgebra generated by $\Lambda
(X^{n,m})$ for $m>0$ in our notation.  However, there is an
automorphism of $\Ecal $ which carries $f(X^{n,m})$ to $f(X^{m,-n})$
for $m>0$, and we use this to identify our right half-plane subalgebra
$\Ecal ^{+}$ with the subalgebra generated by the $u_{km,kn}$ for
$m>0$ in \cite{BurbSchi12}, so that $u_{km,kn}$ corresponds to $\omega
p_{k}(X^{m,-n})$.

The shuffle algebra $\Acal ^{+}$ in \cite{Negut14} is related to our
shuffle algebra as follows.  The parameters $q_{1},q_{2}$ in
\cite{Negut14} are $t,q\in \kk $, and the algebra $\Acal ^{+}$ with
$z_{i}$ replaced by $z_{i}^{-1}$ coincides with our concrete shuffle
algebra $\Scal _{\Gamma }$ for
\begin{equation}\label{e:Negut-Gamma}
\Gamma (w,y) = \frac{(1 - w/y)(1 - q\, t\, w/y)}{(1 - q\, w/y)(1 - t\,
w/y)}.
\end{equation}
This function $\Gamma (w,y)$ is $\omega (y/w)$ in the notation of
\cite[(2.3)]{Negut14}.

The isomorphism $\Upsilon $ in Negut \cite[Theorem~3.1]{Negut14} sends
$u_{1,-a}$ to $z_{1}^{-a}$.  Hence, $\Upsilon ^{-1}$ corresponds in our
notation to an isomorphism $S\cong \Scal _{\Gamma }\xrightarrow
{\simeq } \Ecal ^{+}$ sending $z_{1}^{a}$ to $p_{1}(X^{1,a})\in \Ecal
^{+}$.  Since we defined $\psi \colon S\rightarrow \Ecal ^{+}$ in
Theorem~\ref{thm:shuffle-isomorphism} by $\psi (z_{1}^{a}) = p_{1}[-M
X^{1,a}] = -M p_{1}(X^{1,a})$, we see that $\psi $ differs by a factor
$(-M)^{l}$ on $S^{l}$ from the isomorphism corresponding to $\Upsilon
^{-1}$.
By \cite[Theorem~1.1]{Negut14}, the elements $P_{k,d}\in \Acal ^{+}$
defined by \cite[(1.2)]{Negut14} are given by $P_{k,d} = \Upsilon
(u_{k,d})$.  Using the identification of $u_{km,-kn}$ with $\omega
p_{k}(X^{m,n})$, we have the following diagram of isomorphisms and
corresponding elements.
\newcommand{\vin}{\rotatebox[origin=c]{90}{$\in$}}
\begin{equation}\label{e:Negut-diagram}
\begin{tikzcd}
\Acal ^{+}\ar[r,"\simeq ","z_{i}\mapsto z_{i}^{-1}"']
 &[2ex] S \ar[r,"\simeq ","\psi "'] &[-4ex] \Ecal ^{+}\\[-4ex]
\vin & & \vin\\[-5ex]
P_{km,-kn} \ar[rr,mapsto] & & (-M)^{km}\omega p_{k}(X^{m,n})
\end{tikzcd}
\end{equation}

\section{Cuddly Catalanimals}
\label{s cuddly}

In this section we identify combinatorial conditions which guarantee
that a Catalanimal belongs to $\Scal _{\hGamma }$ and that its image
under $\psi _{\hGamma }$ belongs to one of the distinguished
subalgebras $\Lambda (X^{m,n})$ of the Schiffmann algebra.

\subsection{Tame Catalanimals}
\label{ss:tame-cats}

Negut \cite[Theorem~2.2]{Negut14} gives a criterion based on the wheel
condition of Feigin et al.\ \cite{FeHaHoShYa09} for a symmetric
rational function to belong to the concrete shuffle algebra.  For
$\Scal _{\tGamma }$, Negut's criterion takes the form given by the
theorem below.  Note that, because $\tGamma (w,y)$ in \eqref{e:Gamma
tilde} is a Laurent polynomial, the elements of $\Scal _{\tGamma }$
are symmetric Laurent polynomials and not just rational functions.

A symmetric Laurent polynomial $g(\zz ) \in \kk[z_{1}^{\pm
1},\ldots,z_{l}^{\pm 1}]^{S_l}$ satisfies the {\em wheel condition} if
it vanishes whenever any three of the variables $z_i,z_j,z_k$ are in
the ratio $(z_i:z_j:z_k) = (1:q:q\, t)$ or $(1:t:q\, t)$.  If $l<3$,
the wheel condition holds vacuously.

\begin{thm}[\cite{Negut14}]\label{thm:wheel}
A symmetric Laurent polynomial $g(z_{1},\ldots,z_{l})$ belongs to
$\Scal _{\tGamma }^{l}$ if and only if it satisfies the wheel
condition and vanishes whenever $z_{i} = z_{j}$ for $i\not =j$.
\end{thm}

To connect this with \cite[Theorem~2.2]{Negut14} we remark that, up to
an irrelevant factor, our $g(\zz )$ is the numerator in
\cite[(2.4)]{Negut14} with the variables inverted.  Because $g(\zz )$
is symmetric, the condition that it vanishes whenever $z_{i} = z_{j}$
is equivalent to $\prod _{i<j}(z_{i}-z_{j})^{2}$ dividing $g(\zz )$.
The wheel condition then applies to the factor that remains after
dividing by $\prod _{i<j}(z_{i}-z_{j})^{2}$.

If $A,B\subseteq R_{+}$ are subsets of the positive roots for $\GL
_{l}$, we set $[A,B] = R_{+}\cap (A+B)$.  The reason for this notation
is that if $\gfrak _{C} = \sum _{\alpha \in C}\gfrak _{\alpha }$
denotes a sum of root spaces, then $\gfrak _{[A,B]} = [\gfrak
_{A},\gfrak_{B}]$ in the Lie algebra $\gl _{l}$.

\begin{defn}\label{def:tame-catalanimal}
A Catalanimal $H(R_{q},R_{t},R_{qt},\lambda )$ is {\em tame} if the
root sets $R_{q}$, $R_{t}$, $R_{qt}$ satisfy
\begin{equation}\label{e:tame-condition}
[R_{q},R_{t}]\subseteq R_{qt}.
\end{equation}
\end{defn}

Strictly speaking, tameness is a condition on the root sets rather
than the rational function $H(R_{q},R_{t},R_{qt},\lambda )$, 
but it has the following consequence for the function.

\begin{prop}\label{prop:tame-catalanimal}
If $H = H(R_{q},R_{t},R_{qt},\lambda )$ is a tame Catalanimal, then
$H\in \Scal _{\hGamma }$; hence it corresponds to an element $\psi
_{\hGamma }(H)\in \Ecal ^{+}$ in the Schiffmann algebra.
\end{prop}

\begin{proof}
It's equivalent to show that $g(\zz ) = g(R_{q},R_{t},R_{qt},\lambda
)$ belongs to $\Scal _{\tGamma }$.  Using Theorem~\ref{thm:wheel}, we
need to check that $g(\zz )$ vanishes when any two distinct variables
$z_{i}$, $z_{j}$ are set equal, or when any three are in the ratio
$(1:q:q\, t)$ or $(1 : t :q\, t)$.  We will verify that the $w = 1$
term in \eqref{eq:laurent form g} vanishes under any of these
conditions; this suffices since the conditions are symmetric in the
variables $z_i$.

The factor $\prod _{\alpha \in R_{+}}(1 - \zz ^{\alpha })$ gives the
required vanishing when $z_{i} = z_{j}$.

Suppose $(z_{i}:z_{j}:z_{k}) = (1:q:q\, t)$.  If $\alpha _{ij}\not \in
R_{q}$, the factor $(1 - q\, \zz^{\alpha _{ij}})$ in the product over
$R\setminus R_q$ in \eqref{eq:laurent form g} vanishes, while if
$\alpha _{jk}\not \in R_{t}$, the factor $(1 - t\, \zz ^{\alpha
_{jk}})$ in the product over $R \setminus R_t$ vanishes.  If neither
of these factors appears, we have $\alpha_{ij}\in R_{q}$ and $\alpha
_{jk}\in R_{t}$, hence $\alpha _{ik}\in R_{qt}$ by hypothesis.  Thus,
the product over $R_{qt}$ contains the factor $(1 - q\, t\,
\zz^{\alpha _{ik}})$ which vanishes.  The same reasoning with the
roles of $q$ and $t$ exchanged applies if $(z_{i}:z_{j}:z_{k}) =
(1:t:q\, t)$.
\end{proof}

\subsection{Cuddly Catalanimals}
\label{ss:cuddlies}

\begin{defn}\label{def:mn-cuddly}
We use the abbreviations $I^{c} = [l]\setminus I$ for $I\subseteq [l]
= \{1,\ldots,l \}$ and $A^{I,J} = \{\alpha _{ij} \in A \mid i\in I,\,
j\in J \}$ for $A\subseteq R(GL_{l})$ and $I,J\subseteq [l]$.  We also
write $\sum A = \sum _{\alpha \in A} \alpha $ for the sum of the roots
in $A$.  For any weight $\nu = (\nu _{1},\ldots,\nu _{l})$, we define
$|\nu | = \nu _{1}+\cdots +\nu _{l}$, and let $\nu _{I} = (\nu
_{i_{1}},\ldots,\nu _{i_{k}})$ denote the subsequence with index set
$I = \{i_{1}<\cdots <i_{k} \}\subseteq [l]$.

Let $(m,n) \in \ZZ_+ \times \ZZ$ be a pair of coprime integers.  A
Catalanimal $H(R_q,R_t, R_{qt},\lambda)$ of length $l$ is
\emph{(m,n)-cuddly} if
\begin{enumerate}
\item [(a)] it is tame, that is, $[R_q,R_t] \subseteq R_{qt}$;
\item [(b)] $|\lambda| = l n/m$ (in particular, $m$ must divide $l$);
and
\item [(c)] it satisfies the {\em cuddliness bounds}
\end{enumerate}
\begin{equation}\label{e:cuddly-bound}
\big|\, \lambda[I]_I\, \big| \le \left| I \right| \frac{n}{m} \quad
\text{for all $I \subseteq \{1,\dots, l\}$},
\end{equation}
where
\begin{equation}\label{eqdef:S bracket I}
\lambda[I] = \lambda + \sum R_+^{I,I^c} -\sum R_{q}^{I,I^c} - \sum
R_{t}^{I,I^c}+ \sum R_{qt}^{I,I^c}.
\end{equation}
\end{defn}

\begin{example}
The Catalanimal below is $(1,1)$-cuddly.  It is drawn with the same
conventions as Figure \ref{fig:Catalanimals intro}.
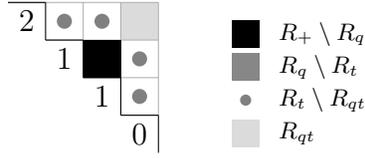
\[
\begin{tikzpicture}[scale = .5]
\draw[draw = none, fill = black!14] (4,-1) rectangle (5,-2);
\draw[draw = none, fill = gray!100] (2+0.5, -1-0.5) circle (.2);
\draw[draw = none, fill = gray!100] (3+0.5, -1-0.5) circle (.2);
\draw[draw = none, fill = gray!100] (4+0.5, -2-0.5) circle (.2);
\draw[draw = none, fill = gray!100] (4+0.5, -3-0.5) circle (.2);
\draw[thin, black!31] (1,-1) -- (5,-1);
\draw[thin, black!31] (2,-1) -- (2,-1);
\draw[thin, black!31] (2,-2) -- (5,-2);
\draw[thin, black!31] (3,-2) -- (3,-1);
\draw[thin, black!31] (3,-3) -- (5,-3);
\draw[thin, black!31] (4,-3) -- (4,-1);
\draw[thin, black!31] (4,-4) -- (5,-4);
\draw[thin, black!31] (5,-4) -- (5,-1);
\draw[draw = none, fill = black!100] (3,-2) rectangle (4,-3);
\draw[thin] (1,-1) -- (2,-1);
\draw[thin] (2,-1) -- (2,-2);
\draw[thin] (2,-2) -- (3,-2);
\draw[thin] (3,-2) -- (3,-3);
\draw[thin] (3,-3) -- (4,-3);
\draw[thin] (4,-3) -- (4,-4);
\draw[thin] (4,-4) -- (5,-4);
\draw[thin] (5,-4) -- (5,-5);
\node at (3/2,-3/2) {$2$};
\node at (5/2,-5/2) {$1$};
\node at (7/2,-7/2) {$1$};
\node at (9/2,-9/2) {$0$};
\end{tikzpicture}
\qquad
\raisebox{0mm}{
\begin{tikzpicture}[scale = .37]
\begin{scope}
\draw[draw = none, fill = black!100] (2,-1) rectangle (3,-2);
\node[anchor = west] at (3,-1.5) {\scriptsize  \  $R_+ \setminus R_q$};
\end{scope}
\begin{scope}[yshift = -34*1]
\draw[draw = none, fill = \mymidgray] (2,-1) rectangle (3,-2);
\node[anchor = west] at (3,-1.5) {\scriptsize  \  $R_q \setminus R_t$};
\end{scope}
\begin{scope}[yshift = -34*2]
\draw[thin, black!0] (2,-1) -- (3,-1);
\draw[thin, black!0] (2,-1) -- (2,-2);
\draw[thin, black!0] (2,-2) -- (3,-2);
\draw[thin, black!0] (3,-2) -- (3,-1);
\draw[draw = none, fill = gray!100] (2+0.5, -1-0.5) circle (.2);
\node[anchor = west] at (3,-1.5) {\scriptsize  \  $R_t \setminus R_{qt}$};
\end{scope}
\begin{scope}[yshift = -34*3]
\draw[draw = none, fill = black!14] (2,-1) rectangle (3,-2);
\node[anchor = west] at (3,-1.5) {\scriptsize \  $R_{qt}$};
\end{scope}
\end{tikzpicture}}
\]
Condition (a) holds since $[R_q,R_t] = R_{qt}$, and $|\lambda|
= l = 4$ verifies (b).  For (c), we must check $|\lambda[I]_I| \le
|I|$ for all subsets $I \subset \{1,2,3,4\}$; this computation is
illustrated below for the subsets of size 2, with bold letters
indicating the subsequence $\lambda[I]_I$ of $\lambda[I]$.
\[
\begin{array}{ccccccccc}
I              &  12   & 13   & 14   & 23   & 24   & 34 \\[1.4mm]
\lambda[I]     &  {\textbf{11}}11 & {\textbf{1}}2{\textbf{0}}1 &
 {\textbf{0}}22{\textbf{0}}
& 2{\textbf{00}}2 & 2{\textbf{2}}0{\textbf{0}} & 21{\textbf{10}} \\[1.4mm] 
|\lambda[I]_I|  &   {2}   &  {1}   &  {0}   &   {0}  &   {2}  & {1} & 
  { \le |I|} 
\end{array}
\]
\end{example}

Below we will show that if $H$ is an $(m,n)$-cuddly Catalanimal then
$\psi_{\hGamma} (H)\in \Lambda (X^{m,n})$.  For this we need the following
criterion for $\psi _{\tGamma }$ to map an element of $\Scal _{\tGamma
}$ into the subalgebra of $\Ecal ^{+}$ generated by the $\Lambda
(X^{m,n})$ for $n/m\leq p$.

\begin{thm}[{\cite{Negut14}}]\label{thm:g-in-Ecal-p}
Let $\Ecal ^{+}_{\leq p}$ be the subalgebra of $\Ecal ^{+}$ generated
by the $\Lambda (X^{m,n})$ for $n/m\leq p$.  Given $g =
g(z_{1},\ldots,z_{l})\in \Scal _{\tGamma }^{l}$, the image $\psi
_{\tGamma }(g)\in (\Ecal ^{+})^{(l,\bullet )}$ belongs to $\Ecal ^{+}
_{\leq p}$ if and only if $g$ is supported on monomials in the $S_{l}$
orbits of dominant weights $\mu$ satisfying $\mu_1 + \cdots + \mu_k
\le 2k(l-k) + kp$ for all $k = 1,\ldots,l$.
\end{thm}

We briefly explain how this follows from \cite{Negut14}.  As noted in
\S \ref{ss:toolkit}, there is an isomorphism $\Acal ^{+}\cong S$
between Negut's shuffle algebra and ours that inverts the variables
$z_{i}$.  There is also an anti-isomorphism that does not invert the
variables, under which the subalgebra generated by the elements
$P_{k,d}\in \Acal ^{+}$ for $d/k\leq p$ corresponds to $\psi
^{-1}(\Ecal ^{+}_{\leq p})\subseteq S$.  The proof of
\cite[Theorem~1.1]{Negut14} shows that this subalgebra is the same as
the one denoted $\Acal ^{p}\subseteq \Acal ^{+}$ in
\cite[Proposition~2.3]{Negut14}.  The criterion on the support of $g$
in Theorem~\ref{thm:g-in-Ecal-p} is a reformulation in terms of $\Scal
_{\tGamma }$ of the condition \cite[(2.7)]{Negut14} that defines
$\Acal ^{p}$.

\begin{prop}\label{prop:cuddly}
Let $H = H(R_q,R_t, R_{qt},\lambda)$ be a Catalanimal of length $l$
and let $g(\zz ) = g(R_q,R_t, R_{qt},\lambda)$ be the corresponding
symmetric Laurent polynomial defined in \eqref{eq:laurent form g}.
\begin{enumerate}
\item[(i)] If $H$ satisfies the cuddliness bounds
\eqref{e:cuddly-bound} for $(m,n)$, then $g(\zz )$ satisfies the
condition in Theorem~\ref{thm:g-in-Ecal-p} for $p =n/m$.
\item[(ii)]
If $H$ is $(m,n)$-cuddly, then $H\in
\Scal _{\hGamma }$ with $\psi _{\hGamma }(H) \in \Lambda(X^{m,n})$, so
$H$ has a cub (Definition \ref{def:cub}).
\end{enumerate}
\end{prop}

\begin{proof}
For (i), we must show that for every monomial $\zz ^{\nu }$ occurring
in $g(\zz )$ and every $I \subseteq [l] $ of size $|I| = k$, we have
$|\nu_I| \leq 2k(l-k)+kn/m$.  We will show that in fact this holds for
each term in \eqref{eq:laurent form g}.  By symmetry, it is enough to
check it for the $w = 1$ term,
\begin{equation}\label{e:w=1-term-of-g}
\zz ^\lambda \! \prod_{\alpha \in R_+} \!\ (1 - \zz ^{\alpha}) \!\!\!
\prod_{\alpha \in R \setminus R_q } \!\!\!  (1 - q\, \zz ^\alpha)
\!\!\! \prod_{\alpha \in R \setminus R_t } \!\!\!  (1 - t\, \zz
^\alpha) \!  \prod_{\alpha \in R_{qt}} \!  (1 - q\, t\, \zz ^\alpha).
\end{equation}
To get a term $\zz ^\nu$ maximizing $|\nu_I|$ in this product, we must
choose the $\zz ^\alpha$ term from the factors with $\alpha\in R^{I,
I^c}$ and the constant term from the factors with $\alpha\in
R^{I^{c},I}$.  For $\alpha \in R^{I,I} \cup R^{I^c,I^c}$, we can
choose either term from the factor in question, since $|\alpha_I|=0$.
Since $|R^{I, I^c}| = k(l-k)$ for any $I$ of size $k$, $|(R \setminus
R_q)^{I,I^c}| = k(l-k) - |R_q^{I,I^c}|$ and $|(R \setminus
R_t)^{I,I^c}| = k(l-k) - |R_t^{I,I^c}|$.  Hence, $\nu$ chosen as
indicated satisfies
\begin{equation}
|\nu_I| = |\lambda[I]_I| +
2k(l-k)
\end{equation}
and the cuddliness bound \eqref{e:cuddly-bound} implies $|\nu_I| \le
2k(l-k)+ k n/m$.

For (ii), we know from Proposition~\ref{prop:tame-catalanimal} that
$H\in \Scal _{\hGamma }$ and from (i) that $\psi _{\hGamma }(H) \in
\Ecal ^{+}_{\leq n/m}$.  The definition of $\Ecal ^{+}_{\leq n/m}$ and
the $\NN \times \ZZ$ grading of $\Ecal ^{+}$ imply that $(\Ecal
^{+}_{\leq n/m})^{(l,ln/m)}\subseteq \Lambda (X^{m,n})$.  Since $H$ is
homogeneous of degree $ |\lambda | = ln/m$, it follows from
Proposition~\ref{pr:grading} that $\psi _{\hGamma }(H) \in
\Lambda(X^{m,n})$.
\end{proof}

\section{The coproduct}
\label{s:coproduct}

The full Schiffmann algebra $\Ecal $ is constructed in
\cite{BurbSchi12} as the Drinfeld double of a subalgebra of the Hall
algebra of coherent sheaves on an elliptic curve.  This subalgebra
corresponds to the algebra denoted $\Acal ^{\geq }$ in Negut
\cite[Proposition 4.1]{Negut14} and, under the identifications in \S
\ref{ss:toolkit}, to the subalgebra $\Ecal ^{\geq }$ of $\Ecal $
generated by $\Ecal ^{+}$ and $\Lambda (X^{0,-1})$ in our notation.
The relations in $\Ecal ^{\geq }$ yield a tensor product decomposition
$\Ecal ^{\geq } = \Lambda (X^{0,-1}) \otimes \Ecal ^{+}$ (as a vector
space).

By \cite[Proposition~4.5]{BurbSchi12}, the Hall algebra realization
gives rise to a geometrically defined coproduct $\Delta $ on $\Ecal
^{\geq }$ taking values in a suitably completed tensor product $\Ecal
^{\geq } \ctens \Ecal ^{\geq }$; the corresponding coproduct on the
shuffle algebra is described in \cite{Negut14}.  Here we will use
properties of $\Delta $ to obtain a combinatorial coproduct formula
for the cub of a cuddly Catalanimal.

\subsection{Leading term}
\label{ss:coprod-leading-term}

When evaluated on $\Lambda (X^{m,n})$, the coproduct $\Delta $ on
$\Ecal ^{\geq }$ has a leading term that coincides with the standard
coproduct on $\Lambda (X^{m,n})$.  To make this precise we define
$\Ecal ^{+}_{<p}$, $\Ecal ^{+}_{>p}$ to be the subalgebras of $\Ecal
^{+}$ generated by the $\Lambda (X^{m,n})$ for $n/m<p$ or $n/m>p$,
respectively (similar to the definition of $\Ecal ^{+}_{\leq p}$ in
Theorem~\ref{thm:g-in-Ecal-p}).  We also define $\Ecal ^{\geq }_{<p}$
to be the subalgebra of $\Ecal ^{\geq }$ generated by $\Ecal
^{+}_{<p}$ and $\Lambda (X^{0,-1})$; it decomposes as $\Ecal ^{\geq
}_{<p} = \Lambda (X^{0,-1}) \otimes \Ecal ^{+}_{<p}$.

The following proposition is a consequence of either \cite[p.~1212,
line 6]{BurbSchi12} or \cite[Lemma 5.3]{Negut14}, with a geometric
proof in \cite{BurbSchi12} and a shuffle algebra proof in
\cite{Negut14}.

\begin{prop}[{\cite{BurbSchi12, Negut14}}]\label{prop:leading-term}
For any $f\in \Lambda $ and coprime integers $m,n$ with $m>0$, the
coproduct in $\Ecal ^{\geq }$ evaluated on $f(X^{m,n})$ has the form
\begin{equation}\label{e:leading-term}
\Delta (f(X^{m,n})) = f[X^{m,n}_{(1)}+X^{m,n}_{(2)}] + (\text{\rm
terms in $\Ecal ^{\geq }_{<n/m} \ctens \Ecal ^{+}_{>n/m}$}),
\end{equation}
where the first term is in $\Lambda (X^{m,n})\otimes \Lambda
(X^{m,n})$, and the subscripts $X^{m,n}_{(1)}$, $X^{m,n}_{(2)}$
distinguish the tensor factors.
\end{prop}

\subsection{Coproduct on the shuffle algebra}
\label{ss:shuffle-coprod}

Recall that for any $\Gamma = \Gamma (w,y)$ satisfying
\eqref{e:Gamma-over-Gamma}, we have an isomorphism $\psi _{\Gamma
}\colon \Scal _{\Gamma }\xrightarrow{\simeq } \Ecal ^{+}$.  
Let
$\Scal _{\Gamma }^{\geq } = \Lambda (X^{0,-1})\otimes \Scal _{\Gamma }$
be the extended algebra
isomorphic to $\Ecal ^{\geq }$ via an isomorphism $\psi _{\Gamma
}^{\geq }\colon \Scal _{\Gamma }^{\geq }\xrightarrow{\simeq } \Ecal
^{\geq }$ that is $\psi _{\Gamma }$ on $\Scal _{\Gamma }$ and the
identity on $\Lambda (X^{0,-1})$.  We then have a coproduct $\Delta
^{\Gamma }$ on $ \Lambda (X^{0,-1})\otimes \Scal _{\Gamma }$
corresponding under $\psi _{\Gamma }^{\geq }$ to the coproduct on
$\Ecal ^{\geq }$.

Negut \cite{Negut14} gives the following formula
(written in our notation) for the component $\Delta ^{\Gamma
}_{k,l-k}$ of $\Delta ^{\Gamma }$ with values in $(\Lambda (X^{0,-1})
\otimes \Scal _{\Gamma }^{k}) \ctens \Scal _{\Gamma }^{l-k}$, when
evaluated on $\Scal _{\Gamma }$. 
The symbol $\ctens $ here indicates that the values are infinite sums
of elements of different degrees in $(\Lambda (X^{0,-1}) \otimes \Scal
_{\Gamma }^{k}) \otimes \Scal _{\Gamma }^{l-k}$, as we explain below
after stating the result.

We remark that, although Negut uses a different $\Gamma $ than we do,
his proof is not specific to the choice.

\begin{prop}[{\cite[Proposition 4.1]{Negut14}}]\label{prop:shuffle-coprod}
Assume that $\Gamma (w,y)$ satisfies \eqref{e:Gamma-over-Gamma} and is
a function of $w/y$.  For any $H(z_{1},\ldots,z_{l})\in \Scal _{\Gamma
}^{l}$, we have
\begin{equation}\label{e:shuffle-coprod}
\Delta^\Gamma_{k,l-k}(H(z_{1},\ldots,z_{l})) =
\omega\Omega[-Y\widehat{M}X^{0,-1}] \frac{H(w_1, \ldots, w_k, y_1,
\ldots, y_{l-k})}{\prod _{i=1}^{k}\prod _{j=1}^{l-k} \Gamma(y_j,w_i)},
\end{equation}
where $Y= y_{1}+\cdots +y_{l-k}$, and we distinguish the factors in
$(\Lambda (X^{0,-1}) \otimes \Scal _{\Gamma }^{k}) \ctens \Scal
_{\Gamma }^{l-k}$ by writing elements of $\Scal _{\Gamma }$ as
functions of $w_{1},\ldots,w_{k}$ in the first tensor factor, or
$y_{1},\ldots,y_{l-k}$ in the second.  Elements of $\Lambda
(X^{0,-1})$ are understood to belong to the first tensor factor.
\end{prop}

More precisely, let $(\Lambda (X^{0,-1}) \otimes \Scal _{\Gamma
}^{k})_{d} = (\psi _{\Gamma }^{\geq })^{-1}((\Ecal ^{\geq })^{(k,d)})$
be the subspace consisting of functions $h(X^{0,-1})f(\ww)$
homogeneous of degree $d$, where $X^{0,-1}$ has degree $-1$, that is,
$h(X^{0,-1})$ has degree $-m$ if $h(X)$ is homogeneous of degree $m$,
and let $(\Scal _{\Gamma }^{l-k})_{d} = \psi _{\Gamma }^{-1}((\Ecal
^{+})^{(l-k,d)})$ be the subspace consisting of functions $g(\yy )$
homogeneous of degree $d$.  The coproduct $\Delta _{k,l-k}^{\Gamma
}(H(\zz ))$ on the left hand side of \eqref{e:shuffle-coprod} is an
infinite sum with components in the spaces $(\Lambda (X^{0,-1})
\otimes \Scal _{\Gamma }^{k})_{d_{1}} \otimes (\Scal _{\Gamma
}^{l-k})_{d_{2}}$ for some set of degrees with $d_{1}$ bounded above
and $d_{2}$ bounded below.

The assumption on $\Gamma (w,y)$ ensures that the right hand side can
be expanded as a formal Laurent series in the $w_{i}^{-1}$ and
$y_{j}$, multiplied by rational functions of $w_{i}/w_{j}$ and
$y_{i}/y_{j}$ (which are thus homogeneous of degree zero in both $\ww
$ and $\yy $).  The meaning of \eqref{e:shuffle-coprod} is that the
component of $\Delta _{k,l-k}^{\Gamma }(H(\zz ))$ in $(\Lambda
(X^{0,-1}) \otimes \Scal _{\Gamma }^{k})_{d_{1}} \otimes (\Scal
_{\Gamma }^{l-k})_{d_{2}}$ is given by the homogeneous component of
these degrees in the series expansion on the right hand side.

Suppose now that $\psi _{\Gamma }(H(\zz )) = f(X^{m,n})$, where $f\in
\Lambda $ is homogeneous of degree $d$.  Then $f(X^{m,n})\in (\Ecal
^{+})^{(dm,dn)}$, and $H(\zz )$ is a function of $l = dm$ variables,
homogeneous of degree $dn$.  Hence, $\Delta ^{\Gamma }_{k,l-k}(H(\zz
))$ has components in $(\Lambda (X^{0,-1}) \otimes \Scal _{\Gamma
}^{k})_{d_{1}} \otimes (\Scal _{\Gamma }^{l-k})_{d_{2}}$ for
$d_{1}+d_{2} = dn$.  Components contributing to the terms of
\eqref{e:leading-term} in $\Ecal ^{\geq }_{<n/m} \ctens \Ecal
^{+}_{>n/m}$
must have $d_{1}<kn/m$, $d_{2}>(l-k)n/m$, while those
contributing to the first term have $d_{1} = kn/m$, $d_{2} =
(l-k)n/m$, necessarily with $k$ a multiple of  $m$.

Since the first term in \eqref{e:leading-term} is in $\Ecal
^+\otimes\Ecal ^+$, all terms of \eqref{e:shuffle-coprod} contributing
to it involve only the constant term in the factor $\omega \Omega [-Y
\widehat{M} X^{0,-1}]$.  These observations yield the following
corollary to Propositions~\ref{prop:leading-term} and
\ref{prop:shuffle-coprod}.

\begin{cor}\label{cor:shuffle-coprod-leading}
Assume that $\Gamma $ satisfies the hypothesis of
Proposition~\ref{prop:shuffle-coprod}.  Suppose that $H(\zz ) =
H(z_{1},\ldots,z_{l})\in \Scal _{\Gamma }^{l}$ has $\psi _{\Gamma
}(H(\zz )) = f(X^{m,n})$, where $f\in \Lambda $ is homogeneous of
degree $d$ and (therefore) $l = dm$.  Given $0\leq e\leq d$, let $k =
e\mkern1mu m$ and let
\begin{equation}\label{e:shuffle-coprod-leading-h}
h(\ww ,\yy ) = h(w_{1},\ldots,w_{k},y_{1},\ldots,y_{l-k}) = \left(
\frac{H(w_1, \ldots, w_k, y_1, \ldots, y_{l-k})}{\prod _{i=1}^{k}\prod
_{j=1}^{l-k} \Gamma(y_j,w_i)} \right)_{\max }
\end{equation}
be the homogeneous component of maximum possible degree $e\mkern1mu n$
in $\ww $ and minimum possible degree $(d-e) n$ in $\yy $.
Regard $h(\ww ,\yy )$ as an element of $\Scal _{\Gamma
}^{k}\otimes \Scal _{\Gamma }^{l-k}$, with variables $\ww $ in
the first tensor factor and $\yy $ in the second.  Then we have
\begin{equation}\label{e:shuffle-coprod-leading}
(\psi _{\Gamma }\otimes \psi _{\Gamma })(h(\ww ,\yy )) = f[X^{m,n}_{
(1)}+X^{m,n}_{(2)}]_{e,d-e}\in \Lambda (X^{m,n})\otimes \Lambda (X^{m,n}),
\end{equation}
where the subscripts $X^{m,n}_{(1)}$, $X^{m,n}_{(2)}$ distinguish the
tensor factors, and $f[X+Y]_{e,d-e}$ designates the homogeneous
component of $f[X+Y]$ of degree $e$ in $X$ and $d-e$ in $Y$.
\end{cor}

\subsection{Coproduct formula for cuddly Catalanimals}
\label{ss:cuddly-coprod}

By Proposition \ref{prop:cuddly}, every $(m,n)$-cuddly Catalanimal has
a cub (Definition~\ref{def:cub}).
Using Corollary~\ref{cor:shuffle-coprod-leading}, we will now obtain a
combinatorial expression for the coproduct of the cub.

We again use the notation $A^{I,J} = \{\alpha_{ij} \in A \mid i \in I,
\, j \in J\}$ from Definition~\ref{def:mn-cuddly}.  Given $A\subseteq
R(\GL _{l})$ and $I\subseteq [l]$ of size $|I| = k$, we also define
$A|_{I}$ to be the set of roots $\{\alpha_{ij} \mid \alpha _{\pi (i)\,
\pi (j)} \in A^{I,I}\}\subseteq R( \GL_{k} )$, where $\pi \colon
[k]\rightarrow I$ is the unique increasing bijection.

\begin{thm}\label{thm:cub-coprod}
Let $H = H(R_q,R_t, R_{qt}, \lambda)$ be an $(m,n)$-cuddly Catalanimal
of length $l = dm$, so its cub has degree $d$.  If $I\subseteq [l]$
attains the cuddliness bound $|\lambda[I]_I| = k n/m$, where
$k = |I|$ is necessarily a multiple of $m$, then the restricted
Catalanimals
\begin{equation}\label{e:H-I}
H_I' = H(R_q|_I, R_t|_I, R_{qt}|_I, \lambda[I]_I), \qquad H_I'' =
H(R_q|_{I^c}, R_t|_{I^c}, R_{qt}|_{I^c}, \lambda[I]_{I^c})
\end{equation}
are $(m,n)$-cuddly, and the coproduct in $\Lambda $ of $f(X) =
\cub(H)$ is given by
\begin{equation}\label{e:cub-coprod}
f[X+Y]= \sum_{I} (-1)^{ |R_+^{I,I^c}| } (-q)^{-|R_q^{I,I^c}|}
(-t)^{-|R_t^{I,I^c}|} (-q\, t)^{|R_{qt}^{I,I^c}|} \cub(H_I')(X)\cdot
\cub(H_I'')(Y),
\end{equation}
where the sum is over subsets $I$ that attain the cuddliness bound.
\end{thm}

\begin{proof}
First, one verifies directly from the definitions the identities
\begin{gather}\label{e:lambda[I]-J}
|(\lambda[I]_I)[J']_{J'}| = |\lambda[J]_J| \qquad \text{for $J \subseteq I$},\\
\label{e:lambda[Ic]-J} |(\lambda[I]_{I^c})[J']_{J'}| = |\lambda[I \cup
J]_{I \cup J}| - |\lambda[I]_I| \qquad \text{for $J \subseteq I^c$},
\end{gather}
where, if $I\subseteq [l]$ has size $|I| = k$, we take $J'\subseteq
[k]$ in \eqref{e:lambda[I]-J} and $J'\subseteq [l-k]$ in
\eqref{e:lambda[Ic]-J} to be the subsets such that $\pi (J') = J$,
where $\pi \colon [k]\rightarrow I$ (resp.\ $\pi \colon
[l-k]\rightarrow I^{c}$) is again the unique increasing bijection.
It then follows from the $(m,n)$-cuddliness of $H$ that $H_I'$ and
$H_I''$ are $(m,n)$-cuddly if $I$ attains the cuddliness bound.

We now apply Corollary~\ref{cor:shuffle-coprod-leading} with the given
Catalanimal $H$ as $H(\zz )$, and $\Gamma = \hGamma (w,y)$, which is
the relevant choice for Catalanimals and their cubs.  Expanding
$\hGamma (y,w)^{-1}$ as a Laurent series in $y/w$ gives
\begin{equation}\label{e:Gamma-Laurent}
\frac{1}{\hGamma (y,w)} = \frac{(1 - w/y) (1 - q\, y/w) (1 - t\,
y/w)}{1- q\, t\, y/w} = -\frac{w}{y}\, (1 + O(y/w)).
\end{equation}
Upon replacing the factor $\prod _{i=1}^{k}\prod _{j=1}^{l-k}
\hGamma(y_{j},w_i)^{-1}$ with its leading term,
\eqref{e:shuffle-coprod-leading-h} simplifies to
\begin{equation}\label{e:h(w,y)-simplified}
h(\ww ,\yy ) = \big( \prod _{i=1}^{k}\prod _{j=1}^{l-k}(-w_{i}/y_{j})
\big)\, H(\ww ,\yy )_{\max }.
\end{equation}
Since $h(\ww ,\yy )$ has degree $kn/m$ in $\ww $ and $(l-k)n/m$ in
$\yy $, the notation $H(\ww ,\yy )_{\max }$ here means the term in
$H(\ww,\yy)$ of degree $kn/m - k(l-k)$ in $\ww $ and $(l-k)n/m+k(l-k)$
in $\yy $.

We now turn to evaluating the leading term $H(\ww ,\yy )_{\max }$ of
$H(\ww ,\yy )$ expanded as a formal Laurent series in the $w_{i}^{-1}$
and $y_{j}$.  Fix a subset $I\subseteq [l]$ of size $|I| = k$, and
consider the terms for which $v(I) = [k]$ in the formula that defines
$H(\ww ,\yy )$, namely,
\begin{equation}\label{e:H(w,y)}
\sum _{v\in S_{l}} v \left(\frac{\zz ^{\lambda }\, \prod _{\alpha \in
R_{qt}}(1 - q\, t\, \zz ^{\alpha })}{\prod _{\alpha \in R_{+}}(1-\zz
^{-\alpha })\, \prod _{\alpha \in R_{q}}(1 - q\, \zz ^{\alpha })\,
\prod _{\alpha \in R_{t}}(1 - t\, \zz ^{\alpha })} \right) \bigg|_{\zz
= (w_{1},\ldots,w_{k},y_{1},\ldots,y_{l-k})}.
\end{equation}
The terms in question are given by evaluating the expression inside
the parentheses with the $z_{i}$ for $i\in I$ specialized to a
permutation of $w_{1},\ldots,w_{k}$, and the $z_{i}$ for $i\in I^{c} =
[l]\setminus I$ to a permutation of $y_{1},\ldots,y_{l-k}$.

When evaluated in this way, each factor $(1 - q\, t\, \zz ^{\alpha })$
with $\alpha \in R_{qt}^{I,I^{c}}$ has leading term $-q\, t\, \zz
^{\alpha }$ since in this case $\zz ^{\alpha }$ evaluates to some
$w_{i}/y_{j}$.  If $\alpha \in R_{qt}^{I,I}\cup R_{qt}^{I^{c},I^{c}}$,
the entire factor becomes homogeneous of degree zero in either $\ww $
or $\yy $.  Otherwise, if $\alpha \in R_{qt}^{I^{c},I}$, the leading
term is $1$.
Similarly, expanding $(1 -q \, \zz ^{\alpha })^{-1}$ as a Laurent
series in $y_{j}/w_{i}$ if $\alpha \not \in R_{q}^{I,I}\cup
R_{q}^{I^{c},I^{c}}$, its leading term is $-q^{-1}\zz ^{-\alpha }$ if
$\alpha \in R_{q}^{I,I^{c}}$, or $1$ otherwise.  The same holds with
$t$ in place of $q$ for factors $(1- t\, \zz ^{\alpha })^{-1}$.
Factors $(1 - \zz ^{-\alpha })^{-1}$ have leading term $1$ if $\alpha
\in R_{+}^{I,I^{c}}$, or $-\zz ^{\alpha }$ if $\alpha \in
R_{+}^{I^{c},I}$.

All of these factors become homogeneous of degree zero if $\alpha \in
R^{I,I}\cup R^{I^{c},I^{c}}$.

Putting all this together, and abbreviating the notation $A^{I,I}$,
$A^{I^{c},I^{c}}$ to $A^{I}$, $A^{I^{c}}$, we find that the
contribution to $H(\ww ,\yy )_{\max }$ from the terms in
\eqref{e:H(w,y)} for a fixed $I$ is given by
\begin{multline}\label{e:I-max-term}
(-1)^{|R_{+}^{I^{c},I}|} (-q)^{-|R_{q}^{I,I^{c}}|}
(-t)^{-|R_{t}^{I,I^{c}}|} (-q\, t)^{|R_{qt}^{I,I^{c}}|}\\
\times \sum _{v(I) = [k]} v \left( \frac{\zz ^{\lambda + \sum
R_{+}^{I^{c},I} - \sum R_{q}^{I,I^{c}} - \sum R_{t}^{I,I^{c}} + \sum
R_{qt}^{I,I^{c}}}\, \prod _{\alpha \in R_{qt}^{I}\cup
R_{qt}^{I^{c}}}(1 - q\, t\, \zz ^{\alpha })}{\prod _{\alpha \in
R_{+}^{I}\cup R_{+}^{I^{c}}}(1-\zz ^{-\alpha })\, \prod _{\alpha \in
R_{q}^{I}\cup R_{q}^{I^{c}}}(1 - q\, \zz ^{\alpha })\, \prod _{\alpha
\in R_{t}^{I} \cup R_{t}^{I^{c}}}(1 - t\, \zz ^{\alpha })}
\right)\bigg|_{\zz = (\ww ,\yy )}
\end{multline}
if the degree of this expression (which is homogeneous) is
$kn/m - k(l-k)$ in $\ww $ and $(l-k)n/m+k(l-k)$ in $\yy $.
Otherwise, its contribution to $H(\ww ,\yy )_{\max} $ is zero.
Observing that
\begin{equation}\label{e:extra-factor}
\begin{aligned}
(-1)^{|R_{+}^{I,I^{c}}| + | R_{+}^{I^{c},I}|} \, v\big( \zz^{\sum
R_{+}^{I,I^{c}} - \sum R_{+}^{I^{c},I}}\big)\Big|_{\zz = (\ww,\yy )} &
= (-1)^{|R^{I,I^{c}}|} \, v(\zz^{\sum R^{I,I^{c}}})\Big|_{\zz = (\ww ,\yy )} \\
& =  \prod_{i=1}^{k}\prod _{j=1}^{l-k}(-w_{i}/y_{j}) ,
\end{aligned}
\end{equation}
we see that the effect of the extra factor in
\eqref{e:h(w,y)-simplified} is to replace $(-1)^{|R_{+}^{I^{c},I}|}$
with $(-1)^{|R_{+}^{I,I^{c}}|}$ and $\zz ^{\sum R_{+}^{I^{c},I}}$ with
$\zz ^{\sum R_{+}^{I,I^{c}}}$ in \eqref{e:I-max-term}, and change the
$\ww $ and $\yy $ degrees back to $kn/m$, $(l-k)n/m$.  The exponent of
$\zz $ then becomes $\lambda [I]$, and \eqref{e:I-max-term} reduces to
\begin{equation}\label{e:restricted-cuddly-term}
(-1)^{|R_{+}^{I,I^{c}}|} (-q)^{-|R_{q}^{I,I^{c}}|}
(-t)^{-|R_{t}^{I,I^{c}}|} (-q\, t)^{|R_{qt}^{I,I^{c}}|} H'_{I}(\ww
)H''_{I}(\yy ),
\end{equation}
which is homogeneous of degree $|\lambda [I]_{I}|$ in $\ww $ and
$|\lambda [I]_{I^{c}}|$ in $\yy $.  These are the maximum (resp.\
minimum) permissible degrees precisely when $I$ attains the cuddliness
bound $|\lambda [I]_{I}| = kn/m$.

Corollary~\ref{cor:shuffle-coprod-leading} now implies that the image
under $\psi _{\hGamma } \otimes \psi _{\hGamma } $ of the expression
in \eqref{e:restricted-cuddly-term}, summed over all $k$ and $I$
attaining the cuddliness bounds, yields
$f[X^{m,n}_{(1)}+X^{m,n}_{(2)}]$.  Identity \eqref{e:cub-coprod} is
this same result expressed in terms of the cubs.
\end{proof}

\section{Principal specialization and evaluating cubs} 
\label{s:principal-specialization}
\setcounter{subsection}{1}

The general strategy to determine the cub $f$ of a cuddly
Catalanimal, justified by the lemma below,
has three steps: show that it exists, determine the inner
terms of its coproduct, and evaluate the specialization $(\omega
f)[1-q]$.  We
completed the first two steps in \S \S \ref{s
cuddly}--\ref{s:coproduct}.  In this section we first give the lemma, 
then establish the required
specialization for the last step, and finally package the resulting
criterion for determining cubs as
Corollary~\ref{cor:evaluating-cubs}.

This strategy is similar to that used by Negut \cite{Negut14} to
establish the special case discussed in Remark \ref{rmk:negut vs
Catalanimal def}.

\begin{lemma}\label{lem:primitive}
If $f \in \Lambda(X)$ is homogeneous of degree $d$, then $f$ is
determined by the terms of the coproduct $f[X_{1}+X_{2}]$ of degree
$k$, $d-k$ in $X_{1}$, $X_{2}$ for $0 < k < d$, together with the
specialization $(\omega f)[1-q]$.
\end{lemma}

\begin{proof}
The terms of degree $k$, $d-k$ for $0 < k < d$ in $f[X_{1}+X_{2}]$
determine $f$ up to adding a primitive (an element $x$ of a Hopf
algebra is primitive if $\Delta(x) = x \otimes 1 + 1 \otimes x$).  By
\cite[Prop. 4.17]{MilMoo65}, the power-sums $p_{j}(X)$ span the vector
space of primitives in $\Lambda (X)$.  Hence, the result follows if the
specialization of $p_{d}(X)$ does not vanish.  But, indeed, $(\omega
p_d)[1-q] = (-1)^{d-1}(1-q^d) \ne 0$.
\end{proof}

Negut \cite[Propositions~6.4,~6.5]{Negut14} gives the value of
elements in the concrete shuffle algebra $\Scal _{\Gamma }$ when
specialized at $\zz =(1,t,\ldots,t^{l-1})$.  The following theorem is
more or less a corollary to Negut's formulas, but we have added what
is needed to express the result in terms of functions $\phi $
representing elements in the abstract shuffle algebra $S$.  Note, in
particular, that since $\phi (z_{1},\ldots,z_{l})$ is not symmetric,
it matters that the powers of $t$ in \eqref{e:ps-main}, below, are in
increasing order.

\begin{thm}\label{thm:principal-spec}
Let $\phi = \phi (z_{1},\ldots,z_{l})\in T^l+I^l$ be a rational
function such that $\psi (\phi ) = f[-MX^{m,n}]$, where $\psi \colon S
= (T+I)/I\rightarrow \Ecal ^{+} $ is the shuffle to Schiffmann
isomorphism in Theorem~\ref{thm:shuffle-isomorphism}, and $f(X)$ is
homogeneous of degree $d$, so $l = dm$.  Assume that the denominator
of $\phi $ does not vanish when evaluated at any permutation of
$(1,t,\ldots,t^{l-1})$.  Then
\begin{equation}\label{e:ps-main}
\phi(1,t, \ldots, t^{l -1}) = \frac{t^a (\omega f)[1-q]}{(1-q)^l},
\end{equation}
where $a = \frac{1}{2}d(d m n - m -n +1)$.
\end{thm}

\begin{proof}
First we show that $\phi (1,t,\ldots,t^{l-1})$ depends only on $\psi
(\phi )$, or equivalently on 
\begin{equation}\label{e:g-for-phi}
g(\zz ) = \psi _{\tGamma }(\phi ) = \sum_{w \in S_{l}} w\big( \phi(\zz
) \prod_{\alpha \in R_+}((1 - \zz ^\alpha)(1 - q\, \zz ^{-\alpha})(1-
t\, \zz ^{-\alpha})(1-q\, t\, \zz ^\alpha)) \big).
\end{equation}
If $w\not =1$, there is some index $i$ such that $j =w^{-1}(i+1)<
w^{-1}(i) = k$.
Then the factor $w(1-t\, z_{k}/z_{j})= (1-t\,
z_{i}/z_{i+1})$ in $w(\prod_{\alpha \in R_+}(1 - t \, \zz
^{-\alpha}))$ vanishes at $\zz =(1,t,\ldots,t^{l-1})$.  By our
assumption on the denominator of $\phi $, the entire $w$ term in
\eqref{e:g-for-phi} vanishes for $w\not =1$, leaving
\begin{equation}\label{e:g-specialized}
g(1,t,\ldots, t^{l -1}) = \phi(1,t,\ldots, t^{l -1}) \prod_{i<j}
\big((1-t^{i-j})(1-q\, t^{j-i})(1-t^{j-i+1})(1-q\, t^{i-j+1})\big).
\end{equation}
The product factor is fixed and non-zero, so $g(1,t,\ldots,t^{l -1})$
determines $\phi(1,t,\ldots,t^{l -1})$.

Next observe that for fixed $m,n$,
\begin{equation}\label{e:a(d)}
a(d) = \frac{1}{2}d(d m n - m -n +1)
\end{equation}
satisfies
\begin{equation}\label{e:a(d)-additivity}
a(d_1+d_2) = a(d_{1})+a(d_{2})+ d_{1}d_{2}mn.
\end{equation}
Given $d_1$, $d_2$ such that $d_1+d_2=d$, suppose that
\eqref{e:ps-main} holds for two functions $\phi
_{1}(z_{1},\ldots,z_{k})$ and $\phi _{2}(z_{1},\ldots,z_{l-k})$ with
corresponding $f_{1}(X)$, $f_{2}(X)$ homogeneous of degrees $d_{1}$,
$d_{2}$.  Then $k = d_{1}m$, $l-k = d_{2}m$, $l=dm$, the functions
$\phi_1$, $\phi_2$ are homogeneous of degrees $d_{1}n$, $d_{2}n$, and
we have
\begin{equation}\label{e:specialization-multiplicativity}
\phi _{1}(1,t,\ldots,t^{k-1})\phi _{2}(t^{k},\ldots,t^{l-1}) =
t^{d_{1}d_{2}mn} \phi _{1}(1,t,\ldots,t^{k-1})\phi
_{2}(1,t,\ldots,t^{l-k-1}).
\end{equation}
The product in $S$ being concatenation, the left hand side of
\eqref{e:specialization-multiplicativity} is $(\phi _{1}\cdot \phi
_{2})(1,t,\ldots,t^{l-1})$.  Using \eqref{e:a(d)-additivity}, it
follows that \eqref{e:ps-main} holds for $\phi = \phi _{1}\cdot \phi
_{2}$ and $f=f_{1}f_{2}$.  Since \eqref{e:ps-main} is also clearly
linear in $\phi $ and $f$, it's enough to prove it when $f = p_{d}$.

Negut \cite[\S 6.3]{Negut14} defines a linear map $\varphi
\colon \Acal ^{+}\rightarrow \kk $ which is characterized by the
property in \cite[Proposition~6.4]{Negut14} and its values $\varphi
(z_{1}^{d}) = (t^{1/2}-t^{-1/2})^{-1}$ on $\Acal ^{+}_{1,d}$.  Using
this one can verify that $\varphi $ corresponds in our notation to the
evaluation map $S\rightarrow \kk $ that sends $\phi \in S^{l}$ to
$\phi (t^{(1-l)/2}, \ldots, t^{(l-1)/2}) / (t^{1/2} - t^{-1/2})^{l}$.
Then, using \cite[Proposition~6.5]{Negut14} and
\eqref{e:Negut-diagram}, one can calculate that for $f = p_{d}$, we
have
$\phi (1,t,\ldots,t^{l-1}) = (-1)^{d-1} t^{a}(1-q^{d})/(1-q)^{l}$.
This agrees with the desired value and completes the proof.
\end{proof}

We can now give a criterion to
determine the cub of a cuddly Catalanimal.

\begin{cor}\label{cor:evaluating-cubs}
Let $H=H(R_q,R_t, R_{qt}, \lambda)$ be an $(m,n)$-cuddly Catalanimal
of length $l = dm$, and let $f \in \Lambda $ be homogeneous of degree
$d$.  To show that $\cub (H) = f$ it suffices to verify the following.

(1) For $0<k<d$, the component $f[X+Y]_{k,d-k}$ of degrees $k$, $d-k$
in $X$, $Y$ is given by
\begin{equation}\label{e:cub-coprod-bis}
\sum_{I} (-1)^{ |R_+^{I,I^c}| } (-q)^{-|R_q^{I,I^c}|}
(-t)^{-|R_t^{I,I^c}|} (-q\, t)^{|R_{qt}^{I,I^c}|} \cub(H_I')(X)\cdot
\cub(H_I'')(Y)\,
\end{equation}
where $H_I'$ and $H_I''$ are as in \eqref{e:H-I}, and the sum is over
index sets $I$ of size $km$ that attain the cuddliness bound $|\lambda
[I]_{I}| = kn$.

(2) The function $\phi(\zz ) = \phi(R_q,R_t,R_{qt}, \lambda)$ defined
in \eqref{eq:phi for Catalanimal} satisfies $\phi(1,t, \ldots, t^{l
-1}) = t^a (\omega f)[1-q]/(1-q)^l$, where $a = \frac{1}{2}d(d m n - m
-n +1)$.
\end{cor}

\begin{proof}
This follows directly from Lemma~\ref{lem:primitive},
Theorem~\ref{thm:cub-coprod}, and Theorem~\ref{thm:principal-spec}.
The only thing to check is the condition on the denominator in
Theorem~\ref{thm:principal-spec}, which clearly holds since the
denominator in this case is a product of factors of the form $(1-q\,
t\, z_{i}/z_{j})$.
\end{proof}

\section{\texorpdfstring{$(1,0)$}{(1,0)}-cuddly LLT Catalanimals}
\label{s:LLT-catalanimals}

For any tuple of skew shapes $\nubold $, we introduce a $(1,0)$-cuddly
Catalanimal $H_{\nubold }$ and prove that its cub is essentially the
LLT polynomial $\Gcal_{\nubold }(X;q)$.  It will be convenient to
establish this case first before turning to the general $(m,n)$ case
in the next section.

\subsection{Definition of the LLT Catalanimals}
\label{ss:def LLT Catalanimals}

We briefly recall the combinatorial concepts used to define LLT
polynomials in \S\ref{ss LLT}.  The adjusted content of a
box $a=(u,v)\in\nu_{(i)}$ in a tuple of skew shapes $\nubold =
(\nu_{(1)},\ldots,\nu_{(k)})$ is $\ctild (a) = u-v+i\epsilon $, and
the reading order is the total ordering on the boxes of $\nubold $
such that $\ctild$ is increasing and boxes on each diagonal increase
from southwest to northeast.  An ordered pair
of boxes $(a,b)$ in $\nubold$ is an attacking pair if $a<b$ in
reading order and $0<\ctild (b)-\ctild(a)<1$.

We let $\nb{1},\ldots,\nb{l}$ denote the boxes of
$\nubold$ in increasing reading order and set
\begin{equation}\label{e:nu(I)}
\nb{I} = \{\nb{i} \mid i\in I \}
\end{equation}
for any subset $I\subseteq [l]$.

\begin{example}\label{ex:LLT Catalanimals}
For the tuple of skew shapes $\nubold = ((32)/(10),(33)/(11))$, the
numbering of boxes in increasing reading order is
\begin{equation*}
\begin{tikzpicture}[scale = .42]
\begin{scope}
\node at (-0.5, 1) {$\Bigg( $};
\draw[thick] (0,1) grid (2,2);
\draw[thick] (1,0) grid (3,1);
\draw[thin, black!60]  (0,0) grid (1,1);
\node at (0.5, 1.5) {1};
\node at (1.5, 1.5) {2};
\node at (1.5, 0.5) {4};
\node at (2.5, 0.5) {7};
\node at (3.34, 0.02) { , };
\end{scope}
\begin{scope}[xshift = 117]
\draw[thick] (1,0) grid (3,2);
\draw[thin, black!60] (0,0) grid (1,2);
\node at (1.5, 1.5) {3};
\node at (2.5, 1.5) {6};
\node at (1.5, 0.5) {5};
\node at (2.5, 0.5) {8};
\node at (3.58, 1) { $ \Bigg)$};
\end{scope}
\end{tikzpicture}
\end{equation*}
\end{example}

The following definition is the special case for $(m,n)=(1,0)$
of a more general construction
in Section~\ref{s mn LLT Catalanimal}.

\begin{defn}
The ($(1,0)$ case) {\em LLT Catalanimal} associated to a tuple of skew shapes
$\nubold $ with a total of $l$ boxes is the length $l$ Catalanimal
$H_{\nubold } = H(R_q,R_t, R_{qt}, \lambda)$, as in
Definition~\ref{def:Catalanimal}, given by the following data:
\begin{align}
\label{ed Rq}
R_q &= \big\{ \alpha_{ij} \in R_+ \mid \ctild(\nb{i}) < \ctild(\nb{j})
 \big\},\\
\label{ed Rt}
R_t &= \big\{ \alpha_{ij} \in R_+ \mid \ctild(\nb{i})+1 \le \ctild(\nb{j})
 \big\}, \\
\label{ed Rqt}
R_{qt} &= \big\{ \alpha_{ij} \in R_+ \mid \ctild(\nb{i})+1 < \ctild(\nb{j})
 \big\}, \\
\label{ed LLT weight}
\lambda_i &=
 \chi(\text{$\diagD{\nb{i}}$ contains the first box in a row})\\
\notag &\quad - \chi(\text{$\diagD{\nb{i}}$ contains the last box in a
row}),
\end{align}
where $\diagD{a}$ is the diagonal of $\nubold $ containing a box $a$,
as in \S \ref{ss LLT},
and $\chi (P)=1$ if $P$ is true, $\chi (P)=0$ otherwise.
\end{defn}

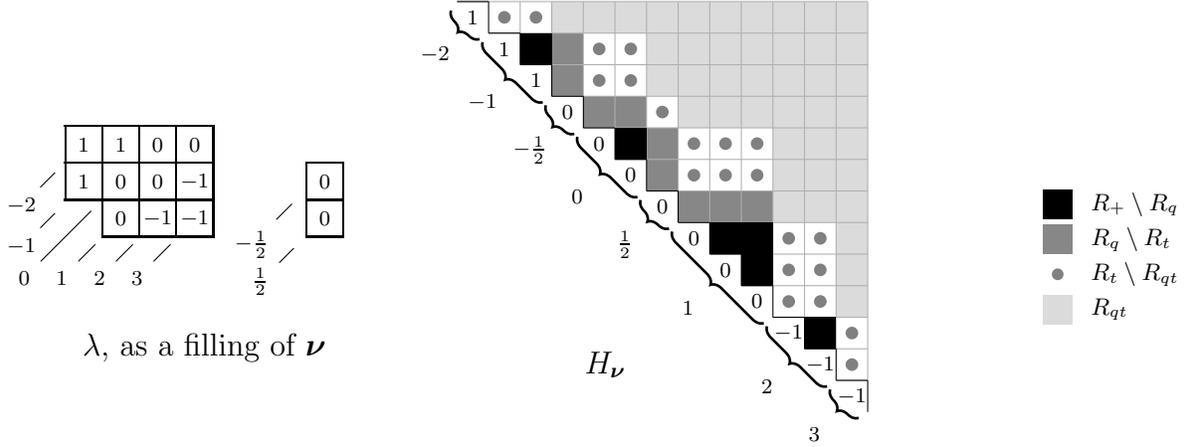
\begin{figure}
\newcounter{ctnt}
\begin{align*}
\begin{tikzpicture}[scale = .42, xshift = 100]
\draw[draw = none, fill = black!14] (4,-1) rectangle (5,-2);
 \draw[draw = none, fill = black!14] (5,-1) rectangle (6,-2);
 \draw[draw = none, fill = black!14] (6,-1) rectangle (7,-2);
 \draw[draw = none, fill = black!14] (7,-1) rectangle (8,-2);
 \draw[draw = none, fill = black!14] (7,-2) rectangle (8,-3);
 \draw[draw = none, fill = black!14] (7,-3) rectangle (8,-4);
 \draw[draw = none, fill = black!14] (8,-1) rectangle (9,-2);
 \draw[draw = none, fill = black!14] (8,-2) rectangle (9,-3);
 \draw[draw = none, fill = black!14] (8,-3) rectangle (9,-4);
 \draw[draw = none, fill = black!14] (8,-4) rectangle (9,-5);
 \draw[draw = none, fill = black!14] (9,-1) rectangle (10,-2);
 \draw[draw = none, fill = black!14] (9,-2) rectangle (10,-3);
 \draw[draw = none, fill = black!14] (9,-3) rectangle (10,-4);
 \draw[draw = none, fill = black!14] (9,-4) rectangle (10,-5);
 \draw[draw = none, fill = black!14] (10,-1) rectangle (11,-2);
 \draw[draw = none, fill = black!14] (10,-2) rectangle (11,-3);
 \draw[draw = none, fill = black!14] (10,-3) rectangle (11,-4);
 \draw[draw = none, fill = black!14] (10,-4) rectangle (11,-5);
 \draw[draw = none, fill = black!14] (11,-1) rectangle (12,-2);
 \draw[draw = none, fill = black!14] (11,-2) rectangle (12,-3);
 \draw[draw = none, fill = black!14] (11,-3) rectangle (12,-4);
 \draw[draw = none, fill = black!14] (11,-4) rectangle (12,-5);
 \draw[draw = none, fill = black!14] (11,-5) rectangle (12,-6);
 \draw[draw = none, fill = black!14] (11,-6) rectangle (12,-7);
 \draw[draw = none, fill = black!14] (11,-7) rectangle (12,-8);
 \draw[draw = none, fill = black!14] (12,-1) rectangle (13,-2);
 \draw[draw = none, fill = black!14] (12,-2) rectangle (13,-3);
 \draw[draw = none, fill = black!14] (12,-3) rectangle (13,-4);
 \draw[draw = none, fill = black!14] (12,-4) rectangle (13,-5);
 \draw[draw = none, fill = black!14] (12,-5) rectangle (13,-6);
 \draw[draw = none, fill = black!14] (12,-6) rectangle (13,-7);
 \draw[draw = none, fill = black!14] (12,-7) rectangle (13,-8);
 \draw[draw = none, fill = black!14] (13,-1) rectangle (14,-2);
 \draw[draw = none, fill = black!14] (13,-2) rectangle (14,-3);
 \draw[draw = none, fill = black!14] (13,-3) rectangle (14,-4);
 \draw[draw = none, fill = black!14] (13,-4) rectangle (14,-5);
 \draw[draw = none, fill = black!14] (13,-5) rectangle (14,-6);
 \draw[draw = none, fill = black!14] (13,-6) rectangle (14,-7);
 \draw[draw = none, fill = black!14] (13,-7) rectangle (14,-8);
 \draw[draw = none, fill = black!14] (13,-8) rectangle (14,-9);
 \draw[draw = none, fill = black!14] (13,-9) rectangle (14,-10);
 \draw[draw = none, fill = black!14] (13,-10) rectangle (14,-11);
 \draw[draw = none, fill = gray!100] (2+0.5, -1-0.5) circle (.2);
\draw[draw = none, fill = gray!100] (3+0.5, -1-0.5) circle (.2);
\draw[draw = none, fill = gray!100] (5+0.5, -2-0.5) circle (.2);
\draw[draw = none, fill = gray!100] (6+0.5, -2-0.5) circle (.2);
\draw[draw = none, fill = gray!100] (5+0.5, -3-0.5) circle (.2);
\draw[draw = none, fill = gray!100] (6+0.5, -3-0.5) circle (.2);
\draw[draw = none, fill = gray!100] (7+0.5, -4-0.5) circle (.2);
\draw[draw = none, fill = gray!100] (8+0.5, -5-0.5) circle (.2);
\draw[draw = none, fill = gray!100] (9+0.5, -5-0.5) circle (.2);
\draw[draw = none, fill = gray!100] (10+0.5, -5-0.5) circle (.2);
\draw[draw = none, fill = gray!100] (8+0.5, -6-0.5) circle (.2);
\draw[draw = none, fill = gray!100] (9+0.5, -6-0.5) circle (.2);
\draw[draw = none, fill = gray!100] (10+0.5, -6-0.5) circle (.2);
\draw[draw = none, fill = gray!100] (11+0.5, -8-0.5) circle (.2);
\draw[draw = none, fill = gray!100] (12+0.5, -8-0.5) circle (.2);
\draw[draw = none, fill = gray!100] (11+0.5, -9-0.5) circle (.2);
\draw[draw = none, fill = gray!100] (12+0.5, -9-0.5) circle (.2);
\draw[draw = none, fill = gray!100] (11+0.5, -10-0.5) circle (.2);
\draw[draw = none, fill = gray!100] (12+0.5, -10-0.5) circle (.2);
\draw[draw = none, fill = gray!100] (13+0.5, -11-0.5) circle (.2);
\draw[draw = none, fill = gray!100] (13+0.5, -12-0.5) circle (.2);
\draw[draw = none, fill = \mymidgray] (4,-2) rectangle (5,-3);
 \draw[draw = none, fill = \mymidgray] (4,-3) rectangle (5,-4);
 \draw[draw = none, fill = \mymidgray] (5,-4) rectangle (6,-5);
 \draw[draw = none, fill = \mymidgray] (6,-4) rectangle (7,-5);
 \draw[draw = none, fill = \mymidgray] (7,-5) rectangle (8,-6);
 \draw[draw = none, fill = \mymidgray] (7,-6) rectangle (8,-7);
 \draw[draw = none, fill = \mymidgray] (8,-7) rectangle (9,-8);
 \draw[draw = none, fill = \mymidgray] (9,-7) rectangle (10,-8);
 \draw[draw = none, fill = \mymidgray] (10,-7) rectangle (11,-8);
 \draw[thin, black!31] (1,-1) -- (14,-1);
\draw[thin, black!31] (2,-1) -- (2,-1);
\draw[thin, black!31] (2,-2) -- (14,-2);
\draw[thin, black!31] (3,-2) -- (3,-1);
\draw[thin, black!31] (3,-3) -- (14,-3);
\draw[thin, black!31] (4,-3) -- (4,-1);
\draw[thin, black!31] (4,-4) -- (14,-4);
\draw[thin, black!31] (5,-4) -- (5,-1);
\draw[thin, black!31] (5,-5) -- (14,-5);
\draw[thin, black!31] (6,-5) -- (6,-1);
\draw[thin, black!31] (6,-6) -- (14,-6);
\draw[thin, black!31] (7,-6) -- (7,-1);
\draw[thin, black!31] (7,-7) -- (14,-7);
\draw[thin, black!31] (8,-7) -- (8,-1);
\draw[thin, black!31] (8,-8) -- (14,-8);
\draw[thin, black!31] (9,-8) -- (9,-1);
\draw[thin, black!31] (9,-9) -- (14,-9);
\draw[thin, black!31] (10,-9) -- (10,-1);
\draw[thin, black!31] (10,-10) -- (14,-10);
\draw[thin, black!31] (11,-10) -- (11,-1);
\draw[thin, black!31] (11,-11) -- (14,-11);
\draw[thin, black!31] (12,-11) -- (12,-1);
\draw[thin, black!31] (12,-12) -- (14,-12);
\draw[thin, black!31] (13,-12) -- (13,-1);
\draw[thin, black!31] (13,-13) -- (14,-13);
\draw[thin, black!31] (14,-13) -- (14,-1);
\draw[draw = none, fill = black!100] (3,-2) rectangle (4,-3);
 \draw[draw = none, fill = black!100] (6,-5) rectangle (7,-6);
 \draw[draw = none, fill = black!100] (9,-8) rectangle (10,-9);
 \draw[draw = none, fill = black!100] (10,-8) rectangle (11,-9);
 \draw[draw = none, fill = black!100] (10,-9) rectangle (11,-10);
 \draw[draw = none, fill = black!100] (12,-11) rectangle (13,-12);
 \draw[thin] (1,-1) -- (2,-1);
\draw[thin] (2,-1) -- (2,-2);
\draw[thin] (2,-2) -- (3,-2);
\draw[thin] (3,-2) -- (3,-3);
\draw[thin] (3,-3) -- (4,-3);
\draw[thin] (4,-3) -- (4,-4);
\draw[thin] (4,-4) -- (5,-4);
\draw[thin] (5,-4) -- (5,-5);
\draw[thin] (5,-5) -- (6,-5);
\draw[thin] (6,-5) -- (6,-6);
\draw[thin] (6,-6) -- (7,-6);
\draw[thin] (7,-6) -- (7,-7);
\draw[thin] (7,-7) -- (8,-7);
\draw[thin] (8,-7) -- (8,-8);
\draw[thin] (8,-8) -- (9,-8);
\draw[thin] (9,-8) -- (9,-9);
\draw[thin] (9,-9) -- (10,-9);
\draw[thin] (10,-9) -- (10,-10);
\draw[thin] (10,-10) -- (11,-10);
\draw[thin] (11,-10) -- (11,-11);
\draw[thin] (11,-11) -- (12,-11);
\draw[thin] (12,-11) -- (12,-12);
\draw[thin] (12,-12) -- (13,-12);
\draw[thin] (13,-12) -- (13,-13);
\draw[thin] (13,-13) -- (14,-13);
\draw[thin] (14,-13) -- (14,-14);
\node at (3/2,-3/2) {\tiny $1 $};
\node at (5/2,-5/2) {\tiny $1 $};
\node at (7/2,-7/2) {\tiny $1 $};
\node at (9/2,-9/2) {\tiny $0 $};
\node at (11/2,-11/2) {\tiny $0 $};
\node at (13/2,-13/2) {\tiny $0 $};
\node at (15/2,-15/2) {\tiny $0 $};
\node at (17/2,-17/2) {\tiny $0 $};
\node at (19/2,-19/2) {\tiny $0 $};
\node at (21/2,-21/2) {\tiny $0 $};
\node at (23/2,-23/2) {\tiny $-1 $};
\node (vv) at (25/2,-25/2) {\tiny $-1 $};
\node at (27/2,-27/2) {\tiny $-1 $};
\node (AA) at (-7,-6.7)
{$ {\fontsize{8pt}{7pt}\selectfont
\tableau{1&1&0&0\\1&0&0&-1\\ \bl &0&-1&-1\\} } \quad
\quad \quad 
{\fontsize{8pt}{7pt}\selectfont \tableau{
\bl \\ 0\\ 0 \\} }  $};
\node[below=1 of AA] {$\lambda$, as a filling of $\nubold  $};
\node[left=21mm of vv] { $H_{\nubold }$};
\setcounter{ctnt}{0};
\foreach \l in {-2,-1}{
  \node[above left = -0.54*\thectnt-1.46 and 0.1cm of AA] (c\thectnt) {\fontsize{8pt}{7pt}\(\l\)};
  \addtocounter{ctnt}{1};
  }
\foreach \l in {0,1,2,3}{
  \node[above left = -2.4 and 0.18-0.5*(\thectnt-2) of AA] (c\thectnt) {\fontsize{8pt}{7pt}\(\l\)};
  \addtocounter{ctnt}{1};
  }
  \foreach \l in {0,1,3,4,5}{
    \draw[thin] (c\l) -- +(45:1.5);
  }
  \draw[thin] (c2)-- +(45:3);
  \node[above right = -2.04 and -1.73 of AA] (c20) {\fontsize{8pt}{7pt}\(-\frac{1}{2}\)};
  \node[above right = -2.54 and -1.5 of AA] (c21) {\fontsize{8pt}{7pt}\(\frac{1}{2}\)};
  \draw[thin] (c20) -- +(45:1.8);
  \draw[thin] (c21) -- +(45:1.5);
  \foreach \l/\bs/\y in {-2/1/1, -1/2/2, 0/2/5, 1/3/8, 2/2/11, 3/1/13}
  {
 \draw[decorate, decoration={brace, mirror, amplitude=4pt,
    raise=4pt},line width=1pt] (\y+0.05,-\y-0.05)--(\y+\bs-0.05,-\y-\bs+0.05)
  node[midway,xshift=-0.5cm,yshift=-0.5cm]{\fontsize{8pt}{7pt}\(\l\)};
  }
 \draw[decorate, decoration={brace, mirror, amplitude=4pt,
    raise=4pt},line width=1pt] (4.05,-4.05)--(4.95,-4.95)
  node[midway,xshift=-0.5cm,yshift=-0.5cm]{\fontsize{8pt}{7pt}\(-\frac{1}{2}\)};
 \draw[decorate, decoration={brace, mirror, amplitude=4pt,
    raise=4pt},line width=1pt] (7.05,-7.05)--(7.95,-7.95)
  node[midway,xshift=-0.5cm,yshift=-0.5cm]{\fontsize{8pt}{7pt}\(\frac{1}{2}\)};
\end{tikzpicture}
\quad \qquad \qquad
\raisebox{16mm}{
\begin{tikzpicture}[scale = .39]
\begin{scope}
\draw[draw = none, fill = black!100] (2,-1) rectangle (3,-2);
\node[anchor = west] at (3,-1.5) {\scriptsize  \  $R_+ \setminus R_q$};
\end{scope}
\begin{scope}[yshift = -34*1]
\draw[draw = none, fill = \mymidgray] (2,-1) rectangle (3,-2);
\node[anchor = west] at (3,-1.5) {\scriptsize  \  $R_q \setminus R_t$};
\end{scope}
\begin{scope}[yshift = -34*2]
\draw[thin, black!0] (2,-1) -- (3,-1);
\draw[thin, black!0] (2,-1) -- (2,-2);
\draw[thin, black!0] (2,-2) -- (3,-2);
\draw[thin, black!0] (3,-2) -- (3,-1);
\draw[draw = none, fill = gray!100] (2+0.5, -1-0.5) circle (.2);
\node[anchor = west] at (3,-1.5) {\scriptsize  \  $R_t \setminus R_{qt}$};
\end{scope}
\begin{scope}[yshift = -34*3]
\draw[draw = none, fill = black!14] (2,-1) rectangle (3,-2);
\node[anchor = west] at (3,-1.5) {\scriptsize \  $R_{qt}$};
\end{scope}
\end{tikzpicture}}
\end{align*}
\caption{\label{fig:LLT Catalanimal}%
The LLT Catalanimal $H_{\nubold}$ for $\nubold = ((444)/(1),(11))$.
We have marked the adjusted contents of $\lambda$ for $\epsilon =
1/2$, along with the corresponding parabolic blocks in the
Catalanimal.}
\end{figure}

\begin{remark}\label{rem:LLT Catalanimal def}
(i) The root sets in \eqref{ed Rq}--\eqref{ed LLT weight} satisfy $R_+
\supseteq R_q \supseteq R_t \supseteq R_{qt}$.  It is convenient to
think of these as providing a partition of the positive roots into the
following four subsets defined by the combinatorial features of the
LLT polynomial $\Gcal_{\nubold }$:
\begin{align*}
R_+ \setminus R_q & = \{\alpha _{ij} \in R_+ \mid \text{$\nb{i}$,
 $\nb{j}$ are on the same diagonal} \},\\
R_q \setminus R_t & = \{\alpha _{ij} \in R_+ \mid \text{$\nb{i}$,
 $\nb{j}$ form an attacking pair} \},\\
R_t \setminus R_{qt} & = \{\alpha _{ij} \in R_+ \mid \text{$\nb{i}$,
 $\nb{j}$ are on adjacent diagonals} \},\\
R_{qt} & = \{\text{all other } \alpha _{ij} \in R_+ \}.
\end{align*}
We say that diagonals in $\nubold $ are {\em adjacent} if their
adjusted contents differ by $1$, that is, they are in the same skew
shape $\nu _{(r)}$ and their ordinary contents differ by $1$.

(ii) The data $R_{q}$, $R_{t}$, $R_{qt}$ and $\lambda $ are constant
on diagonals of $\nubold $, in the sense that $\lambda _{i}$ depends
only on $\diagD{\nb{i}}$,
and whether or not $\alpha _{ij}$ belongs to $R_{q}$ depends only on 
$\diagD{\nb{i}} $ and $ \diagD{\nb{j}}$; likewise for $R_{t}$ and $R_{qt}$.

(iii) One way to picture the weight $\lambda $ is as a filling of
$\nubold $ with $\lambda _{i}$ in box $\nb{i}$, so that $\lambda $ is
the list of labels in the filling in reading order.  Viewed as a
filling, $\lambda $ is constant on each diagonal $D\subseteq \nubold
$, with value $\pm 1$ or $0$ depending on whether the boxes at the
southwest and northeast ends of $D$ are the first or last boxes in
their row.  This is illustrated in Figures~\ref{fig:Catalanimal 444
minus 1}~and~\ref{fig:LLT Catalanimal}.
\end{remark}

\begin{figure}
\[
\begin{tikzpicture}[scale=.42]
\draw[draw = none, fill = black!14] (4,-1) rectangle (5,-2);
 \draw[draw = none, fill = black!14] (5,-1) rectangle (6,-2);
 \draw[draw = none, fill = black!14] (6,-1) rectangle (7,-2);
 \draw[draw = none, fill = black!14] (6,-2) rectangle (7,-3);
 \draw[draw = none, fill = black!14] (6,-3) rectangle (7,-4);
 \draw[draw = none, fill = black!14] (7,-1) rectangle (8,-2);
 \draw[draw = none, fill = black!14] (7,-2) rectangle (8,-3);
 \draw[draw = none, fill = black!14] (7,-3) rectangle (8,-4);
 \draw[draw = none, fill = black!14] (8,-1) rectangle (9,-2);
 \draw[draw = none, fill = black!14] (8,-2) rectangle (9,-3);
 \draw[draw = none, fill = black!14] (8,-3) rectangle (9,-4);
 \draw[draw = none, fill = black!14] (9,-1) rectangle (10,-2);
 \draw[draw = none, fill = black!14] (9,-2) rectangle (10,-3);
 \draw[draw = none, fill = black!14] (9,-3) rectangle (10,-4);
 \draw[draw = none, fill = black!14] (9,-4) rectangle (10,-5);
 \draw[draw = none, fill = black!14] (9,-5) rectangle (10,-6);
 \draw[draw = none, fill = black!14] (10,-1) rectangle (11,-2);
 \draw[draw = none, fill = black!14] (10,-2) rectangle (11,-3);
 \draw[draw = none, fill = black!14] (10,-3) rectangle (11,-4);
 \draw[draw = none, fill = black!14] (10,-4) rectangle (11,-5);
 \draw[draw = none, fill = black!14] (10,-5) rectangle (11,-6);
 \draw[draw = none, fill = black!14] (11,-1) rectangle (12,-2);
 \draw[draw = none, fill = black!14] (11,-2) rectangle (12,-3);
 \draw[draw = none, fill = black!14] (11,-3) rectangle (12,-4);
 \draw[draw = none, fill = black!14] (11,-4) rectangle (12,-5);
 \draw[draw = none, fill = black!14] (11,-5) rectangle (12,-6);
 \draw[draw = none, fill = black!14] (11,-6) rectangle (12,-7);
 \draw[draw = none, fill = black!14] (11,-7) rectangle (12,-8);
 \draw[draw = none, fill = black!14] (11,-8) rectangle (12,-9);
 \draw[draw = none, fill = gray!100] (2+0.5, -1-0.5) circle (.2);
\draw[draw = none, fill = gray!100] (3+0.5, -1-0.5) circle (.2);
\draw[draw = none, fill = gray!100] (4+0.5, -2-0.5) circle (.2);
\draw[draw = none, fill = gray!100] (5+0.5, -2-0.5) circle (.2);
\draw[draw = none, fill = gray!100] (4+0.5, -3-0.5) circle (.2);
\draw[draw = none, fill = gray!100] (5+0.5, -3-0.5) circle (.2);
\draw[draw = none, fill = gray!100] (6+0.5, -4-0.5) circle (.2);
\draw[draw = none, fill = gray!100] (7+0.5, -4-0.5) circle (.2);
\draw[draw = none, fill = gray!100] (8+0.5, -4-0.5) circle (.2);
\draw[draw = none, fill = gray!100] (6+0.5, -5-0.5) circle (.2);
\draw[draw = none, fill = gray!100] (7+0.5, -5-0.5) circle (.2);
\draw[draw = none, fill = gray!100] (8+0.5, -5-0.5) circle (.2);
\draw[draw = none, fill = gray!100] (9+0.5, -6-0.5) circle (.2);
\draw[draw = none, fill = gray!100] (10+0.5, -6-0.5) circle (.2);
\draw[draw = none, fill = gray!100] (9+0.5, -7-0.5) circle (.2);
\draw[draw = none, fill = gray!100] (10+0.5, -7-0.5) circle (.2);
\draw[draw = none, fill = gray!100] (9+0.5, -8-0.5) circle (.2);
\draw[draw = none, fill = gray!100] (10+0.5, -8-0.5) circle (.2);
\draw[draw = none, fill = gray!100] (11+0.5, -9-0.5) circle (.2);
\draw[draw = none, fill = gray!100] (11+0.5, -10-0.5) circle (.2);
\draw[thin, black!31] (1,-1) -- (12,-1);
\draw[thin, black!31] (2,-1) -- (2,-1);
\draw[thin, black!31] (2,-2) -- (12,-2);
\draw[thin, black!31] (3,-2) -- (3,-1);
\draw[thin, black!31] (3,-3) -- (12,-3);
\draw[thin, black!31] (4,-3) -- (4,-1);
\draw[thin, black!31] (4,-4) -- (12,-4);
\draw[thin, black!31] (5,-4) -- (5,-1);
\draw[thin, black!31] (5,-5) -- (12,-5);
\draw[thin, black!31] (6,-5) -- (6,-1);
\draw[thin, black!31] (6,-6) -- (12,-6);
\draw[thin, black!31] (7,-6) -- (7,-1);
\draw[thin, black!31] (7,-7) -- (12,-7);
\draw[thin, black!31] (8,-7) -- (8,-1);
\draw[thin, black!31] (8,-8) -- (12,-8);
\draw[thin, black!31] (9,-8) -- (9,-1);
\draw[thin, black!31] (9,-9) -- (12,-9);
\draw[thin, black!31] (10,-9) -- (10,-1);
\draw[thin, black!31] (10,-10) -- (12,-10);
\draw[thin, black!31] (11,-10) -- (11,-1);
\draw[thin, black!31] (11,-11) -- (12,-11);
\draw[thin, black!31] (12,-11) -- (12,-1);
\draw[draw = none, fill = black!100] (3,-2) rectangle (4,-3);
 \draw[draw = none, fill = black!100] (5,-4) rectangle (6,-5);
 \draw[draw = none, fill = black!100] (7,-6) rectangle (8,-7);
 \draw[draw = none, fill = black!100] (8,-6) rectangle (9,-7);
 \draw[draw = none, fill = black!100] (8,-7) rectangle (9,-8);
 \draw[draw = none, fill = black!100] (10,-9) rectangle (11,-10);
 \draw[thin] (1,-1) -- (2,-1);
\draw[thin] (2,-1) -- (2,-2);
\draw[thin] (2,-2) -- (3,-2);
\draw[thin] (3,-2) -- (3,-3);
\draw[thin] (3,-3) -- (4,-3);
\draw[thin] (4,-3) -- (4,-4);
\draw[thin] (4,-4) -- (5,-4);
\draw[thin] (5,-4) -- (5,-5);
\draw[thin] (5,-5) -- (6,-5);
\draw[thin] (6,-5) -- (6,-6);
\draw[thin] (6,-6) -- (7,-6);
\draw[thin] (7,-6) -- (7,-7);
\draw[thin] (7,-7) -- (8,-7);
\draw[thin] (8,-7) -- (8,-8);
\draw[thin] (8,-8) -- (9,-8);
\draw[thin] (9,-8) -- (9,-9);
\draw[thin] (9,-9) -- (10,-9);
\draw[thin] (10,-9) -- (10,-10);
\draw[thin] (10,-10) -- (11,-10);
\draw[thin] (11,-10) -- (11,-11);
\draw[thin] (11,-11) -- (12,-11);
\draw[thin] (12,-11) -- (12,-12);
\node at (3/2,-3/2) {\tiny $1 $};
\node at (5/2,-5/2) {\tiny $1 $};
\node at (7/2,-7/2) {\tiny $1 $};
\node at (9/2,-9/2) {\tiny $0 $};
\node at (11/2,-11/2) {\tiny $0 $};
\node at (13/2,-13/2) {\tiny $0 $};
\node at (15/2,-15/2) {\tiny $0 $};
\node at (17/2,-17/2) {\tiny $0 $};
\node at (19/2,-19/2) {\tiny $-1 $};
\node (vv) at (21/2,-21/2) {\tiny $-1 $};
\node at (23/2,-23/2) {\tiny $-1 $};
\node (AA) at (-7,-5.7) {$ {\fontsize{8pt}{7pt}\selectfont
\tableau{1&1&0&0\\1&0&0&-1\\ \bl &0&-1&-1\\} } $}; \node[below=0mm of
AA] {$\lambda$, as a filling of $\nubold $}; \node[left=14.2mm of vv]
{ $H_{\nubold } $};
\end{tikzpicture}
\quad \qquad \qquad \qquad
\raisebox{16mm}{
\begin{tikzpicture}[scale = .39]
\begin{scope}
\draw[draw = none, fill = black!100] (2,-1) rectangle (3,-2);
\node[anchor = west] at (3,-1.5) {\scriptsize  \  $R_+ \setminus R_q$};
\end{scope}
\begin{scope}[yshift = -34*1]
\draw[draw = none, fill = \mymidgray] (2,-1) rectangle (3,-2);
\node[anchor = west] at (3,-1.5) {\scriptsize  \  $R_q \setminus R_t$};
\end{scope}
\begin{scope}[yshift = -34*2]
\draw[thin, black!0] (2,-1) -- (3,-1);
\draw[thin, black!0] (2,-1) -- (2,-2);
\draw[thin, black!0] (2,-2) -- (3,-2);
\draw[thin, black!0] (3,-2) -- (3,-1);
\draw[draw = none, fill = gray!100] (2+0.5, -1-0.5) circle (.2);
\node[anchor = west] at (3,-1.5) {\scriptsize  \  $R_t \setminus R_{qt}$};
\end{scope}
\begin{scope}[yshift = -34*3]
\draw[draw = none, fill = black!14] (2,-1) rectangle (3,-2);
\node[anchor = west] at (3,-1.5) {\scriptsize \  $R_{qt}$};
\end{scope}
\end{tikzpicture}}
\]
\caption{\label{fig:Catalanimal 444 minus 1}%
The LLT Catalanimal $H_{\bm{\nu}}$ for $\bm{\nu}$ the single skew
shape $(444)/(1)$, drawn with same conventions as in Figure
\ref{fig:Catalanimals intro}.  }
\end{figure}
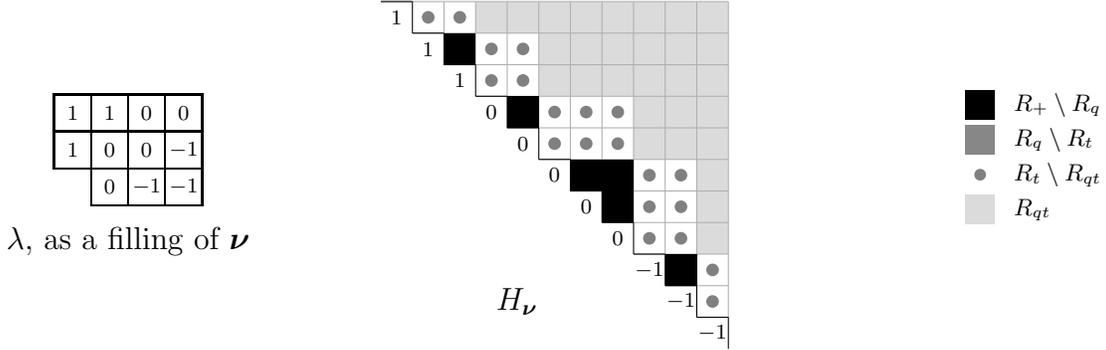

\begin{example}\label{ex:some-LLT-cats}
(i) If $\nubold $ is a single straight shape $((111))$, then $R_q =
R_t = R_+ = \{\alpha_{12}, \alpha_{13}, \alpha_{23}\}$, $R_{qt} =
\{\alpha_{13}\}$, and $\lambda = (000)$.  This gives
\[
H_{\nubold } =
\sum_{w \in S_3} w \Big(\frac{ (1 - q\, t\, z_{1} /z_3)}
{\prod_{1\le i < j \le 3}(1-z_j/z_i) (1 - q\, z_{i}/ z_{j})
(1 - t\,z_{i}/ z_{j})} \Big),\quad \text{drawn as}\quad 
\begin{tikzpicture}[scale = .44,baseline=-1.2cm]
\draw[draw = none, fill = black!14] (3,-1) rectangle (4,-2);
\draw[draw = none, fill = gray!100] (2+0.5, -1-0.5) circle (.2); 
\draw[draw = none, fill = gray!100] (3+0.5, -2-0.5) circle (.2); 
\draw[thin, black!31] (1,-1) -- (4,-1);
\draw[thin, black!31] (2,-1) -- (2,-1);
\draw[thin, black!31] (2,-2) -- (4,-2);
\draw[thin, black!31] (3,-2) -- (3,-1);
\draw[thin, black!31] (3,-3) -- (4,-3);
\draw[thin, black!31] (4,-3) -- (4,-1);
\draw[thin] (1,-1) -- (2,-1);
\draw[thin] (2,-1) -- (2,-2);
\draw[thin] (2,-2) -- (3,-2);
\draw[thin] (3,-2) -- (3,-3);
\draw[thin] (3,-3) -- (4,-3);
\draw[thin] (4,-3) -- (4,-4);
\node at (3/2,-3/2) {\tiny $0 $}; 
\node at (5/2,-5/2) {\tiny $0 $}; 
\node at (7/2,-7/2) {\tiny $0 $}; 
\end{tikzpicture}.
\]

(ii) If $\nubold $ is a single ribbon skew shape $(\theta )$, then
$R_q = R_t = R_+$ and $R_{qt} = [R_+,R_+] = \{\alpha_{ij} \in R_+ \mid
j-i > 1\}$.  The weight $\lambda$, viewed as a filling of $\theta $ as
in Remark~\ref{rem:LLT Catalanimal def} (iii), is obtained by writing
$1,0, \dots, 0, -1$ across each row, or $0$ if the row has one box.
This case relates to prior work of Negut as will be 
discussed in Remark~\ref{rmk:negut vs
Catalanimal def}.

(iii) More generally, if $\nubold $ is an arbitrary single skew shape
$(\theta )$, then $R_q = R_t$ is the set of roots corresponding to
positions above the diagonal blocks in a block matrix with block sizes
equal to the lengths of the diagonals in $\theta $, and $R_{qt}$
corresponds to the subset of these blocks obtained by removing the
blocks immediately above the diagonal, as shown for $\theta =
(444)/(1)$ in Figure \ref{fig:Catalanimal 444 minus 1}.  See also
Figure \ref{fig:Catalanimals intro} (i) for the case $\theta = (433)$
(but note that the LLT Catalanimal $H^{1,1}_{\nubold }$ shown there is
obtained from that of $H_{\nubold }$ by adding $1^l$ to the weight, as
will be explained in \S\ref{s mn LLT Catalanimal}).

(iv) The case $\nubold = ((444)/(1),(11))$ in Figure \ref{fig:LLT
Catalanimal} illustrates the general construction.
\end{example}

\begin{defn}\label{d:prec order}
Define a partial order $\prec $ on boxes in $\nubold $ by setting $a
\preceq b$ if boxes $a,b$ belong to the same skew shape $\nu_{(i)}$
and
$a$ is weakly southwest of $b$, i.e.,
$a\leq b$ in the usual product order on $\NN \times \NN$.  Let
$R_+^\prec = \{\alpha_{ij} \in R_+ \mid \nb{i} \prec \nb{j}\}$.
\end{defn}

\begin{prop}\label{p:weight defs}
The weight $\lambda$ in \eqref{ed LLT weight} has the following
alternative descriptions:
\begin{multline}\label{ed LLT weight 2}
\lambda_i = \chi(\text{\rm $\diagD{\nb{i}}$ does not contain the last
box in a row})\\
 - \chi(\text{\rm $\diagD{\nb{i}}$ does not contain the first box in a
row});
\end{multline}
\begin{equation}\label{ed LLT weight 3}
\lambda = - \sum (R_+\setminus R_q) + \sum ((R_t \setminus R_{qt})
\cap R_+^\prec)
\end{equation}
\end{prop}

\begin{proof}
The first description is just a reformulation of \eqref{ed LLT
weight}.  For the second, let $\mu$ denote the right side of \eqref{ed LLT
weight 3}.  We compute the contribution to $\mu_i$ from each of the
sums.
Let $D = \diagD{\nb{i}}$ and let $C$, $E$ be the (possibly
empty) adjacent diagonals with adjusted contents $\ctild(\nb{i})-1$,
$\ctild(\nb{i})+1$, respectively.
Then, setting $\nubold _{\succ i} = \{a \in \nubold \mid a \succ
\nb{i} \}$, $\nubold _{\prec i} = \{a \in \nubold \mid a \prec \nb{i}
\}$, we have
\begin{gather}\label{e:llt-weight-3a}
\mu_i = -|D \cap \nubold _{\succ i}|+|D \cap \nubold _{\prec i}| -| C
\cap \nubold _{\prec i}|+|E \cap \nubold _{\succ i}|.
\end{gather}
The middle two terms sum to $-\chi(\text{$D$ does not contain the
first box in a row})$, while the first and last terms sum to
$\chi(\text{$D$ does not contain the last box in a row})$, so this matches
\eqref{ed LLT weight 2}.
\end{proof}

\subsection{Statistics on  \texorpdfstring{$\nubold$}{nu}}
\label{ss statitistics on nu}

In preparation for determining the cubs of LLT Catalanimals, we
require the following statistics associated to a tuple of skew shapes
$\nubold $:
\begin{align}\label{eqdef:gamma nu}
 \gamma(\nubold ) &  \defeq  \text{sequence of lengths of
  diagonals in $\nubold $, in increasing reading order}; \\
\label{eqdef:n gamma}
  n'(\gamma) &  \defeq  \sum_{i} \binom{\gamma_i}{2}
  \quad \text{for any $\gamma$, but chiefly used for
   $n'(\gamma (\nubold ))$}; \\
\label{eqdef:magic num} p(\nubold ) & \defeq \sum_{\text{diagonals
$D\subseteq \nubold $}}
 \chi(\text{$D$ does not contain the first box in a row})\cdot |D|\, ; \\
A(\nubold ) &  \defeq  \text{number of attacking pairs in $\nubold $}.
\end{align}
We also refer to $p(\nubold )$ as the \emph{magic number} of
$\nubold $.

\begin{example}\label{ex:LLT Cat stats} The statistics associated to
the LLT Catalanimals of Figures \ref{fig:Catalanimals intro},
\ref{fig:Catalanimal 444 minus 1}, \ref{fig:LLT Catalanimal} are as
follows.

\ytableausetup{boxsize=0.325em,mathmode,centertableaux}

\smallskip

\noindent Figure \ref{fig:Catalanimals intro} (i), $\nubold = \left(\,
\ydiagram{3,3,4}\, \right)$: $p(\nubold ) = 4$, \ $\gamma(\nubold ) =
( 1, 2, 3, 2, 1, 1)$, \ $n'(\gamma(\nubold )) = 5$, \ $A(\nubold ) =
0$.

\smallskip

\noindent Figure \ref{fig:Catalanimals intro} (ii), $\nubold = \left(
\,\ydiagram{2,3}*[*(black)]{0,1}\,,\ydiagram{3,3}*[*(black)]{1,1}\,
\right)$: $p(\nubold ) = 3 $, \ $\gamma(\nubold ) = ( 1, 1, 1, 1, 2,
1,1 ) $, \ $n'(\gamma(\nubold )) = 1 $, \ $A(\nubold ) \, = 7 $.

\smallskip

\noindent Figure \ref{fig:Catalanimal 444 minus 1}, $\nubold = \left(
\,\ydiagram{4,4,4}*[*(black)]{0,0,1}\, \right)$: $p(\nubold ) = 5$, \
$\gamma(\nubold ) = ( 1, 2, 2,3, 2, 1 )$, \ $n'(\gamma(\nubold )) = 6
$, \ $A(\nubold ) = 0$.

\smallskip

\noindent Figure \ref{fig:LLT Catalanimal}, $\nubold = \left(
\,\ydiagram{4,4,4}*[*(black)]{0,0,1}\,, \ydiagram{0,1,1} \, \right)$:
$p(\nubold ) = 5$, \ $\gamma(\nubold ) = ( 1, 2,1, 2,1, 3, 2, 1 )$, \
$n'(\gamma(\nubold )) = 6 $, \ $A(\nubold ) = 9$.
\end{example}

\begin{lemma}\label{lem:magic-num}
For each diagonal $D$ in $\nubold $, let $D_{+}$ denote the (possibly
empty) adjacent diagonal 
southeast of $D$, that is, the adjusted
contents of boxes $d\in D$ and $e\in D_{+}$ satisfy $\ctild (e) =
\ctild (d)+1$. Then
\begin{equation}\label{e:eta-plus-magic}
n'(\gamma (\nubold ))+p(\nubold ) = \sum _{D} |\{(d,e)\in D\times
D_{+}\mid d \prec e \}|,
\end{equation}
summed over all diagonals $D$ in $\nubold $.
\end{lemma}

\begin{proof}
Letting $e' = \south (d)$, the sum in \eqref{e:eta-plus-magic} is
almost the same as the number $n'(\gamma (\nubold ))$ of pairs of
boxes $(e',e)$ such that $e'$ and $e$ are on the same diagonal and
$e'\prec e$.  The only difference is that \eqref{e:eta-plus-magic}
also counts pairs $(d,e)$ for which $\south (d)$ is not in $\nubold $.
But this means that the diagonal $D_{+}$ containing $e$ does not
contain the first box in a row, and that $d = \west (e_{1})$, where
$e_{1}$ is the first box of $D_{+}$ in reading order.  Hence, the
number of pairs not counted by $n'(\gamma (\nubold ))$ is $p(\nubold
)$.
\end{proof}

\subsection{Proof of cuddliness and determining the cubs}
\label{ss:LLT-cuddliness-and-cubs}

Before determining the cubs of LLT Catalanimals, we observe
that some formulas involved in Theorem~\ref{thm:cub-coprod} and
Corollary~\ref{cor:evaluating-cubs} simplify when
$R_{q}\supseteq R_{t}\supseteq R_{qt}$, as is the case for LLT
Catalanimals. 
Here and below it will be convenient to use the
abbreviations
\begin{equation}\label{e:bar-Rq}
\barRq = R_{+}\setminus R_{q},\quad \barRt = R_{q}\setminus
R_{t},\quad \barRqt = R_{t}\setminus R_{qt}.
\end{equation}
The adjusted weight $\lambda [I]$ in \eqref{eqdef:S bracket I}
simplifies to
\begin{equation}\label{e:lambda-I-simple}
\lambda[I] = \lambda + \sum \barRq ^{I,I^c} - \sum \barRqt ^{I,I^c}
\end{equation}
and the quantity $|\lambda [I]_{I}|$ in the cuddliness bound
\eqref{e:cuddly-bound} to
\begin{equation}\label{e:cuddly-weight}
|\lambda [I]_{I}| = \sum_{i \in I}\lambda_i + |\barRq ^{I,I^c}| -
|\barRqt ^{I,I^c}|.
\end{equation}

Formula \eqref{e:cub-coprod-bis} in
Corollary~\ref{cor:evaluating-cubs} part (1), which is also the right
hand side of \eqref{e:cub-coprod}, can be rewritten as
\begin{equation}\label{e:cub-coprod-q-style}
\sum_{I} (-1)^{ |\barRq ^{I,I^c}| + |\barRqt ^{I,I^c}| }
(qt)^{-|\barRqt ^{I,I^c}|} q^{-|\barRt ^{I,I^c}|}
\cub(H_I')(X) \cdot \cub(H_I'')(Y),
\end{equation}
where the sum is still over index sets $I\subseteq [l]$ of size $km$
attaining the cuddliness bound $|\lambda [I]_{I}| = kn$, and $H_{I}'$,
$H_{I}''$ are as in \eqref{e:H-I}.

To invoke Corollary~\ref{cor:evaluating-cubs} part (2) we will need
the following lemma, which is the special case for $(m,n) = (1,0)$ of
Lemma~\ref{lem:LLT ps}, below---see Remark~\ref{rem:LLT-1-0-ps}.

\begin{lemma}\label{lem:LLT-1-0-ps}
The principal specialization of the function $\phi(\zz ) =
\phi(R_q,R_t,R_{qt}, \lambda)$ in \eqref{eq:phi for Catalanimal}
associated to the LLT Catalanimal $H_{\nubold } = H(R_q,R_t,R_{qt},
\lambda)$ is given by
\begin{equation}\label{e:LLT-1-0-ps}
\phi(1,t, \ldots, t^{l -1}) = (\omega f)[1-q]/(1-q)^l,
\end{equation}
where $f = (-1)^{p(\nubold )}(qt)^{-p(\nubold )-n'(\gamma (\nubold ))}
q^{-A(\nubold )} \Gcal_{\nubold }(X;q)$.
\end{lemma}

\begin{thm}
\label{t LLT cuddly} For any tuple of skew shapes $\nubold $, the LLT
Catalanimal $H_{\nubold }= H(R_q,R_t,R_{qt},\lambda)$ is
$(1,0)$-cuddly with cub related to the corresponding LLT polynomial by
\begin{align}\label{eq:LLT cuddly}
\cub(H_{\nubold }) = (-1)^{p}(q\, t)^{-p-n'(\gamma)} q^{-A}
\Gcal_{\nubold }(X;q),
\end{align}
where $A = A(\nubold )$ is the number of attacking pairs in $\nubold
$, $p = p(\nubold )$ is the magic number, and $\gamma = \gamma(\nubold
)$ are the diagonal lengths.
\end{thm}

\begin{proof} We first verify that $H_{\nubold }$ is $(1,0)$-cuddly,
then use Corollary~\ref{cor:evaluating-cubs} to determine its cub.

\smallskip

\noindent {\em Checking $(1,0)$-cuddliness.}  We start by checking the
tameness condition $[R_q,R_t] \subseteq R_{qt}$.  For $\alpha_{ij} \in
R_q$ and $\alpha_{jk} \in R_t$, \eqref{ed Rq} and \eqref{ed Rt} give
$\ctild(\nb{i}) < \ctild(\nb{j}) \le \ctild(\nb{k})-1$, which ensures
$\alpha_{ik} \in R_{qt}$ by \eqref{ed Rqt}, as desired.  A similar
argument works with the roles of $R_q$ and $R_t$ interchanged.

We next prove that for all $I \subseteq [l]$, we have
\begin{equation}\label{e LLT cuddly induction}
\text{
$|\lambda [I]_{I}|\le 0$,
with equality if and only if $\nb{I}$ is a lower order ideal for $\prec$}.
\end{equation}
This establishes the cuddliness bounds \eqref{e:cuddly-bound} for
$H_{\nubold }$ and completes the proof that it is $(1,0)$-cuddly.  The
second part of the claim will help later on to determine the cub.

In verifying \eqref{e LLT cuddly induction} it will be
convenient to use the abbreviation
\begin{equation}\label{e:r-I}
r_{I} = |\lambda [I]_{I}|.
\end{equation}
Let $I \subseteq [l]$.  First suppose there is a pair $x-1,\, x \in
[l]$ such that $x \in I$, $x-1 \notin I$, and $\nb{x-1}$ and $\nb{x}$
belong to the same diagonal.  Set $J = I \cup \{x-1\} \setminus
\{x\}$.
Then by \eqref{e:cuddly-weight},
\begin{equation}\label{e:rJ-rI}
r_{J} -r_I = \lambda_{x-1}-\lambda_x + |\barRq ^{J,J^c}| - |\barRq
^{I,I^c}| - |\barRqt ^{J,J^c}| + |\barRqt ^{I,I^c}| = 1,
\end{equation}
where for the second equality, Remark~\ref{rem:LLT Catalanimal def} (i)--(ii)
imply that
$\lambda_{x-1}=\lambda_{x}$,  $|\barRq^{J,J^{c}} | -
|\barRq^{I,I^{c}} | = 1$, and 
$|\barRqt ^{J,J^c}| = |\barRqt ^{I,I^c}|$.
Hence, if
such a pair $x-1,\, x$ exists, then $r_{J} > r_I$.  Thus, it suffices
to prove \eqref{e LLT cuddly induction} for sets $I$ having no such
pair, that is, such that
\begin{equation}\label{eq:reduction}
\text{for each diagonal $D \subseteq \nubold$, $D\cap \nb{I}$ is a
lower order ideal of $(D, \prec)$}.
\end{equation}
We prove this restricted statement---that the claim in \eqref{e LLT
cuddly induction} holds for $I$ satisfying \eqref{eq:reduction}---by
induction on $|I|$.  The base case $I = \varnothing$ is trivial.  Now
suppose $I \ne \varnothing$.  Choose an element $x \in I$ such that
$\nb{x}$ is $\prec$-maximal in $\nb{I}$; this means $\north(\nb{x})$,
$\east(\nb{x}) \notin \nb{I}$.  Set $K = I \setminus \{x\}$.  Let $D =
\diagD{\nb{x}}$ and let $C$, $E$ be the (possibly empty) adjacent
diagonals with adjusted contents $\ctild(\nb{x})-1$,
$\ctild(\nb{x})+1$.  We compute
\begin{align*}
r_{K} - r_I  & = |\barRq ^{K,K^c}| - |\barRq ^{I,I^c}| - |\barRqt ^{K,K^c}|
 + |\barRqt ^{I,I^c}| - \lambda_x\\ 
& = |\barRq ^{K,x}| - |\barRq ^{x,I^c}| - |\barRqt ^{K,x}|
 + |\barRqt ^{x,I^c}| - \lambda_x\\
& = |D \cap \nb{K}| - |D \cap \nb{I^c}| - |C \cap \nb{I}|
 + |E \cap \nb{I^c}| - \lambda_x \\
& = |D \cap \nb{K}| - |C \cap \nb{I}|
 + \chi({\text{$D$ does not contain a row start}}) \\
&\qquad - |D \cap \nb{I^c}| + |E \cap \nb{I^c}|
 - \chi({\text{$D$ does not contain a row end}}),
\end{align*}
where we have used \eqref{ed LLT weight 2} for the last equality.
Since $C\cap \nb{I}$ and $D\cap \nb{K}$ are lower order ideals in $C$
and $D$, and $\north(\nb{x}) \notin \nb{I}$,
\[
|D \cap \nb{K}| - |C \cap \nb{I}| + \chi({\text{$D$ does not contain a row
start}})\ge 0
\]
with equality if and only if $C \cap \nb{I} = \{b \in C \mid b \prec
\north(\nb{x})\}$.  Similarly, 
\[
- |D \cap \nb{I^c}| + |E \cap \nb{I^c}| - \chi({\text{$D$ does not
contain a row end}}) \ge 0
\] 
with equality if and only if $E \cap \nb{I^c} = \{b\in E \mid  \east(\nb{x})
\preceq b \}$.

Thus, we have shown that
\begin{equation}\label{eq:reduction2}
\text{$r_{K} = r_I$ \, if $\nubold \cap \{\south(\nb{x}),
\west(\nb{x})\} \subseteq \nb{I}$, \, and otherwise $r_{K} > r_{I}$.}
\end{equation}
If $\nb{I}$ is a lower order ideal, then so is $\nb{K}$.  In this case
we have $r_{K} = 0$ by induction and $r_{I} = r_{K} = 0$ by
\eqref{eq:reduction2}.  If $\nb{I}$ is not a lower order ideal then
either $\nubold \cap \{\west(\nb{x}), \south(\nb{x})\} \not\subseteq
\nb{I}$ or $\nb{K}$ is not a lower order ideal; by induction we have
$0 \ge r_K > r_I$ if the former holds and $0 > r_K \ge r_I$ if the
latter holds.

\smallskip

\noindent {\em Determining the cub.}  We now prove \eqref{eq:LLT
cuddly} by verifying the two conditions in
Corollary~\ref{cor:evaluating-cubs}.  The condition in part (2) holds
by Lemma~\ref{lem:LLT-1-0-ps}.  For the condition in part (1) we can
assume by induction on the number of boxes in $\nubold$ that
Theorem~\ref{t LLT cuddly} applies to the two Catalanimals $H_I'$ and
$H_I''$ in \eqref{e:cub-coprod-q-style}.  Our remaining task is to
relate the coproduct formula \eqref{e:cub-coprod-q-style} to the
coproduct formula \eqref{e:LLT coprod} for LLT polynomials.  We first
address the coefficients in \eqref{e:cub-coprod-q-style}.  Let $I$ be
a subset of $[l]$ appearing in \eqref{e:cub-coprod-q-style}; by
\eqref{e LLT cuddly induction}, this is equivalent to $\nb{I}$ being a
lower order ideal.  Let $p_I$, $\gamma_I$, $A_I$ denote the magic
number, diagonal lengths, and number of attacking pairs of $\nb{I}$;
let $p_{I^c}$, $\gamma_{I^c}$, $A_{I^c}$ denote the corresponding data
for $\nb{I^c}$.

We begin by computing $|\barRt ^{I,I^c}|$.  Recall from Remark
\ref{rem:LLT Catalanimal def} (i) that $\{ (\nb{i},\nb{j}) \mid
\alpha_{ij} \in \barRt\}$ is the set of attacking pairs in $\nubold $.
Hence, $\{ (\nb{i},\nb{j}) \mid \alpha_{ij} \in \barRt ^{I^c,I}\}$ is
the set of attacking pairs going from $\nb{I^c}$ to $\nb{I}$, which
has size $A(\nb{I^c}, \nb{I})$ in the notation of \eqref{e:LLT
coprod}.  Thus,
\begin{equation}\label{eq:q stat}
|\barRt ^{I,I^c}| = |\barRt | - |\barRt ^{I^c,I}| -
|\barRt ^{I,I}| - |\barRt ^{I^c,I^c}| \\
 = A - A(\nb{I^c}, \nb{I}) - A_I - A_{I^c}.
\end{equation}

Next we compute $|\barRqt ^{I,I^c}|$. Recall that $\barRqt $ is the
set of roots $\alpha_{ij} \in R_+$ with $\nb{i}$ and $\nb{j}$ in
consecutive diagonals.  Then since $\nb{I}$ is a lower order ideal, we
have $\barRqt ^{I,I^c} \subseteq R_+^\prec $ and $\barRqt ^{I^c,I}
\cap R_+^\prec = \varnothing$, with $R_{+}^{\prec }$ as in
Definition~\ref{d:prec order}, allowing us to write
\begin{equation}\label{eq:t statistic0}
|\barRqt ^{I,I^c}| = |\barRqt \cap R_+^\prec| - |\barRqt ^{I,I} \cap
R_+^\prec| - |\barRqt ^{I^c,I^c} \cap R_+^\prec|.
\end{equation}
Using Lemma \ref{lem:magic-num}, this becomes
\begin{equation}\label{eq:t statistic}
|\barRqt ^{I,I^c}| = p + n'(\gamma) - p_I - n'(\gamma_I) - p_{I^c} -
n'(\gamma_{I^c}).
\end{equation}

For the sign in \eqref{e:cub-coprod-q-style}, it remains to compute
$|\barRq ^{I,I^c}|$.  Recall that $\barRq $ is the set of roots
$\alpha_{ij} \in R_+$ with $\nb{i}$ and $\nb{j}$ in the same diagonal,
hence $|\barRq | =
n'(\gamma)$.  Since $\barRq \subseteq R_+^\prec$, we have $\barRq ^{I^c,
I} = \varnothing$, and we conclude that
\begin{align}
\label{eq:sign}
|\barRq ^{I,I^c}| = |\barRq | - |\barRq ^{I,I}|
- |\barRq ^{I^c,I^c}| = n'(\gamma) - n'(\gamma_I)- n'(\gamma_{I^c}).
\end{align}
Combining \eqref{eq:q stat}, \eqref{eq:t statistic}, and
\eqref{eq:sign}, we can express the coefficient in the term for $I$ in
\eqref{e:cub-coprod-q-style} as
\begin{multline}\label{e:q-style-coef}
(-1)^{ |\barRq ^{I,I^c}| + |\barRqt ^{I,I^c}| }
(q\, t)^{-|\barRqt ^{I,I^c}|} q^{-|\barRt ^{I,I^c}|} = \\
q^{A(\nb{I^c}, \nb{I})}\cdot (-1)^{p} (q\, t)^{-p-n'(\gamma)} q^{-A}
\cdot (-1)^{p_I} (q\, t)^{p_I+n'(\gamma_I)} q^{A_I} \cdot
(-1)^{p_{I^c}} (q\, t)^{p_{I^c}+n'(\gamma_{I^c})} q^{A_{I^c}}.
\end{multline}

We next consider the restricted Catalanimals $(H_{\nubold })_I' =
H(R_q|_I, R_t|_I, R_{qt}|_I, \lambda[I]_I)$ and $(H_{\nubold })_I'' =
H(R_q|_{I^c}, R_t|_{I^c}, R_{qt}|_{I^c}, \lambda[I]_{I^c})$ in
\eqref{e:cub-coprod-q-style}.  We claim that $(H_{\nubold })_I' =
H_{\nb{I}}$ and $(H_{\nubold })_I'' = H_{\nb{I^c}}$.  It is clear from
Remark \ref{rem:LLT Catalanimal def} (i) that $R_q|_I, R_t|_I,
R_{qt}|_I$ are the root sets defining the Catalanimal $H_{\nb{I}}$ and
similarly for $I^c$.  It remains to consider the weights.  By
\eqref{ed LLT weight 3} and \eqref{e:lambda-I-simple}, we have
\begin{equation}\label{eq:lambda I comp}
\lambda[I] = - \sum \barRq + \sum (\barRqt \cap R_+^\prec) + \sum
\barRq ^{I,I^c} - \sum \barRqt ^{I,I^c}\,.
\end{equation}
Using the same reasoning that gave \eqref{eq:t statistic0} and
\eqref{eq:sign} to compute $\sum \barRq ^{I,I^c} - \sum \barRq$ and $
\sum (\barRqt \cap R_+^\prec) - \sum \barRqt ^{I,I^c}$
yields
\begin{equation}\label{eq:lambda I comp2}
\lambda[I] =
- \sum \barRq ^{I,I}
+ \sum (\barRqt ^{I,I} \cap R_+^\prec)
- \sum \barRq ^{I^c,I^c}
+ \sum (\barRqt ^{I^c,I^c}\cap R_+^\prec).
\end{equation}
By \eqref{ed LLT weight 3} again, $\lambda[I]_I$ is the weight for
$H_{\nb{I}}$ and $\lambda[I]_{I^c}$ the weight for $H_{\nb{I^c}}$.
Thus, $(H_{\nubold })_I' = H_{\nb{I}}$ and $(H_{\nubold })_I'' =
H_{\nb{I^c}}$, as asserted.

Combining this with \eqref{e:q-style-coef}, the term indexed by $I$ in
\eqref{e:cub-coprod-q-style} becomes
\begin{multline*}
(-1)^{p} (q\, t)^{-p-n'(\gamma)} q^{-A} q^{A(\nb{I^c}, \nb{I})}\\ 
\times \Big( (-1)^{p_I} (q\, t)^{p_I+n'(\gamma_I)} q^{A_I}
\cub(H_{\nb{I}})(X)\Big) \Big( (-1)^{p_{I^c}}
(q\, t)^{p_{I^c}+n'(\gamma_{I^c})} q^{A_{I^c}}
\cub(H_{\nb{I^{c}}})(Y)\Big).
\end{multline*}
For $0 < |I| < l$, this is equal to $(-1)^{p} (q\, t)^{-p-n'(\gamma)}
q^{-A} q^{A(\nb{I^c}, \nb{I})} \Gcal_{\nb{I}}(X;q)
\Gcal_{\nb{I^c}}(Y;q)$ by induction.  Hence, by \eqref{e:LLT coprod},
the symmetric function $f = (-1)^{p} (q\, t)^{-p-n'(\gamma)} q^{-A}
\Gcal_{\nubold }(X;q)$ satisfies the condition in Corollary
\ref{cor:evaluating-cubs} part (1).
\end{proof}

\section{LLT Catalanimals}
\label{s mn LLT Catalanimal}

We now generalize the results of the previous section by constructing,
for any tuple of skew shapes $\nubold $, an $(m,n)$-cuddly Catalanimal
whose cub is a scalar multiple of the LLT polynomial $\Gcal_{\nubold
}$.  As a corollary, we obtain a raising operator formula for
$\nabla^m \Gcal_{\nubold }$.

\subsection{\texorpdfstring{$m$}{m}-stretching}

Let $\theta$ be a skew shape and $m \in \ZZ_+, o \in \ZZ$.  
We construct a new skew shape
$\theta^{m,o}$ by stretching $\theta$ vertically by a factor of $m$ as
follows: for each box $x$ of content $c$ in $\theta$, place $m$ boxes
of contents $o+mc,o+mc-1, \dots, o+mc-m+1$ in the same column as $x$.
The set of $m$ boxes arising from $x$ in this way is called a
\emph{stretched box}, denoted $\stre(x) \subseteq \theta^{m,o}$.  Thus,
the stretched boxes
partition the boxes of $\theta^{m,o}$ into
$|\theta|$ sets of size $m$.  For example, for $m=3$ and $o =-4$,
\begin{align*}
\theta = (3,2)/(0,0) = 
\raisebox{-2.7mm}{
\begin{tikzpicture}[scale = .39]
\begin{scope}
\draw[thick] (0,0) grid (2,2);
\draw[thick] (0,0) grid (3,1);
\end{scope}
\end{tikzpicture}}
\qquad \text{yields} \qquad \theta^{m,o}
= \theta ^{3,-4} =
\raisebox{-1.12cm}{
\begin{tikzpicture}[scale = .34]
\draw[draw = none, fill = blue!33] (2.27,0.31) rectangle (3-0.27,3-0.31);
\draw[draw = none, fill = blue!33] (1.27,2.31) rectangle (2-0.27,5-0.31);
\draw[draw = none, fill = blue!33] (1.27,5.31) rectangle (2-0.27,8-0.31);
\draw[draw = none, fill = blue!33] (0.27,4.31) rectangle (1-0.27,7-0.31);
\draw[draw = none, fill = blue!33] (0.27,7.31) rectangle (1-0.27,10-0.31);
\draw[thick] (0,4) grid (1,10);
\draw[thick] (1,2) grid (2,8);
\draw[thick] (2,0) grid (3,3);
\draw[thin, black!60] (0,0) -- (0,10);
\draw[thin, black!60] (0,0) -- (3,0);
\end{tikzpicture}}
\end{align*}
where the shaded rectangles indicate the stretched boxes.

\begin{defn}
Given a tuple of skew shapes $\nubold = (\nu_{(1)}, \dots,
\nu_{(k)})$, $m \in \ZZ_+$, and $\obold = (o_1,\dots, o_k) \in
\ZZ^k$ satisfying
\begin{align}
\label{eq:offset assumption}
 o_1 \le o_2 \le \cdots \le o_k < m + o_1,
\end{align}
the associated \emph{$m$-stretching of $\nubold $} is $\nubold ^m =
\nubold (m,\obold ) = (\nu_{(1)}^{m,o_1}, \ldots,
\nu_{(k)}^{m,o_k}).$ We often use the abbreviation $\nubold ^m$ even
though
it actually
depends on $\obold $
as well.
\end{defn}

We use the same notation $\nbm{i}$, $\nbm{I}$ as in \eqref{e:nu(I)}
for the boxes of $\nubold ^{m}$
numbered in reading order, and for the set of boxes corresponding to a
set of indices $I\subseteq [l]$, where $l = dm = |\nubold ^m|$ if $d =
|\nubold |$.

The assumption \eqref{eq:offset assumption} on the offsets $\obold $
allows us to relate attacking pairs in $\nubold $ to attacking pairs
in $\nubold ^m$, as follows.  Let $\Abold (\nubold )$ denote the set
of attacking pairs in $\nubold $ (with the pair in increasing reading
order, as always).  For $(\nb{i},\nb{j}) \in \Abold (\nubold )$, with
$\nb{i}\in \nu _{(r)}$ and $\nb{j}\in \nu _{(s)}$, we set
\begin{equation}\label{eq:aij def}
\begin{aligned}
a_{ij} & \defeq \big|\big\{(a,b) \in \Abold (\nubold ^m) \mid a \in
\stre(\nb{i}), b \in \stre(\nb{j}) \big\}\big| \\
& =
\begin{cases}
  m+ o_r - o_s & \text{if  $r < s$} \\
  1+ o_r - o_s & \text{if  $r > s$}.
\end{cases}
\end{aligned}
\end{equation}
Note that \eqref{eq:offset assumption} implies $a_{ij} \ge 1$, and
furthermore,
\begin{equation}\label{eq:aij 2}
a_{ij}-1 = \big|\big\{(a,b) \in \Abold (\nubold ^m) \mid a \in
\stre(\nb{j}), b \in \stre(\nb{i}) \big\}\big|,
\end{equation}
as illustrated in
Example \ref{ex:attack stretch}.
Finally, for any
$1 \le i < j \le d$,
\begin{equation}\label{eq:aij 3}
(\nb{i},\nb{j}) \notin \Abold (\nubold ) 
\Rightarrow \big\{(a,b) \in \Abold (\nubold ^m) \mid a,b \in
\stre(\nb{i}) \cup \stre(\nb{j}) \big\} = \varnothing.
\end{equation}

\begin{example}\label{ex:attack stretch}
Let $\nubold = ((32)/(00),(2)/(0))$, $m=3$, and $\obold = (-4,-2)$.
The tuples $\nubold$ and  $\nubold^m$ are shown below with boxes numbered
in reading order, along with another
drawing of $\nubold^m$ to illustrate additional features.
There are five attacking pairs in $\nubold $, with $a_{24} =a_{34}
=a_{56} = 1$, $a_{45}=a_{67} = 3$.  The three attacking pairs in
$\nubold ^{m}$ counted by $a_{45}$ are indicated by the solid arrows
and the two attacking pairs counted by $a_{45}-1$ as in
\eqref{eq:aij 2} are indicated by the dashed arrows.
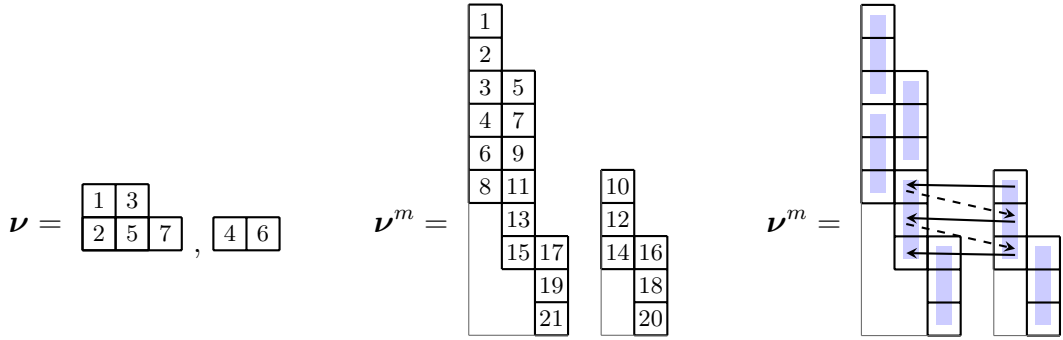
\[
\nubold  = 
\raisebox{-2.7mm}{
\begin{tikzpicture}[scale = .44]
\begin{scope}
\draw[thick] (0,0) grid (2,2);
\draw[thick] (0,0) grid (3,1);
\node at (0.5, 1.5) {\scriptsize 1};
\node at (1.5, 1.5) {\scriptsize 3};
\node at (0.5, 0.5) {\scriptsize 2};
\node at (1.5, 0.5) {\scriptsize 5};
\node at (2.5, 0.5) {\scriptsize 7};
\end{scope}
\end{tikzpicture}
,
\begin{tikzpicture}[scale = .44]
\begin{scope}
\draw[thick] (0,0) grid (2,1);
\node at (0.5, 0.5) {\scriptsize 4};
\node at (1.5, 0.5) {\scriptsize 6};
\end{scope}
\end{tikzpicture}}
\qquad \quad
\nubold^m = 
\raisebox{-1.4cm}{
\begin{tikzpicture}[scale = .44]
\tikzstyle{vertex}=[inner sep=0pt, outer sep=2.7pt]
\begin{scope}
\draw[thick] (0,4) grid (1,10);
\draw[thick] (1,2) grid (2,8);
\draw[thick] (2,0) grid (3,3);

\draw[thick] (0+4,2) grid (1+4,5);
\draw[thick] (1+4,0) grid (2+4,3);

\draw[thin, black!60] (0,0) -- (0,10);
\draw[thin, black!60] (0,0) -- (3,0);
\draw[thin, black!60] (0+4,0) -- (0+4,5);
\draw[thin, black!60] (0+4,0) -- (2+4,0);

\foreach \x/\y/\nn in 
{1/10/1, 1/9/2, 1/8/3, 1/7/4, 1/6/6, 1/5/8,
2/8/5, 2/7/7, 2/6/9, 2/5/11, 2/4/13, 2/3/15,
3/3/17, 3/2/19, 3/1/21}
{\node at (\x-0.5,\y-0.5) {\scriptsize \nn};}
\foreach \x/\y/\nn in 
{1/5/10, 1/4/12, 1/3/14,
2/3/16, 2/2/18, 2/1/20}
{\node at (\x-0.5+4,\y-0.5) {\scriptsize \nn};}
\end{scope}
\end{tikzpicture}}
\qquad \quad
\nubold^m = 
\raisebox{-1.4cm}{
\begin{tikzpicture}[scale = .44]
\tikzstyle{vertex}=[inner sep=0pt, outer sep=2.7pt]
\begin{scope}
\draw[draw = none, fill = blue!20] (2.27,0.31) rectangle (3-0.27,3-0.31);
\draw[draw = none, fill = blue!20] (1.27,2.31) rectangle (2-0.27,5-0.31);
\draw[draw = none, fill = blue!20] (1.27,5.31) rectangle (2-0.27,8-0.31);
\draw[draw = none, fill = blue!20] (0.27,4.31) rectangle (1-0.27,7-0.31);
\draw[draw = none, fill = blue!20] (0.27,7.31) rectangle (1-0.27,10-0.31);
\draw[draw = none, fill = blue!20] (0.27+4,2.31) rectangle (1-0.27+4,5-0.31);
\draw[draw = none, fill = blue!20] (1.27+4,0.31) rectangle (2-0.27+4,3-0.31);

\draw[thick] (0,4) grid (1,10);
\draw[thick] (1,2) grid (2,8);
\draw[thick] (2,0) grid (3,3);

\draw[thick] (0+4,2) grid (1+4,5);
\draw[thick] (1+4,0) grid (2+4,3);

\draw[thin, black!60] (0,0) -- (0,10);
\draw[thin, black!60] (0,0) -- (3,0);
\draw[thin, black!60] (0+4,0) -- (0+4,5);
\draw[thin, black!60] (0+4,0) -- (2+4,0);

\node[vertex] (L2) at (1+.1,2+.5) {};
\node[vertex] (L3u) at (1+.1,3+.5+.067) {  };
\node[vertex] (L4u) at (1+.1,4+.5+.067) {  };
\node[vertex] (L3d) at (1+.1,3+.5-.067) {  };
\node[vertex] (L4d) at (1+.1,4+.5-.067) {  };

\node[vertex] (R2u) at (4+.9,2+.5+.067) { };
\node[vertex] (R3u) at (4+.9,3+.5+.067) { };
\node[vertex] (R2d) at (4+.9,2+.5-.067) { };
\node[vertex] (R3d) at (4+.9,3+.5-.067) { };
\node[vertex] (R4) at (4+.9,4+.5) { };

\draw[thick, dashed,-stealth] (L4d) to (R3u);
\draw[thick, dashed,-stealth] (L3d) to (R2u);
\draw[thick,-stealth] (R2d) to (L2);
\draw[thick,-stealth] (R3d) to (L3u);
\draw[thick,-stealth] (R4) to (L4u);
\end{scope}
\end{tikzpicture}}
\]
\end{example}

\subsection{Definition of LLT Catalanimals}
\label{ss:mn LLT Catalanimal}

For $(m,n) \in \ZZ_+ \times \ZZ$ (not necessarily coprime), we define
a vector of integers $\bbold (m, n)$ of length $m$ by
\begin{equation}\label{eq:def south steps}
\bbold (m, n)_{i} = \lceil i n/m \rceil -\lceil (i-1) n/m \rceil
\qquad (i=1,\ldots,m).
\end{equation}
More pictorially, if $n\geq 0$, then $\bbold (m,n)_{i}$ is the number
of south steps on the line $x=i-1$ in the highest south/east lattice
path weakly below the line segment from $(0,n)$ to $(m,0)$.  Note that
$\bbold (dm,dn)$ is the concatenation of $d$ copies of $\bbold
(m,n)$.

\begin{defn}\label{def:mn LLT Catalanimal}
Let $(m,n)\in \ZZ _{+}\times \ZZ$ be a pair of coprime
integers, let $\nubold $ be a tuple of skew shapes with $d =|\nubold
|$ boxes, and let $\nubold ^m = \nubold (m,\obold )$ be an
$m$-stretching of $\nubold $.  Let $l = dm = |\nubold ^{m}|$ be the
number of boxes in $\nubold ^{m}$.

We define the \emph{LLT Catalanimal} $H^{m,n}_{\nubold ^m} =
H(R_q,R_t, R_{qt}, \lambda)$ as follows.  The root sets $R_{q}$,
$R_{t}$, $R_{qt}$ are defined as in \eqref{ed Rq}--\eqref{ed Rqt} but
with $\nubold ^{m}$ in place of $\nubold $.  The weight $\lambda$ is
defined by
\begin{equation}\label{e:lambda-vs-hat-lambda}
\lambda =\hat{\lambda} + \bboldtild,
\end{equation}
where $\hat{\lambda}$ is the weight given by \eqref{ed LLT weight}
with $\nubold ^{m}$ in place of $\nubold $, and $\bboldtild \in \ZZ^l$
is given by
\begin{equation}\label{e:btilde}
\bboldtild_i = \bbold (m,n)_{\operatorname{mod}_m(c-o_s)}, \quad \text{where
$\nbm{i}$ is a box of content $c$ in $\nu _{(s)}^{m,o_{s}}$.}
\end{equation}
Here $\operatorname{mod}_m(c-o_s)$ denotes the 
integer $j\in [m]$ such that $c-o_{s}\equiv j\pmod{m}$, 
and $\bbold(m,n)$ is defined above.

Viewed as a filling of $\nubold ^m$, as in Remark~\ref{rem:LLT
Catalanimal def} (iii), the weight $\lambda$ is obtained from
$\hat{\lambda}$ by adding the vector $\bbold (m,n)$ to each stretched
box from north to south as in Example \ref{ex:hat zeta plus b} (i).
\end{defn}

\begin{example}\label{ex:hat zeta plus b}
(i) Let $\nubold = (\, \partition{~&~ \\ ~&~ &~}\, ,\,\,
\partition{\\~&~}\, )$, $m = 3$, $\obold = (-4,-2)$ be as in Example
\ref{ex:attack stretch} and $n=2$.  We have $\bbold (m,n) = (1,1,0)$.
The weights $\hat{\lambda}, \bboldtild, \lambda$ are drawn below as
fillings of $\nubold (m, \obold )$.  For the weight $\lambda$, entries
of the second shape are shown in bold to help translate between its
filling and vector depictions.
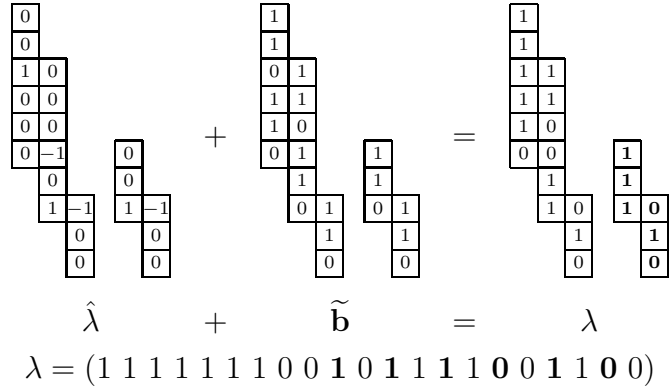
\begin{gather*}
\begin{array}{ccccc}
{\TTiny \text{ $\mytableau{ 0 \\0 \\ 1&0 \\ 0&0 \\0&0 \\0&-1\\ &0& \\ & 1&-1\\ & &0 \\ & &0     }  \ \ \ \mytableau{ \\ \\ \\ \\ \\ 0\\ 0 \\ 1&-1\\ & 0 \\ & 0 }$  }}
& + &
{\TTiny \text{ $\mytableau{ 1 \\1 \\ 0&1 \\ 1&1 \\1&0 \\0&1 \\ &1& \\ & 0&1 \\ & &1 \\ & &0     }  \ \ \ \mytableau{ \\ \\ \\ \\ \\ 1\\ 1 \\ 0&1 \\ & 1 \\ & 0 }$  }}
& = &
{\TTiny \text{ $\mytableau{ 1 \\1 \\ 1&1 \\ 1&1 \\1&0 \\0&0 \\ &1& \\ & 1&0 \\ & &1 \\ & &0     }  \ \ \ \mytableau{ \\ \\ \\ \\ \\ \textbf{1}\\ \textbf{1} \\ \textbf{1}&\textbf{0} \\ & \textbf{1} \\ & \textbf{0} }$  }}
\\
\\[-1ex]
\hat{\lambda} & + & \bboldtild & = & \lambda 
\end{array}\\
\lambda = (1\ 1\ 1\ 1\ 1\ 1\ 1\ 0\ 0\ \textbf{1} \ 0\ \textbf{1} \ 1\ \textbf{1} \ 1\ \textbf{0} \ 0\ \textbf{1} \ 1\ \textbf{0} \ 0)
\end{gather*}

(ii) The LLT Catalanimal $H^{3,2}_{\nubold(3,\obold)}$ for 
$\nubold  = ((2),(1))$, $\obold  = (-2,-2)$ is shown in
Figure \ref{fig:mn LLT Catalanimal}.
\end{example}

\begin{figure}
\[
\begin{tikzpicture}[scale = .55]
\draw[draw = none, fill = black!14] (4,-1) rectangle (5,-2);
 \draw[draw = none, fill = black!14] (5,-1) rectangle (6,-2);
 \draw[draw = none, fill = black!14] (5,-2) rectangle (6,-3);
 \draw[draw = none, fill = black!14] (6,-1) rectangle (7,-2);
 \draw[draw = none, fill = black!14] (6,-2) rectangle (7,-3);
 \draw[draw = none, fill = black!14] (6,-3) rectangle (7,-4);
 \draw[draw = none, fill = black!14] (7,-1) rectangle (8,-2);
 \draw[draw = none, fill = black!14] (7,-2) rectangle (8,-3);
 \draw[draw = none, fill = black!14] (7,-3) rectangle (8,-4);
 \draw[draw = none, fill = black!14] (7,-4) rectangle (8,-5);
 \draw[draw = none, fill = black!14] (8,-1) rectangle (9,-2);
 \draw[draw = none, fill = black!14] (8,-2) rectangle (9,-3);
 \draw[draw = none, fill = black!14] (8,-3) rectangle (9,-4);
 \draw[draw = none, fill = black!14] (8,-4) rectangle (9,-5);
 \draw[draw = none, fill = black!14] (8,-5) rectangle (9,-6);
 \draw[draw = none, fill = black!14] (8,-6) rectangle (9,-7);
 \draw[draw = none, fill = black!14] (9,-1) rectangle (10,-2);
 \draw[draw = none, fill = black!14] (9,-2) rectangle (10,-3);
 \draw[draw = none, fill = black!14] (9,-3) rectangle (10,-4);
 \draw[draw = none, fill = black!14] (9,-4) rectangle (10,-5);
 \draw[draw = none, fill = black!14] (9,-5) rectangle (10,-6);
 \draw[draw = none, fill = black!14] (9,-6) rectangle (10,-7);
 \draw[draw = none, fill = black!14] (9,-7) rectangle (10,-8);
 \draw[draw = none, fill = gray!100] (3+0.5, -1-0.5) circle (.2);
\draw[draw = none, fill = gray!100] (4+0.5, -2-0.5) circle (.2);
\draw[draw = none, fill = gray!100] (5+0.5, -3-0.5) circle (.2);
\draw[draw = none, fill = gray!100] (6+0.5, -4-0.5) circle (.2);
\draw[draw = none, fill = gray!100] (7+0.5, -5-0.5) circle (.2);
\draw[draw = none, fill = gray!100] (8+0.5, -7-0.5) circle (.2);
\draw[draw = none, fill = gray!100] (9+0.5, -8-0.5) circle (.2);
\draw[draw = none, fill = \mymidgray] (2,-1) rectangle (3,-2);
 \draw[draw = none, fill = \mymidgray] (3,-2) rectangle (4,-3);
 \draw[draw = none, fill = \mymidgray] (4,-3) rectangle (5,-4);
 \draw[draw = none, fill = \mymidgray] (5,-4) rectangle (6,-5);
 \draw[draw = none, fill = \mymidgray] (6,-5) rectangle (7,-6);
 \draw[draw = none, fill = \mymidgray] (7,-6) rectangle (8,-7);
 \draw[thin, black!31] (1,-1) -- (10,-1);
\draw[thin, black!31] (2,-1) -- (2,-1);
\draw[thin, black!31] (2,-2) -- (10,-2);
\draw[thin, black!31] (3,-2) -- (3,-1);
\draw[thin, black!31] (3,-3) -- (10,-3);
\draw[thin, black!31] (4,-3) -- (4,-1);
\draw[thin, black!31] (4,-4) -- (10,-4);
\draw[thin, black!31] (5,-4) -- (5,-1);
\draw[thin, black!31] (5,-5) -- (10,-5);
\draw[thin, black!31] (6,-5) -- (6,-1);
\draw[thin, black!31] (6,-6) -- (10,-6);
\draw[thin, black!31] (7,-6) -- (7,-1);
\draw[thin, black!31] (7,-7) -- (10,-7);
\draw[thin, black!31] (8,-7) -- (8,-1);
\draw[thin, black!31] (8,-8) -- (10,-8);
\draw[thin, black!31] (9,-8) -- (9,-1);
\draw[thin, black!31] (9,-9) -- (10,-9);
\draw[thin, black!31] (10,-9) -- (10,-1);
\draw[thin] (1,-1) -- (2,-1);
\draw[thin] (2,-1) -- (2,-2);
\draw[thin] (2,-2) -- (3,-2);
\draw[thin] (3,-2) -- (3,-3);
\draw[thin] (3,-3) -- (4,-3);
\draw[thin] (4,-3) -- (4,-4);
\draw[thin] (4,-4) -- (5,-4);
\draw[thin] (5,-4) -- (5,-5);
\draw[thin] (5,-5) -- (6,-5);
\draw[thin] (6,-5) -- (6,-6);
\draw[thin] (6,-6) -- (7,-6);
\draw[thin] (7,-6) -- (7,-7);
\draw[thin] (7,-7) -- (8,-7);
\draw[thin] (8,-7) -- (8,-8);
\draw[thin] (8,-8) -- (9,-8);
\draw[thin] (9,-8) -- (9,-9);
\draw[thin] (9,-9) -- (10,-9);
\draw[thin] (10,-9) -- (10,-10);
\node at (3/2,-3/2) {\tiny $1 $};
\node at (5/2,-5/2) {\tiny $1 $};
\node at (7/2,-7/2) {\tiny $1 $};
\node at (9/2,-9/2) {\tiny $1 $};
\node at (11/2,-11/2) {\tiny $1 $};
\node at (13/2,-13/2) {\tiny $0 $};
\node at (15/2,-15/2) {\tiny $0 $};
\node (vv) at (17/2,-17/2) {\tiny $1 $};
\node at (19/2,-19/2) {\tiny $0 $};
\node[anchor=west] at (-9.6,-10.5+3*1.1) {$p(\nubold ^m) = 1$};
\node[anchor=west] at (-9.6,-10.5+2*1.1) {$\gamma(\nubold ^m) = (1, 1, 1, 1, 1, 1, 1, 1, 1)$};
\node[anchor=west] at (-9.6,-10.5+1*1.1) {$n'(\gamma(\nubold ^m)) =  0$};
\node[anchor=west] at (-9.6,-10.5) {$\Asum = \sum_{(i,j) \in \Abold (\nubold )} a_{ij} = 4$};
\node[left=12.6mm of vv] {$H^{m,n}_{\nubold^m}$};
\node at (-5.5,-2.5) {$ {\fontsize{8pt}{6pt}\selectfont \tableau{
1\\1\\1&0\\ \bl &1\\ \bl &0\\} } \quad
{\fontsize{8pt}{6pt}\selectfont \tableau{
1\\1\\0\\ \bl \\ \bl \\} }  $};
\node at (-5.5,-5.6) {$\lambda$, as a filling of $\nubold ^m$};
\end{tikzpicture}
\qquad \qquad 
\raisebox{23mm}{
\begin{tikzpicture}[scale = .39]
\begin{scope}
\draw[draw = none, fill = black!100] (2,-1) rectangle (3,-2);
\node[anchor = west] at (3,-1.5) {\scriptsize  \  $R_+ \setminus R_q$};
\end{scope}
\begin{scope}[yshift = -34*1]
\draw[draw = none, fill = \mymidgray] (2,-1) rectangle (3,-2);
\node[anchor = west] at (3,-1.5) {\scriptsize  \  $R_q \setminus R_t$};
\end{scope}
\begin{scope}[yshift = -34*2]
\draw[thin, black!0] (2,-1) -- (3,-1);
\draw[thin, black!0] (2,-1) -- (2,-2);
\draw[thin, black!0] (2,-2) -- (3,-2);
\draw[thin, black!0] (3,-2) -- (3,-1);
\draw[draw = none, fill = gray!100] (2+0.5, -1-0.5) circle (.2);
\node[anchor = west] at (3,-1.5) {\scriptsize  \  $R_t \setminus R_{qt}$};
\end{scope}
\begin{scope}[yshift = -34*3]
\draw[draw = none, fill = black!14] (2,-1) rectangle (3,-2);
\node[anchor = west] at (3,-1.5) {\scriptsize \  $R_{qt}$};
\end{scope}
\end{tikzpicture}}
\]
\caption{\label{fig:mn LLT Catalanimal}%
The LLT Catalanimal $H^{m,n}_{\nubold ^m}$ for $m =3$, $n=2$,
$\nubold  = ((2),(1))$, $\obold  = (-2,-2)$, and its associated
statistics. By Theorem~\ref{t mn LLT cuddly}, $\psi _{\hGamma }(-q^5t
H^{m,n}_{\nubold ^m}) =\Gcal_{\nubold }[-MX^{m,n}]$.}
\end{figure}

\begin{remark}\label{rmk:product of Catalanimals}
Define a binary operation  $\uplus$ on sets of roots as follows:
for  $A \subseteq R_+(\GL_l)$ and $B \subseteq R_+(\GL_{l'})$,
\begin{equation}\label{e:root-join}
A\uplus B \defeq A \sqcup \big\{ (i+l,j+l) \mid (i,j)\in B \big\}
\sqcup \big\{ (i,j) \mid
1\leq i\leq l<j\leq l+l'
\big\} \subseteq R_+(\GL_{l+l'}).
\end{equation}
The product of two Catalanimals $H = H(R_q, R_t, R_{qt}, \lambda)$ and
$H' = H(R_q', R_t', R_{qt}', \lambda')$ in the concrete shuffle
algebra $\Scal _\hGamma$ is another Catalanimal,
\begin{equation}\label{eq:Catalanimal prod}
H H' = H(R_q \uplus R_q', R_t \uplus R_t', R_{qt} \uplus R_{qt}',
(\lambda ;\lambda')),
\end{equation}
where $(\lambda ;\lambda')$ denotes the concatenation of $\lambda$ and
$\lambda'$.

The definition of the LLT Catalanimals $H^{m,n}_{\nubold (m, \obold
)}$ interacts well with this product in the following sense.  If
$\nubold $ decomposes as $\nubold = \nubold ' \sqcup \nubold ''$
(meaning that $\nubold _{(i)} = \nubold _{(i)}' \sqcup \nubold _{(i)}
''$ for each $i$), and $\ctild(a)+1 < \ctild(b)$ for all boxes $a \in
\nubold '$, $b \in \nubold ''$, then the root sets and weights of
$H^{m,n}_{\nubold (m,\obold )}$ are constructed from those of
$H^{m,n}_{\nubold '(m,\obold )}, H^{m,n}_{\nubold ''(m,\obold )}$ as
in \eqref{eq:Catalanimal prod}, giving $H^{m,n}_{\nubold (m, \obold )}
= H^{m,n}_{\nubold '(m,\obold )} H^{m,n}_{\nubold ''(m,\obold )}$ in
$\Scal _\hGamma$.
\end{remark}

\subsection{Determining the cubs}
\label{ss:trap-the-cub}

We now come to our main theorem giving a Catalanimal formula for LLT
polynomials in any of the subalgebras $\Lambda (X^{m,n})$ of $\Ecal
^{+}$.

\begin{thm}\label{t mn LLT cuddly}
For any tuple of skew shapes $\nubold $ and any $m$-stretching
$\nubold ^m = \nubold (m,\obold )$, the LLT Catalanimal
$H^{m,n}_{\nubold ^m} = H(R_q,R_t,R_{qt},\lambda)$ is $(m,n)$-cuddly
with cub given by
\begin{equation}\label{eq:mn LLT cuddly}
\cub(H^{m,n}_{\nubold ^m}) = (-1)^{p}(q\, t)^{-p-n'(\gamma)}
q^{-\Asum} \Gcal_{\nubold }(X;q),
\end{equation}
where $p = p(\nubold ^m)$, $\gamma = \gamma(\nubold ^m)$ are the magic
number and diagonal lengths of $\nubold ^m$, and $\Asum = \sum_{(i,j)
\in \Abold (\nubold )} a_{ij}$, with $\Abold (\nubold )$ and $a_{ij}$
as in \eqref{eq:aij def}.
\end{thm}

\begin{remark}\label{rem:intro-form}
Equation \eqref{eq:mn LLT cuddly} can also be written as the following
more precise form of the formula \eqref{e:main-preview} mentioned in
the introduction:
\begin{equation}\label{e:intro-form}
\psi _{\hGamma}(H^{m,n}_{\nubold ^m}(\zz )) = (-1)^{p}(q\,
t)^{-p-n'(\gamma)} q^{-\Asum} \Gcal_{\nubold }[-MX^{m,n}],
\end{equation}
where $\psi _{\hGamma}$ is the isomorphism from the shuffle algebra to
the Schiffmann algebra defined in \S \ref{ss:shuffle-iso}.
\end{remark}

The proof will be given below after some further remarks and
preliminary lemmas.

\begin{remark}
(i) One can check that in fact $p(\nubold ^m) = p(\nubold )$ and
$n'(\gamma(\nubold ^m)) = m \cdot n'(\gamma(\nubold ))$.

(ii) For constant offsets $\obold = (c,c,\dots, c)$, $\Asum =
\sum_{(i,j) \in \Abold (\nubold )} a_{ij}$ is $m$ times the number of
attacking pairs in $\nubold $ in which both boxes have the same
content, plus the number of attacking pairs in which the boxes have
different contents.

(iii) If $\nubold $ is a single skew shape $(\theta )$, the only
effect of the offset is to translate the $m$-stretched diagram
$\nubold ^{m} = (\theta ^{m,o})$ vertically.  In this case there are
no attacking pairs, so $\Acal = 0$, and the Catalanimal
$H^{m,n}_{\nubold ^{m}}$, whose cub is $(-1)^{p}(q\, t)^{-p-n'(\gamma
)} s_{\theta }(X)$, does not depend on the offset.
\end{remark}

\begin{remark}\label{rmk:negut vs Catalanimal def}
If $\nubold =(\theta )$, where $\theta $ is a ribbon skew shape of
size $d$, then $\theta ^{m,o}$ is a ribbon of size $l=dm$, and the
function $\phi (\zz ) = \phi (R_{q},R_{t},R_{qt},\lambda )$ in
\eqref{eq:phi for Catalanimal} such that $H^{m,n}_{\nubold ^{m}} =
\sigma _{\hGamma }(\phi (\zz ))$ is given by
\begin{equation}\label{e:phi-for-ribbon}
\phi (\zz ) = \frac{\zz^{\lambda } }{\prod _{i=1}^{l-1} (1 - q\, t\,
z_{i}/z_{i+1})},
\end{equation}
where $\lambda = \bbold (dm,dn)+\sum _{i\in I} \alpha _{mi,mi+1}$,
with $I\subseteq [d-1]$ the set of indices such that boxes $\nb{i}$
and $\nb{i+1}$ of $\theta $ are in the same row.

Negut \cite[Proposition~6.1]{Negut14} showed that image under 
$\sigma_\hGamma$ of the rational
function in \eqref{e:phi-for-ribbon} 
lies in the shuffle algebra $\Scal _\hGamma$
 for any weight $\lambda $. 
For the specific
weight $\lambda $ occurring here, translating
\cite[Proposition~6.7]{Negut14} into our notation using
\eqref{e:Negut-diagram} gives $\psi (\phi (\zz )) = (-q t)^{-|I|}
s_{\theta }[-MX^{m,n}]$, which agrees with \eqref{eq:mn LLT cuddly},
because for ribbons we have $n'(\gamma ) = 0$, and the magic number
$p$ is equal to $|I|$.  This result of Negut is thus the special case
of Theorem~\ref{t mn LLT cuddly} for $\nubold $ a single ribbon skew
shape.
\end{remark}

\begin{lemma}\label{lem:LLT 1-q}
If each component of $\nubold$ is a disjoint union of ribbon skew
shapes, and $\nubold $ has no attacking pairs, then $\Gcal _{\nubold
}(X;q)$ is a product of ribbon skew Schur functions, with no
dependence on $q$.  Otherwise, $\omega\Gcal_{\nubold } [1-q] = 0$.
\end{lemma}

\begin{proof}
The first statement is clear from the definition of $\Gcal _{\nubold
}$.  For the second, note that $\omega\Gcal_{\nubold } [1-q] =
(-1)^{|{\nubold}|}\Gcal_{\nubold }[q-1]$, and that $\Gcal_{\nubold
}[q-1]$ is given by evaluating $\omega _{Y}\Gcal_{\nubold} [X+Y]$ at
$X = q$ and $Y = y_{1}$, followed by setting $y_{1} = -1$.  It then
follows from Lemma \ref{lem:super-G} that $\Gcal_{\nubold }[q-1] =
\sum_{T \in \Tcal } q^{\inv(T)} q^{\#1 \text{\,'s}}
(-1)^{\#\overline{1} \text{\,'s}}$, where the sum is over the set
$\Tcal $ of super tableaux on ${\nubold }$ in the alphabet $1
<\overline{1}$.

If $\nubold $ is not a disjoint union of ribbon skew shapes,
the set $\Tcal $ is empty and we are done.

It remains to show that $\Gcal_{\nubold }[q-1] = 0$ if $\nubold $ has
an attacking pair.  If so, adapting the argument in the proof of
\cite[Lemma~5.1]{HaHaLo05}, we construct an involution $\Psi$ on
$\Tcal $ that cancels the terms $q^{\inv(T)} q^{\#1 \text{\,'s}}
(-1)^{\#\overline{1} \text{\,'s}}$.  To define $\Psi $, let $b$ be the
last box in reading order that is part of an attacking pair and let
$a$ be the last box in reading order such that $(a,b)$ is an attacking
pair; then we let $\Psi \, T$ be the tableau obtained from $T$ by
changing the sign of $T(a)$.  It is not hard to see that $a$ is
necessarily the southeast corner of the ribbon containing it, and
therefore that $\Psi \, T$ is indeed a super tableau on $\nubold $.
Since $a$ and $b$ only depend on $\nubold $, it is clear that $\Psi\,
\Psi \, T = T$.  For $T$ with $T(a)=1$, changing this $1$ to a
$\overline{1}$ adds one to the number of inversions and subtracts one
from the number of 1's, hence the contributions of $T$ and $\Psi T$ to
$\Gcal_{\nubold }[q-1]$ cancel.
\end{proof}

\begin{lemma}\label{lem:a-two-ways}
Given coprime integers $(m,n)\in \ZZ _{+}\times \ZZ$ and an integer
$d>0$, the vector $\bbold (dm,dn)$ satisfies
\begin{equation}\label{e:a-two-ways}
\sum _{i=1}^{dm} (i-1)\bbold (dm,dn)_{i} = a \defeq \frac{1}{2}d (d m
n - m - n +1),
\end{equation}
where $a$ is the exponent in Theorem~\ref{thm:principal-spec} and
Corollary~\ref{cor:evaluating-cubs}.
\end{lemma}

\begin{proof}
Adding $rm$ to $n$ adds a constant vector $(r,r,\ldots,r)$ to $\bbold
(dm,dn)$ and thus increases both sides of \eqref{e:a-two-ways} by $r
\binom{dm}{2}$.  It therefore suffices to prove \eqref{e:a-two-ways}
for $n\geq 0$.

The left hand side of \eqref{e:a-two-ways} is then the area under the
highest lattice path weakly below the line segment from $(0,dn)$ to
$(dm,0)$, or the number of complete lattice squares below the diagonal
in a $d n \times d m$ rectangle.  Call this number $b$.  Then
$d^{2}mn-2b$ is the number of lattice squares cut by the diagonal.

If $d=1$, the cut squares form a ribbon of size $m+n-1$. In general,
they are a union of $d$ copies of this ribbon. Hence,
$b = \frac{1}{2}(d^2mn-d(m+n-1)) = a$.
\end{proof}

\begin{lemma}\label{lem:LLT ps}
The element $\phi(\zz )= \phi(R_q,R_t,R_{qt}, \lambda)$ defined in
\eqref{eq:phi for Catalanimal} corresponding to the LLT Catalanimal
$H^{m,n}_{\nubold ^m} = H(R_q,R_t,R_{qt}, \lambda)$ has principal
specialization
\begin{equation}\label{eq:LLT ps}
\phi(1,t, \ldots, t^{l -1}) = t^a (\omega f)[1-q]/(1-q)^l,
\end{equation}
where $f = (-1)^{p}(q\, t)^{-p-n'(\gamma)} q^{-\Asum} \Gcal_{\nubold
}(X;q)$ is the right hand side of \eqref{eq:mn LLT cuddly}, $l =
|\nubold ^m|$, and $a = \frac{1}{2}d (d m n - m - n +1)$ with $d=
|\nubold|$.
\end{lemma}

\begin{proof}
If ${\nubold }$ is not a disjoint union of mutually non-attacking
ribbon shapes, then any $m$-stretching of ${\nubold }$ contains two
successive boxes in reading order that are either on the same diagonal
(if some component is not a ribbon) or form an attacking pair.  The
Catalanimal $H^{m,n}_{\nubold ^m}$ then has at least one simple root
$\alpha_{i, i+1} \notin R_t$.  The factor $\prod_{\alpha \in R_+
\setminus R_t}(1-t\, \zz ^{\alpha})$ in $\phi(\zz )$ includes $(1-t\,
z_{i}/z_{i+1})$, so $\phi(1,t,\ldots, t^{l -1}) = 0$. By Lemma
\ref{lem:LLT 1-q}, $\omega \Gcal_{\nubold } [1-q] = 0$ as well.

When $\nubold = (\theta )$ is a single ribbon shape, $\phi (\zz )$ is given
by \eqref{e:phi-for-ribbon} and specializes to
\begin{equation}\label{e:phi-for-ribbon-spec}
\phi (1,t,\ldots,t^{l-1}) = t^{a - p} /(1-q)^{l-1},
\end{equation}
using Lemma~\ref{lem:a-two-ways} and the fact that the magic number
$p$ for this $\nubold $ is equal 
to $|I|$, where $I$ is as defined after \eqref{e:phi-for-ribbon}.
On the right hand side of \eqref{eq:LLT ps} we have $t^{a}
(-q\, t)^{-p} (\omega s_{\theta })[1-q]/(1-q)^{l}$.  This agrees with
\eqref{e:phi-for-ribbon-spec} because $(\omega s_{\theta })[1-q] =
(-q)^{p}(1-q)$, noting that $p$ is one less than the number of columns
of $\theta $.

Finally, when ${\nubold }$ is a disjoint union of non-attacking
ribbons, we use the fact that, as in the proof of
Theorem~\ref{thm:principal-spec}, \eqref{e:a(d)-additivity} and
\eqref{e:specialization-multiplicativity} imply that $t^{-a}\phi
(1,t,\ldots,t^{l-1})$ is multiplicative for elements $\phi \in S$ such
that $\psi (\phi )\in \Lambda (X^{m,n})$.  By Remark~\ref{rmk:product
of Catalanimals}, the Catalanimal $H^{m,n}_{\nubold ^{m}}$ in this
case is the product in $\Scal _{\hGamma }$ of the Catalanimals for the
individual ribbons $\theta $, so $\phi (\zz )$ is the product in $S$
of the corresponding functions $\phi _{\theta }$.  On the right hand
side of \eqref{eq:LLT ps}, the function $f$ is the product of the
functions $(- q t)^{-p_\theta } s_{\theta }(X)$, so \eqref{eq:LLT ps} reduces
to the single ribbon case.
\end{proof}

\begin{remark}\label{rem:LLT-1-0-ps}
When $(m,n) = (1,0)$, the $m$-stretching is trivial, so $p$ and
$n'(\gamma )$ in Lemma~\ref{lem:LLT ps} are the same as for $\nubold
$.  The offset $\obold $ is necessarily constant, so the number
$a_{ij}$ in \eqref{eq:aij def} is equal to $1$ for every attacking
pair, and $\Acal = A(\nubold )$ is the number of attacking pairs.  The
exponent $a$ is zero.  Hence, Lemma~\ref{lem:LLT ps} reduces to
Lemma~\ref{lem:LLT-1-0-ps} in this case.
\end{remark}

\begin{lemma}\label{l path underline basics}
Let $(m,n)\in \ZZ _{+}\times \ZZ $ be a coprime pair.

\noindent (i) For any interval $J = \{a, a+1, \ldots, b\} \subseteq
[dm]$, we have
\begin{equation}\label{e:b(dm,dn)-inequality}
\sum_{j \in J} \bbold (dm, dn)_{j} < \left| J \right| \frac{n}{m} + 1.
\end{equation}

\noindent (ii) For  $J = \{a, a+1, \ldots,d  m\}$, we have
\begin{equation}
\sum_{j \in J} \bbold (dm,dn)_{j} \le \left| J \right| \frac{n}{m}
\end{equation}
with equality if
and only if $a$ is one more than a multiple of $m$.
\end{lemma}

\begin{proof}
Computing directly from the definition \eqref{eq:def south steps} of
$\bbold (dm, dn)$, for any interval $J= \{a, a+1, \ldots, b\}$, there
holds $\left| J \right| \frac{n}{m} - \sum_{j \in J} \bbold (dm,
dn)_{j} = [(a-1) \frac{n}{m}]- [b \frac{n}{m}]$, using the notation
$[x] = \lceil x \rceil - x$.  Both parts of the lemma follow.
\end{proof}

\begin{proof}[Proof of Theorem \ref{t mn LLT cuddly}]
We verify that $H^{m,n}_{\nubold ^{m}}$ is $(m,n)$-cuddly and
use Corollary~\ref{cor:evaluating-cubs} to determine its cub.

Let $d = |\nubold |$ and $l = dm = |\nubold ^{m}|$.  For
$I\subseteq [l]$, we again use the abbreviation
\begin{equation}\label{e:r-I-mn}
r_I = |\lambda [I]_{I}| =  \sum_{i \in I} \lambda_i + |\barRq
^{I,I^c}| -|\barRt ^{I,I^c}|,
\end{equation}
now with the root sets and weight for $H^{m,n}_{\nubold ^{m}}$.
Replacing $\lambda $ with the weight $\hat{\lambda }$ for the
$(1,0)$-cuddly LLT Catalanimal $H_{\nubold ^{m}}$, we also set
\begin{equation}\label{e:r-hat-I-mn}
\hat{r}_I = |\hat{\lambda } [I]_{I}| = \sum_{i \in I} \hat{\lambda }_i
+ |\barRq ^{I,I^c}| -|\barRt ^{I,I^c}|,
\end{equation}
so that $r_I = \hat{r}_I + \sum_{i \in I}\bboldtild_i$, by
\eqref{e:lambda-vs-hat-lambda}.

\smallskip

\noindent {\em Checking the cuddliness conditions.}  The tameness
condition $[R_q,R_t] \subseteq R_{qt}$ holds by the same argument as
in Theorem \ref{t LLT cuddly}, applied to $\nubold ^m$ instead of
$\nubold $.
For the cuddliness bound,
it will be convenient to prove the following stronger claim:
\begin{equation}\label{e mn LLT cuddly induction}
\parbox{.85\linewidth}{$r_I \le \left| I \right| \frac{n}{m}$ for all
$I\subseteq [l]$, with equality if and only if $\nbm{I}$ is a lower
order ideal in $(\nubold ^m, \prec)$ and $\nbm{I}$ is a union of
stretched boxes.}
\end{equation}

If $\nbm{I}$ is a lower order ideal for $(\nubold ^m, \prec)$, then
$\hat{r}_I = 0$ by~\eqref{e LLT cuddly induction} in the proof of
$(1,0)$-cuddliness for $H_{\nubold ^m}$.  Hence,
\begin{equation}\label{e:r-I-for-order-ideal}
r_I= \sum_{i\in I}\bboldtild_i = \sum_{x \in \nubold }
\sum_{\substack{\nbm{i} \in\\
\stre(x) \cap \nbm{I}}} \bboldtild_{i} \,.
\end{equation}
Since $\nbm{I}$ is a lower order ideal, each sum $\sum_{ \nbm{i} \in
\stre(x) \cap \nbm{I}} \bboldtild_{i}$ equals $\sum_{i = j}^m \bbold
(m,n)_i$ for some $j$.  Thus, it follows from Lemma \ref{l path
underline basics} (ii) that $r_I \leq \left| I \right|\frac{n}{m}$,
with equality if and only if $\nbm{I}$ is a union of stretched boxes.

By the same argument as in the proof of \eqref{e LLT cuddly
induction}, if there is a pair $x-1,\, x \in [l]$ such that $x \in I$,
$x-1 \notin I$, and $\nbm{x-1}$ and $\nbm{x}$ belong to the same
diagonal, then $r_{J} > r_I$, where $J = I \cup \{x-1\} \setminus
\{x\}$.  Here we are using the fact that the weight $\lambda$ is
constant on diagonals of $\nubold ^m$.
It therefore suffices to prove \eqref{e mn LLT cuddly induction} for
$I$ satisfying
\begin{equation}\label{e:mn-reduction}
\text{for each diagonal $D \subseteq \nubold ^{m}$, $D\cap \nbm{I}$ is
a lower order ideal of $(D, \prec)$}.
\end{equation}
We prove this restricted statement by induction on $|I|$.
We can assume that $\nbm{I}$ is not a lower order ideal, since
we have already dealt with the case when it is.
Choose a column $C_0$ in one of the skew shapes of $\nubold ^m$ such
that 
$C_{0}$ contains at least one box of $\nbm{I}$, and no box of
$\nbm{I}$ is in the column immediately east of $C_{0}$.

If every box $b\in \nbm{I}\cap C_{0}$ has $\nubold^m \cap \{\south
(b),\west (b) \}\subseteq \nbm{I}$, let $C = \nbm{I}\cap C_{0}$; call
this Case 1.  Otherwise, let $\nbm{y}$ be the northernmost box of
$\nbm{I}\cap C_{0}$ such that $\nubold^m \cap \{\south (\nbm{y}),\west
(\nbm{y}) \}\not \subseteq \nbm{I}$, and let $C$ be the set of boxes
of $\nbm{I}\cap C_{0}$ north of and including $\nbm{y}$; call this
Case 2.  Let $K$ be the set of indices such that $\nbm{K} =
\nbm{I}\setminus C$.  If we remove the boxes of $C$ one at a time from
north to south, each of these boxes is $\prec $-maximal in the set
remaining just before we remove it.  Using \eqref{eq:reduction2} with
$\nubold ^{m}$ in place of $\nubold $ at each step, we obtain
$\hat{r}_I = \hat{r}_{K \cup \{y\}} < \hat{r}_{K}$ in Case 2, and
$\hat{r}_I = \hat{r}_{K}$ in Case 1.  Note that the boxes of $C$ are
contiguous in $C_{0}$ in both cases.  Hence, since the entries
$\bboldtild_{i}$ for $\nbm{i} \in C_0$ form the sequence
$\bbold(rm,rn)$ for some $r$, the entries $\bboldtild_{i}$ for
$\nbm{i} \in C$ form an interval in this sequence.  In Case 1, $C$
contains the southernmost box of $C_0$.  Then using $\hat{r}_I =
\hat{r}_{K}$ and Lemma \ref{l path underline basics} (ii) we obtain
\begin{equation*}
\left| I \right|\frac{n}{m} - r_I - (\left| K \right|\frac{n}{m} -
r_{K}) = \left| C \right|\frac{n}{m}- \hat{r}_I + \hat{r}_{K} -
\sum_{\nbm{i} \in C}\bboldtild_{i} = \left| C \right|\frac{n}{m} -
\sum_{\nbm{i} \in C}\bboldtild_{i} \ge 0.
\end{equation*}
Since $\nbm{I}$ is not a lower order ideal, there is some box $b\in
\nbm{I}$ such that $\nubold^m \cap \{\south (b),\west (b) \}\not
\subseteq \nbm{I}$.  Since we are in Case 1, such a box $b$ does not
belong to $C$, so $b\in \nbm{K}$.  This shows that $\nbm{K}$ is not a
lower order ideal.  Thus, by induction, $\left| K \right|\frac{n}{m} -
r_{K} > 0$, hence $\left| I \right|\frac{n}{m} - r_I >0$.

In Case 2, using $\hat{r}_I + 1 \le \hat{r}_{K}$ and 
Lemma \ref{l path underline basics} (i) we obtain
\begin{equation*}
\left| I \right|\frac{n}{m} - r_I - (\left| K \right|\frac{n}{m} -
r_{K}) = \left| C \right|\frac{n}{m}- \hat{r}_I + \hat{r}_{K} -
\sum_{\nbm{i} \in C}\bboldtild_{i} \ge \left| C \right|\frac{n}{m} -
\sum_{\nbm{i} \in C}\bboldtild_{i}+ 1 > 0.
\end{equation*}
By induction, $\left| K \right|\frac{n}{m} - r_{K} \ge 0$, hence
$\left| I \right|\frac{n}{m} - r_I > 0$.  This completes the proof of
\eqref{e mn LLT cuddly induction}.

\smallskip

\noindent {\em Determining the cub.} We now prove \eqref{eq:mn LLT
cuddly} using Corollary~\ref{cor:evaluating-cubs}, assuming by
induction that Theorem~\ref{t mn LLT cuddly} holds for smaller shapes
$\nubold $.

By \eqref{e mn LLT cuddly induction}, the subsets $I$ indexing the
summands in \eqref{e:cub-coprod-q-style} are characterized by the
property that $\nbm{I}$ is a lower order ideal and a union of
stretched boxes.  Given such an $I$, let $J\subseteq [d]$ be the set
of indices such that $\nbm{I}$ is the $m$-stretching $\nb{J}^{m}$ of
the lower order ideal $\nb{J}$ in $\nubold $, with the same offsets
$\obold $ as for $\nubold ^m = \nubold (m,\obold )$.  Then we also
have $\nbm{I^{c}} = \nb{J^{c}} ^{m}$.

Our first task is to show that $H_{I}'$ and $H_{I}''$ in
\eqref{e:cub-coprod-q-style} for $H = H^{m,n}_{\nubold ^{m}}$ are
given by $H_{I}' = H^{m,n}_{\nb{J}^{m}}$, $H_{I}'' =
H^{m,n}_{\nb{J^{c}}^{m}}$.
The proof is similar to the proof of the $(m,n) = (1,0)$ case in
Theorem~\ref{t LLT cuddly}.  In particular, the root sets clearly
agree.  We only discuss the adjustment needed to see that the
Catalanimals $(H^{m,n}_{\nubold ^m})_I'$ and $H^{m,n}_{\nb{J}^m}$ have
the same weight.  The adjustment for $(H^{m,n}_{\nubold ^m})_I''$ and
$H^{m,n}_{\nb{J^{c}}^m}$ is similar.

Since $\nbm{I}$ is a union of stretched boxes, $\lambda[I]_I-
\hat{\lambda}[I]_I = \bboldtild_I$ consists of copies of the
vector $\bbold (m,n)$ in the $m$ indices corresponding to the
individual boxes in each stretched box.  The proof of Theorem~\ref{t
LLT cuddly} shows that the weight of the $(1,0)$-cuddly Catalanimal
$H_{\nbm{I}} = H_{\nb{J}^{m}}$ is $\hat{\lambda}[I]_I$.  By
construction the weight of $H^{m,n}_{\nb{J}^m}$ is obtained by adding
$\bboldtild _{I}$ to this, so its weight is $\lambda [I]_{I}$, as
desired.

Now we turn to
the coefficients in \eqref{e:cub-coprod-q-style}.  Since 
$H^{m,n}_{\nubold ^m}$ and $H_{\nubold ^m}$
have the same root sets, the computation
is almost the same 
as in the proof of Theorem \ref{t LLT cuddly}, except that the
statistic in the exponent of $q$ is more complicated.  Let $p_I =
p(\nbm{I})$, $\gamma_I = \gamma(\nbm{I})$, $p_{I^c} = p(\nbm{I^c})$,
$\gamma_{I^c} = \gamma(\nbm{I^c})$ denote the magic number and
diagonal lengths of $\nbm{I}$ and $\nbm{I^c}$.  Just as in (\ref{eq:t
statistic}, \ref{eq:sign}), we have
\begin{align}
\label{eq:t statistic mn}
|\barRqt ^{I,I^c}| &= p+n'(\gamma) - p_I-n'(\gamma_I) -
p_{I^c}-n'(\gamma_{I^c}), \\ 
\label{eq:sign mn}
|\barRq ^{I,I^c}| & = n'(\gamma) - n'(\gamma_I)- n'(\gamma_{I^c}).
\end{align}

Next, by Remark \ref{rem:LLT Catalanimal def} (i), $|\barRt ^{I,I^c}|
= A(\nbm{I}, \nbm{I^c})$.  We need to relate this to attacking pairs
in $\nubold $.  Let $\Abold = \Abold (\nubold )$ be the set of
attacking pairs in $\nubold $, and let $\Abold ^{X,Y}$ denote the set
of attacking pairs from $\nb{X}$ to $\nb{Y}$ for any subsets $X,Y
\subseteq [d]$.  Using \eqref{eq:aij def}--\eqref{eq:aij 3}, we have
\begin{equation}\label{eq:q stat mn}
\begin{aligned}
|\barRt ^{I,I^c}| = A(\nbm{I}, \nbm{I^c}) &=
\sum_{(i,j) \in \Abold ^{J, J^c}}  a_{ij} + \sum_{(i,j) \in \Abold ^{J^c, J}}
  (a_{ij}-1)  \\
&= \sum_{(i,j) \in \Abold } a_{ij} - \sum_{(i,j) \in \Abold ^{J,J}}
a_{ij} - \sum_{(i,j) \in \Abold ^{J^c,J^c}} a_{ij} - |\Abold ^{J^c, J}|\\
& = \Asum - \Asum_J - \Asum_{J^c} - A(\nb{J^c}, \nb{J}),
 \end{aligned}
\end{equation}
where 
$\Asum_J = \sum_{(i,j) \in \Abold ^{J,J}} a_{ij}$, $\Asum_{J^c} =
\sum_{(i,j) \in \Abold ^{J^c,J^c}} a_{ij}$.

Combining \eqref{eq:t statistic mn}--\eqref{eq:q stat mn}, the term
indexed by $I$ in \eqref{e:cub-coprod-q-style} becomes
\begin{multline*}
(-1)^{p} (q\, t)^{-p-n'(\gamma)} q^{-\Acal } q^{A(\nb{J^c}, \nb{J})}\\ 
\times \Big( (-1)^{p_I} (q\, t)^{p_I+n'(\gamma_I)} q^{\Acal _{J}}
\cub(H^{m,n}_{\nb{J}^{m}})(X)\Big) \Big( (-1)^{p_{I^c}}
(q\, t)^{p_{I^c}+n'(\gamma_{I^c})} q^{\Acal _{J^c}}
\cub(H^{m,n}_{\nb{J^{c}}^{m}})(Y)\Big).
\end{multline*}
The desired formula \eqref{eq:mn LLT cuddly} now follows from
\eqref{e:LLT coprod}, Lemma \ref{lem:LLT ps}, and Corollary
\ref{cor:evaluating-cubs} just as in the proof of Theorem \ref{t LLT
cuddly}.
\end{proof}

\subsection{Formulas for \texorpdfstring{$\nabla$}{nabla} on LLT polynomials}
\label{ss:nabla-formulas}

Combining Theorem~\ref{t mn LLT cuddly} and
Proposition~\ref{prop:nabla-on-cuddlyish}, we obtain the following
formulas.  Recall that $\sigmabold $ denotes the Weyl symmetrization
operator in \eqref{e:Weyl-symmetrization}.

\begin{cor}
\label{cor:nabla on LLT} For any tuple of skew shapes $\nubold $ with
$|\nubold | = l$, we have the following raising operator style formula
for $\nabla$ applied to the associated LLT polynomial:
\begin{equation*}
(\omega \nabla \Gcal_{\nubold })(z_1,\dots, z_l) =
(-1)^p(q\, t)^{p+n'(\gamma)}q^{A} \, \sigmabold \Big(\frac{ \zz
^{\lambda+(1,\dots,1)} \, \prod_{\alpha\in R_{qt}}(1 - q\, t\, \zz^{\alpha})}
{\prod_{\alpha \in R_q}(1 - q\, \zz^\alpha)\prod_{\alpha \in
R_t}(1 - t\, \zz^\alpha)} \Big)_{\pol},
\end{equation*}
where $R_q$, $R_t$, $R_{qt}, \lambda$ are as in \eqref{ed
Rq}--\eqref{ed LLT weight}, $A = A(\nubold )$ is the number of
attacking pairs in $\nubold $, and $p = p(\nubold )$, $\gamma =
\gamma(\nubold )$ are the magic number and diagonal lengths of
$\nubold $.
\end{cor}

\begin{cor}
\label{cor:nabla on LLT2} More generally, $\nabla^m$ on any LLT
polynomial $ \Gcal_{\nubold }$ is given by
\begin{equation*}
(\omega \nabla^m \Gcal_{\nubold })(z_1,\dots, z_l) = (-1)^p(q\, t)^{p+m
\, n'(\gamma)}q^{\Asum} \, \sigmabold \Big(\frac{ \zz ^\lambda
\prod_{\alpha\in R_{qt}}(1 - q\, t\, \zz^\alpha)} {\prod_{\alpha \in
R_q}(1 - q\, \zz^\alpha)\prod_{\alpha \in R_t}(1 - t\, \zz^\alpha)}
\Big)_{\pol},
\end{equation*}
where $l= m|\nubold |$; $R_q$, $R_t$, $R_{qt}, \lambda$ are as in
Definition \ref{def:mn LLT Catalanimal} for $n=1$ and the
$m$-stretching $\nubold (m, \obold)$ of $\nubold $ with constant
offsets $\obold = (c,c,\ldots,c)$; 
$p =
p(\nubold )$, $\gamma = \gamma(\nubold )$, and $\Asum$ is $m$ times
the number of attacking pairs in $\nubold $ with the same content,
plus the number of attacking pairs with different contents.
\end{cor}

\begin{remark}
(i) If a Catalanimal $H(R_q, R_t, R_{qt}, \lambda)$ is $(1,0)$-cuddly, 
then $H(R_q, R_t, R_{qt}, \lambda+(1,\dots,1))$ is $(1,1)$-cuddly.
This explains why 
when $m=1$, the $\lambda $ in Corollary~\ref{cor:nabla on LLT2} becomes
$\lambda +(1,\ldots ,1)$ in Corollary~\ref{cor:nabla on LLT}. 

(ii) Corollary \ref{cor:nabla on LLT2} also holds for any $m$-stretching
$\nubold (m, \obold)$ of $\nubold $, but now with $\Asum = \sum_{(i,j)
\in \Abold (\nubold )} a_{ij}$, where $\Abold (\nubold )$ and $a_{ij}$
are as in \eqref{eq:aij def}.
\end{remark}

\begin{example}
(i) Continuing Example \ref{ex:some-LLT-cats} (i),
Corollary~\ref{cor:nabla on LLT} gives
\begin{align*}
& (\omega \nabla e_3) (z_1,z_2,z_3) = \big(z_1z_2z_3 \,
H_{((111))}\big)_{\pol} = \sigmabold \Big( \frac{z_1z_2z_3 \, (1 - q\,
t\, z_{1} /z_3)} {\prod_{1\le i < j \le 3}((1 - q\, z_{i}/ z_{j})(1 -
t\,z_{i}/ z_{j}))} \Big)_{\pol} \\[1mm] & = s_{111} + (q+t+q^2+q\,
t+t^2)s_{21}+ (q\, t + q^3+ q^2\, t + q\, t^2+t^3)s_{3}.
\end{align*}

(ii) When $\nubold $ is the single shape $((433))$,
$\Gcal_{\nubold}=s_{433}$. In this case, the LLT Catalanimal
$H_{((433))}^{1,1} = z_1 \cdots z_l \, H_{((433))}$ is shown in Figure
\ref{fig:Catalanimals intro} (i), and its associated statistics in
Example~\ref{ex:LLT Cat stats}. By Corollary~\ref{cor:nabla on LLT},
we have $(\omega \nabla s_{433})(z_1,\dots, z_l) = (q\,
t)^{9}\big(z_1\cdots z_l\, H_{((433))}\big)_{\pol}$.

(iii) For $\nubold = ((32)/(10),(33)/(11))$, the LLT Catalanimal
$H_{\nubold }^{1,1} = z_1 \cdots z_l \, H_{\nubold }$ is shown in
Figure \ref{fig:Catalanimals intro} (ii), and its associated
statistics in Example~\ref{ex:LLT Cat stats}.  By
Corollary~\ref{cor:nabla on LLT}, we have $(\omega \nabla
\Gcal_{\nubold})(z_1,\dots, z_l) = -(q\, t)^{4}q^7\big( z_1\cdots z_l
\, H_{\nubold }\big)_{\pol}$.
\end{example}

\bibliographystyle{hamsplain}
\bibliography{nablaonLLT}

\end{document}